\documentclass{article}
\usepackage{latexsym,amsthm,amssymb,amscd,amsmath}
\usepackage{geometry}
\newtheorem{thm}{Theorem}[section]

\newtheorem{lem}[thm]{Lemma}

\newcommand{\be}{\begin{equation}}
\newcommand{\ee}{\end{equation}}
\newcommand{\ben}{\begin{enumerate}}
\newcommand{\een}{\end{enumerate}}
\newcommand{\beq}{\begin{eqnarray}}
\newcommand{\eeq}{\end{eqnarray}}
\newcommand{\beqn}{\begin{eqnarray*}}
\newcommand{\eeqn}{\end{eqnarray*}}

\vspace{7cm} \textwidth 14cm \textheight 21.6cm

\title{Matsumoto Metrics of Reversible Curvature\footnote{Accepted in Acta  Mathematica  Academiae  Paedagogicae  Ny\'{i}regyh\'{a}ziensis, {\bf 31}(3) (2015) }}
\author{A. Tayebi and T. Tabatabaeifar}
\begin{document}

\maketitle
\begin{abstract}
In this paper, we study the reversibility of Riemann curvature and Ricci curvature for the  Matsumoto metric and prove three global results.  First, we prove that a Matsumoto metric is $R$-reversible if and only if  it is $R$-quadratic. Then we show that a Matsumoto metric is Ricci-reversible if and only if  it is Ricci-quadratic. Finally, we prove that every weakly Einstein Matsumoto metric is Ricci-reversible.
\end{abstract}
\section{Introduction}
The class of $(\alpha, \beta)$-metrics was introduced by Matsumoto as extension of Randers and Kropina metrics \cite{Mat5}. An $(\alpha, \beta)$-metric is a Finsler metric on $M$ defined by $F:=\alpha\phi(s)$, where $s=\beta/\alpha$,  $\phi=\phi(s)$ is a $C^\infty$ function on the $(-b_0, b_0)$ with certain regularity, $\alpha$ is a Riemannian metric and $\beta$ is a 1-form on $M$ \cite{LP}\cite{TS}\cite{TTP}.

Recently, Crampin proved that a Randers metric $F=\alpha+\beta$ has reversible geodesics if and only if $\beta$ is parallel with respect to $\alpha$ \cite{C}. Then Masca-Sabau-Shimada investigate $(\alpha, \beta)$-metrics with reversible geodesics and projectively reversible geodesics \cite{MSS1}\cite{MSS2}. In general, the Finsler metrics  might not be reversible. In spite of the non-reversibility of Finsler metrics, the geodesics and curvatures might be reversible.

In \cite{SYa}, Shen-Yang introduced the notions of R-reversibility and Ricci-reversibility. They proved that Randers metrics are R-reversible or Ricci-reversible if and only if they are R-quadratic or Ricci-quadratic, respectively. In this paper, we are going to study the reversibility of Riemann curvature and Ricci curvature for a Matsumoto metric $F=\alpha^2/(\alpha-\beta)$ which  is called by Matsumoto's slope-of-a-mountain metric, also.  This metric was introduced by Matsumoto as a realization of Finsler's idea  ``a slope measure of a mountain with respect to a time measure" \cite{SS}\cite{TPS}. He gave an exact formulation of a Finsler surface to measure the time on the slope of a hill and introduced the Matsumoto metrics in \cite{Mat}\cite{TTP}.

The Riemann curvature ${\bf R}_y : T_xM \to T_xM$ is a family of linear maps on tangent spaces. A Finsler metric $F$  on a manifold $M$ is said to be R-quadratic if its Riemann curvature ${\bf R}_y$ is quadratic in $y\in T_xM$ \cite{LS}\cite{SYa}. The notion of R-quadratic metric was introduced by Shen \cite{Shrq}.  $F$ is called R-reversible if ${\bf R}_y={\bf R}_{-y}$ (see \cite{SYa}). A Finsler metric $F$ is called Ricci-quadratic if its Ricci curvature ${\bf Ric}_y$,  is quadratic in $y\in T_xM$. $F$ is called Ricci-reversible if ${\bf Ric}_y={\bf Ric}_{-y}$ (see \cite{SYa}).  In this paper, we have
the two following theorems.
\begin{thm}\label{theorem2}
A Matsumoto metric is $Ricci$-reversible if and only if it is $Ricci$-quadratic.
\end{thm}
\begin{thm}\label{theorem1}
A Matsumoto metric  is $R$-reversible if and only if it is $R$-quadratic.
\end{thm}
A Finsler metric on an $n$-dimensional manifold $M$ is said to be of weakly Einsteinian if the Ricci curvature be in the  form ${\bf Ric}=(n-1)[3\theta F+\sigma F^2]$, where $\theta=\theta_i(x)y^i$ is a 1-form and $\sigma=\sigma(x)$ is a scalar function on $M$.
\begin{thm}\label{mainthm3}
Every weakly Einstein Matsumoto metric is Ricci reversible.
\end{thm}
\section{Preliminaries}

An important class of Finsler metrics is so called $(\alpha,\beta)$-metrics, which are expressed in the form of $F=\alpha \phi(s)$, $s =\beta/\alpha$ where $\alpha=\sqrt{a_{ij}(x)y^iy^j}$ is a Riemannian metric and $\beta =b_i(x)y^i$ is a 1-form. Put
\be
r_{ij}:=\frac{1}{2}(b_{i|j}+b_{j|i}),\qquad s_{ij}:=\frac{1}{2}(b_{i|j}-b_{j|i}),
\ee
where $'|'$ denotes the covariant derivative with respect to the Levi-Civita connection of $\alpha$. Let
\begin{align*}
r^i_j&:=a^{im}r_{mj},\ \ s^i_j:=a^{im}s_{mj}, \ \ r_j:=b^mr_{mj},\ \ s_j:=b^ms_{mj}, \ \  r:=b^ir_i,
\end{align*}
where $a^{ij}:=(a_{ij})^{-1}$, $b^i:=a^{ij}b_j$. We
the subscript $0$ means the contraction by $y^i$.
Let $G^i$ and $\bar{G}^i$ be the geodesic coefficients of $F$ and $\alpha$ respectively. Then we
have the following
\begin{lem}\rm{(\cite{shen})}
Let $G^i=G^i(x,y)$ and $\bar{G}^i_{\alpha}=\bar{G}^i_{\alpha}(x,y)$ denote the
coefficients  of $F$ and  $\alpha$ respectively in the same coordinate system. Then, we have
\be
G^i = \bar{G}^i+\alpha Qs^i_{~0}+\Psi(r_{00}-2\alpha Qs_0)b^i+\frac{1}{\alpha}\Theta(r_{00}-2\alpha Qs_0)y^i,\label{a3}
\ee
where
\begin{align*}
Q&:=\frac{\phi'}{\phi -s\phi'},\ \ \ \Psi :=\frac{\phi^{''}}{2[\phi -s\phi' +(b^2-s^2)\phi^{''}]}\ \ \
\Theta :=\frac{\phi\phi' -s(\phi\phi^{''}+\phi'^2}{2\phi[\phi -s\phi' +(b^2-s^2)\phi^{''}]}.
\end{align*}
\end{lem}
For a Matsumoto metric $F=\frac{\alpha^2}{\alpha-\beta}$ and by a quite long computational procedure using Maple program, we obtain the following (see Proposition 3.1 in \cite{CST}).
\begin{lem}
Let $F=\alpha^2/(\alpha -\beta)$ be a Matsumoto metric on a manifold $M$. Then the Riemannian curvature of F is given by
\begin{equation}\label{a4}
R^i_{~j}=\bigg(\frac{1}{4\alpha^4(\alpha -3\beta +2b^2\alpha)^4(\alpha -2\beta)^3}\bigg)\sum_{k=0}^{13}t_k\alpha^k,
\end{equation}
where $\bar{R}^i_{~j}$ is the Riemannian curvature of $\alpha$ and $t_k, k=0, 1, ..., 13$ are as follows
\begin{equation*}
\left\{
\begin{array}{rl}
                       t_0&:=2880y^iy_jr_{00}^2\beta^7,\\
                       &\vdots \\
                       t_{13}&:=4(2b^2+1)^3(-2s^i_ks^k_jb^2-s^i_ks^k_j-2s^is_j+2b^is_ms^m_j)
   \end{array} \right.
\end{equation*}
All the coefficients of $t_i$ are tedious, listed in Appendix 1.
\end{lem}
By \cite{CST},  we have the following.
\begin{lem}
 The Ricci curvature of  Matsumoto metric $F=\alpha^2/(\alpha -\beta)$ is given by
\begin{equation}\label{a8}
\textbf{Ric}=\bigg(\frac{1}{4\alpha^2(\alpha -3\beta +2b^2\alpha)^4(\alpha -2\beta)^3}\bigg)\sum_{k=0}^{11}d_k\alpha^k
\end{equation}
where ${\bf \overline{Ric}}:=\bar{R}^m_{~m}$ and  and $d_k, k=0, 1, ..., 11$ are as follows
\begin{equation*}
\left\{
\begin{array}{rl}
                       d_0 &:=-288r_{00}^2(8n-11)\beta^7.\\
                       &\vdots \\
                       d_{11}&:=-4(1+2b^2)^3(2s_m^is_i^mb^2+s_m^is_i^m+4s_ms^m)
\end{array} \right.
\end{equation*}
All the coefficients of $d_i$ are tedious, listed in Appendix 2.
\end{lem}

\section{Proof of Theorem \ref{theorem2}}
{\bf Proof of Theorem \ref{theorem2}:} Let the Ricci curvature of $F$ be reversible, i.e., $\textbf{Ric}(y)=\textbf{Ric}(-y)$. Then by contracting both sides of  (\ref{a8}) with $4\alpha^2(\alpha -3\beta +2b^2\alpha)^4(\alpha -2\beta)^3$ and by a quite long computational procedure using Maple program, we obtain
\begin{align}\label{a1}
\sum_{i=0}^5d'_{2i}\alpha^{2i} =0,
\end{align}
where
\begin{equation*}
d'_0:=288(8n-11)\beta^7r_{00}^2.
\end{equation*}
All the coefficients of $d'_{2i}$ are listed in Appendix 3.
From (\ref{a1}), we know that $\alpha^2$ divides $d'_0$. Since $\alpha^2$ is an irreducible polynomial in $y$, it must be the case that $\alpha^2$ divides $r_{00}$. Thus we have
\[
r_{00}=c\alpha^2
\]
for some function $c=c(x)$, i.e., $\beta$ is a conformal form with respect to $\alpha$. So it is easy to get the following
\begin{align}
r_{00}&=c\alpha^2, \ \ r_{0j}=cy_j, \ \  r_i=cb_i, \ \ r=cb^2, \ \  r^i_j=c\delta^i_j,\label{a15}\\
\nonumber r_{0k}s^k_{~0}&=0, \ \ r_{0k}s^k=cs_0, \ \  s^k_{~0}r_k=cs_0, \ \ r_{00|0}=c_0\alpha^2.
\end{align}
Plugging above equations into (\ref{a1}) imply that
\begin{align}
d^{''}_0+d^{''}_2\alpha^2 +\cdots +d^{''}_8\alpha^8=0, \label{a2}
\end{align}
where
\[
d^{''}_0:=-2592(\textbf{Ric}-{\bf \overline{Ric}})\beta^7
\]
and other coefficients of $d^{''}_i$ are listed in Appendix 4. By (\ref{a2}), we get
\[
d^{''}_0=k\alpha^2
\]
where $k=k(x)$ is a scalar function on $M$.  So we have
\[
\textbf{Ric}={\bf \overline{Ric}}-k\alpha^2.
\]
This shows that $F$ is Ricci-quadratic.

Conversely, assume that the Ricci curvature of $F$ is quadratic. Then by (\ref{a8}) we have
\begin{align}
d'_0&+d'_2\alpha^2 +\cdots +d'_{10}\alpha^{10}=0, \label{a13}\\
\nonumber d'_1&+d'_3\alpha^2 +\cdots +d'_{11}\alpha^{10}=0,
\end{align}
where $d'_{2i}$ are the same as Ricci-reversibility (since we only use coefficients of $d'_{2i}$, then we relinquish the coefficients of $d'_{2i+1}$). Since (\ref{a13}) is the same as (\ref{a1}), by the same method we get
\[
\textbf{Ric}={\bf \overline{Ric}}-k\alpha^2
\]
which shows that $F$ is Ricci-reversible.
\qed

\section{Proof of Theorem \ref{theorem1}}
{\bf Proof of Theorem \ref{theorem1}:} Let $F$ be $R$-reversible, $\mathbf{R}(y)=\mathbf{R}(-y)$. Then by contracting both side (\ref{a4}) in $4\alpha^4(\alpha -3\beta +2b^2\alpha)^4(\alpha -2\beta)^3$ and using (\ref{a15}) and by a quite long computational procedure using Maple program, we obtain
\begin{equation}\label{a5}
\sum_{i=0}^6t'_{2i}\alpha^{2i}=0,
\end{equation}
where
\begin{eqnarray*}
t'_0:=\!\!\!\!\!\!&-&\!\!\!\!\!\!\ 288[10 y^iy_jc^2\beta^2 +(12y^iy_jcs_0\beta +17y^iy_jc_0\beta +16y^iy_jc_0b^2\beta -16y^iy_jc^2\beta^2 +4y^iy_jc^2\beta^2) \\
\!\!\!\!\!\!&+&\!\!\!\!\!\!\ (32y^iy_js_0^2b -20y^iy_jcs_0\beta -8y^iy_js_{0|0}b^2 -6y^iy_js_0^2 -10y^iy_js_{0|0} -4y^icy_js_0\beta -8y^icy_jc\beta^2 \\
\!\!\!\!\!\!&-&\!\!\!\!\!\!\ 24y^i c_0y_j b^2\beta -24s_{0|0}^iy_jb^2\beta  +9(R^i_j-\overline{R}^i_j)\beta^2 -2979s^i_0s_{0j}\beta^2 -33y^ic_0y_j\beta -60s_{0|0}^iy_j\beta \\
\!\!\!\!\!\!&+&\!\!\!\!\!\!\ 32y^iy_jcs_0b\beta   +32y^icy_js_0b\beta+ 32y^ic^2y_jb\beta^2 )]\beta^5
\end{eqnarray*}
and other coefficients of $t'_i$ are listed in Appendix 5. By (\ref{a5}), it follows that $\alpha^2$ divides $t'_0$, which is impossible. Therefore
\[
t'_0=0
\]
and we have
\begin{eqnarray}
R^i_{~j}:=\overline{R}^i_{~j} \!\!\!\!\!\!&-&\!\!\!\!\!\!\ \frac{1}{9}\big[10 y^iy_jc^2 -16y^iy_jc^2 +4y^iy_jc^2 -8y^iy_jc^2 -2979s^i_0s_{0j} +32y^ic^2y_jb\label{a16}\\
\nonumber  \!\!\!\!\!\!&+&\!\!\!\!\!\!\  \frac{1}{\beta}(12y^iy_jcs_0 +17y^iy_jc_0 +16y^iy_jc_0b^2 -20y^iy_jcs_0 -4y^icy_js_0  -24y^i c_0y_j b^2 \\
\nonumber  \!\!\!\!\!\!&-&\!\!\!\!\!\!\  24s_{0|0}^iy_jb^2+32y^iy_jcs_0b -33y^ic_0y_j -60s_{0|0}^iy_j +32y^icy_js_0b) \\
\nonumber  \!\!\!\!\!\!&+&\!\!\!\!\!\!\ \frac{1}{\beta^2} (32y^iy_js_0^2b -8y^iy_js_{0|0}b^2 - 6y^iy_js_0^2 -10y^iy_js_{0|0})\big].
\end{eqnarray}
This means that $F$ is $R$-quadratic.

Conversely, let $F$ be R-quadratic. Then by the same method we have
\begin{align}
t'_0&+\cdots +t'_{12}\alpha^{10}=0, \label{a14}\\
t'_1&+\cdots +t'_{13}\alpha^{12}=0,\label{a17}
\end{align}
where coefficients of $t'_{i}$ are listed in Appendix 5. (since we only use coefficients of $t'_{2i}$, then we relinquish the coefficients of $t'_{2i+1}$). By (\ref{a14}), we get  (\ref{a16}) which proves that $F$ is R-reversible.\qed

\section{Proof of Theorem \ref{mainthm3}}
{\bf Proof of Theorem \ref{mainthm3}:} Let  $F$ be weakly Einstein
\be
{\bf Ric}=(n-1)\big[3\theta F+\sigma F^2\big].\label{E}
\ee
By substituting  (\ref{a8}) in (\ref{E}) and with a long computational procedure, using Maple program, we get
\begin{equation}\label{s3}
\bigg[\frac{1}{4\alpha^2(\alpha -3\beta +2b^2\alpha)^4(\alpha -2\beta)^3(\alpha-\beta)^2}\bigg]\sum_{k=0}^{13}A_k\alpha^k =0,
\end{equation}
where
\begin{eqnarray*}
A_0:=-288\beta^9r_{00}^2(8n-11),
\end{eqnarray*}
 and other coefficients of $A_{i}$ are listed in Appendix 6. (\ref{s3}) is equivalent to the following equations
\begin{align}
A_0+A_2\alpha^2+\cdots +A_{12}\alpha^{12}&=0, \label{s4}\\
\nonumber A_1+A_3\alpha^2+\cdots +A_{13}\alpha^{12}&=0.
\end{align}
By (\ref{s4}),  it follows that $\alpha^2$ divides $A_0$. Since $\alpha^2$ is an irreducible polynomial in $y$, it must be the case that $\alpha^2$ divides $r_{00}$. Thus
\[
r_{00}=c\alpha^2
\]
for some function $c=c(x)$, i.e., $\beta$ is a conformal form with respect to $\alpha$. So it is easy to get
\begin{align*}
&&r_{00}=c\alpha^2, \ \ r_{0j}=cy_j, \ \  r_i=cb_i, \ \ r=cb^2, \ \ r^i_{~j}=c\delta^i_j,\\
&&r_{0k}s^k_{~0}=0, \ \ r_{0k}s^k=cs_0,  \ \ s^k_{~0}r_k=cs_0, \ \ r_{00|0}=c_0\alpha^2.
\end{align*}
Plugging the above equations into (\ref{s4}) yields
\begin{equation}\label{s5}
\bigg[\frac{1}{4\alpha^2(\alpha -3\beta +2b^2\alpha)^4(\alpha -2\beta)^3(\alpha-\beta)^2}\bigg]\sum_{k=0}^{5}A'_{2k}\alpha^{2k} =0,
\end{equation}
where
\[
A'_0:=2592\beta^9\big[{\bf Ric}-{\bf \overline{Ric}}\big]
\]
and other coefficients of $A'_{2i}$ are listed in Appendix 7.

By (\ref{s5}), it follows that $\alpha^2$ must divide $A'_0$. So we have
\[
\textbf{Ric}={\bf \overline{Ric}}-k\alpha^2,
\]
where $k=k(x)$ is a scalar function on $M$. Therefore,   $F$ is Ricci-reversible.
\qed

\section{Appendix 1: Coefficients in (\ref{a4})}

\begin{eqnarray}
\nonumber t_0 \!\!\!\!&:=&\!\!\!\! 2880y^iy_jr_{00}^2\beta^7,
\\
\nonumber t_1\!\!\!\!&:=&\!\!\!\! -192y^iy_j\beta^6(9r_{00|0}\beta +28r_{00}^2b^2 +35r_{00}^2)
\\
\nonumber t_2 \!\!\!\!&:=&\!\!\!\! 48\beta^5(24y^ir_{j0}r_{00}\beta^2 -48\delta_j^ir_{00}^2\beta^2 +18b^iy_jr_{00}^2\beta +72y^iy_jr_{00}s_0\beta +102y^iy_jr_{00|0}\beta \\
\nonumber  \!\!\!\!\!\!&-&\!\!\!\!\!\!\ 72y^iy_jr_{00}r_0\beta -36y^ib_jr_{00}^2\beta +96y^iy_jr_{00|0}b^2\beta +131y^iy_jr_{00}^2 +64y^iy_jr_{00}^2b^4 +222y^iy_jr_{00}^2b^2)
\end{eqnarray}
\begin{eqnarray}
\nonumber t_3 \!\!\!\!&:=&\!\!\!\!\ -48\beta^4 \big[154y^iy_jr_{00}s_0\beta -170y^iy_jr_{00}r_0\beta +230y^iy_jr_{00|0}b^2\beta +80y^iy_jr_{00|0}b^4\beta -48y^ib_jr_{00}^2b^2\beta \\
\nonumber  \!\!\!\!\!\!&+&\!\!\!\!\!\!\ 122y^iy_jr_{00|0}\beta +144y^iy_jr_{0m}s^m_0\beta^2 -36y^iy_js_{0|0}\beta^2 +21b^iy_jr_{00}^2\beta +18b^iy_jr_{00|0}\beta^2 \\
\nonumber  \!\!\!\!\!\!&-&\!\!\!\!\!\!\ 84y^ib_jr_{00}^2\beta -80\delta_j^ir_{00}^2b^2\beta^2 +68y^ir_{j0}r_{00}\beta^2 -54y^is_jr_{00}\beta^2 -112y^iy_jr_{00}r_0b^2\beta \\
\nonumber  \!\!\!\!\!\!&-&\!\!\!\!\!\!\ 128y^iy_jr_{00}s_0b^3\beta +128y^iy_jr_{00}s_0b^2\beta -40y^iy_jr_{00}s_0b\beta -40y^iy_jr_{00}r_0b\beta -128y^iy_jr_{00}r_0b^3\beta \\
\nonumber  \!\!\!\!\!\!&+&\!\!\!\!\!\!\ 72y^ir_{00|j}\beta^3 -72y^ir_{j0|0}\beta^3 -54s_{0|0}^iy_j\beta^3 -124\delta_j^ir_{00}^2\beta^2 -36\delta_j^ir_{00|0}\beta^3 +24b^iy_jr_{00}^2b^2\beta \\
\nonumber  \!\!\!\!\!\!&+&\!\!\!\!\!\!\  16y^ir_{j0}r_{00}b^2\beta^2 +18r^i_0y_jr_{00}\beta^2 +63y^iy_jr_{00}^2 +120y^iy_jr_{00}^2b^4 +183y^iy_jr_{00}^2b^2\big]
\\
\nonumber t_4 \!\!\!\!&:=&\!\!\!\!\ 4\beta^3(1152\delta_j^ir_{00}s_0b\beta^3+1152\delta_j^ir_{00}r_0b\beta^3+1270y^iy_jr_{00}s_0\beta -1964y^iy_jr_{00}r_0\beta +2700y^iy_jr_{00|0}b^2\beta \\
\nonumber  \!\!\!\!\!\!&+&\!\!\!\!\!\!\ 1920y^iy_jr_{00|0}b^4\beta +2304y^iy_js_0^2b\beta^2 -1440y^iy_jr_0s_0\beta^2 +4032y^iy_jr_{0m}s^m_0b^2\beta^2 -576y^iy_js_{0|0}b^2\beta^2 \\
\nonumber  \!\!\!\!\!\!&-&\!\!\!\!\!\!\ 384y^ib_jr_{00}^2b^4\beta -648y^ib_jr_{00}^2b^2\beta -720y^ib_jr_{00}s_0\beta^2 +288y^ib_jr_{00}r_0\beta^2-288y^ib_jr_{00|0}b^2\beta^2 \\
\nonumber  \!\!\!\!\!\!&-&\!\!\!\!\!\!\ 864y^is_jr_{00}b^2\beta^2+956y^iy_jr_{00|0}\beta -432y^iy_js_0^2\beta^2 +3744y^iy_jr_{0m}s^m_0\beta^2 -720y^iy_js_{0|0}\beta^2 \\
\nonumber  \!\!\!\!\!\!&+&\!\!\!\!\!\!\ 54b^iy_jr_{00}^2\beta +468b^iy_jr_{00|0}\beta^2 -72b^ib_jr_{00}^2\beta^2 -894y^ib_jr_{00}^2\beta -36y^ib_jr_{00|0}\beta^2 +576y^ir_jr_{00}\beta^3 \\
\nonumber  \!\!\!\!\!\!&+&\!\!\!\!\!\!\ 720y^is_jr_{00}\beta^3 -384\delta_j^ir_{00}^2b^4\beta^2 -2352\delta_j^ir_{00}^2b^2\beta^2 -576\delta_j^ir_{00}s_0\beta^3 -288y^ir_{j0}s_0\beta^3 +960y^ir_{j0}r_{00}\beta^2\\
\nonumber  \!\!\!\!\!\!&-&\!\!\!\!\!\!\ 576y^ir_{j0}r_0\beta^3+1728y^ir_{00|j}b^2\beta^3-1728y^ir_{j0|0}b^2\beta^3-1728s_{0|0}^iy_jb^2\beta^3 -144b^ir_{j0}r_{00}\beta^3
\\
\nonumber  \!\!\!\!\!\!&-&\!\!\!\!\!\!\  1404y^is_jr_{00}\beta^2 -648\overline{R}^i_j\beta^4 -214488s^i_0s_{0j}\beta^4 -2672y^iy_jr_{00}r_0b^2\beta-3328y^iy_jr_{00}s_0b^3\beta \\
\nonumber  \!\!\!\!\!\!&+&\!\!\!\!\!\!\ 2368y^iy_jr_{00}s_0b^2\beta-976y^iy_jr_{00}s_0b\beta -976y^iy_jr_{00}r_0b\beta -3328y^iy_jr_{00}r_0b^3\beta+2304y^iy_jr_0s_0b\beta^2\\
\nonumber  \!\!\!\!\!\!&-&\!\!\!\!\!\!\ 1024y^iy_jr_{00}s_0b^5\beta +1024y^iy_jr_{00}s_0b^4\beta+576\delta_j^ir_{00}r_0\beta^3 -864\delta_j^ir_{00|0}b^2\beta^3 +1128y^iy_jr_{00}^2b^4 \\
\nonumber  \!\!\!\!\!\!&+&\!\!\!\!\!\!\ 2376y^ir_{00|j}\beta^3-2376y^ir_{j0|0}\beta^3-2160s_{0|0}^iy_j\beta^3 -1638\delta_j^ir_{00}^2\beta^2-1188\delta_j^ir_{00|0}\beta^3 +1050y^iy_jr_{00}^2b^2
\\
\nonumber  \!\!\!\!\!\!&-&\!\!\!\!\!\!\ 144b^iy_jr_{00}r_0\beta^2+288b^iy_jr_{00|0}b^2\beta^2+256y^iy_jr_{00|0}b^6\beta -2304y^ir_jr_{00}b\beta^3 +504b^iy_jr_{00}^2b^2\beta\\
\nonumber  \!\!\!\!\!\!&+&\!\!\!\!\!\!\ 144b^iy_jr_{00}s_0\beta^2 -1024y^iy_jr_{00}r_0b^5\beta-512y^iy_jr_{00}r_0b^4\beta -1152b^iy_jr_{00}s_0b\beta^2 -1152b^iy_jr_{00}r_0b\beta^2 \\
\nonumber  \!\!\!\!\!\!&-&\!\!\!\!\!\!\  1152y^ib_jr_{00}s_0b\beta^2-1152y^ib_jr_{00}r_0b\beta^2+336y^ir_{j0}r_{00}b^2\beta^2 +2304y^ir_{j0}s_0b\beta^3+2304y^ir_{j0}r_0b\beta^3\\
\nonumber  \!\!\!\!\!\!&-&\!\!\!\!\!\!\ 2304y^is_jr_{00}b\beta^3+288r^i_0y_jr_{00}b^2\beta^2 +468r^i_0y_jr_{00}\beta^2 +486\beta^2s^i_0y_jr_{00}+207y^iy_jr_{00}^2)
\\
\nonumber t_5 \!\!\!\!&:=&\!\!\!\!\ -4\beta^2\big[972\beta^2s^i_0y_jr_{00}b^2 -1152\delta_j^ir_{00}s_0b^2\beta^3 +2784\delta_j^ir_{00}s_0b\beta^3 +1536\delta_j^ir_{00}s_0b^3\beta^3 +768\delta_j^ir_{00}r_0b^2\beta^3 \\
\nonumber  \!\!\!\!\!\!&+&\!\!\!\!\!\!\ 2784\delta_j^ir_{00}r_0b\beta^3 +293y^iy_jr_{00}s_0\beta -982y^iy_jr_{00}r_0\beta +1374y^iy_jr_{00|0}b^2\beta +1488y^iy_jr_{00|0}b^4\beta \\
\nonumber  \!\!\!\!\!\!&-&\!\!\!\!\!\!\ 576y^iy_js_0^2b^2\beta^2 +3360y^iy_js_0^2b\beta^2 +1536y^iy_js_0^2b^3\beta^2 -2568y^iy_jr_0s_0\beta^2 +7248y^iy_jr_{0m}s^m_0b^2\beta^2 \\
\nonumber  \!\!\!\!\!\!&+&\!\!\!\!\!\!\ 3072y^iy_jr_{0m}s^m_0b^4\beta^2 -696y^iy_js_{0|0}b^2\beta^2-192y^iy_js_{0|0}b^4\beta^2 -3072y^ir_jr_{00}b^3\beta^3 -192b^ib_jr_{00}^2b^2\beta^2 \\
\nonumber  \!\!\!\!\!\!&-&\!\!\!\!\!\!\ 480y^ib_jr_{00}^2b^4\beta +84y^ib_jr_{00}^2b^2\beta -1848y^ib_jr_{00}s_0\beta^2 +648y^ib_jr_{00}r_0\beta^2 -672y^ib_jr_{00|0}b^2\beta^2 \\
\nonumber  \!\!\!\!\!\!&-&\!\!\!\!\!\!\ 384y^ib_jr_{00|0}b^4\beta^2 -1260y^is_jr_{00}b^2\beta^2 -288y^is_jr_{00}\beta^2b^4 +367y^iy_jr_{00|0}\beta -504y^iy_js_0^2\beta^2
\\
\nonumber  \!\!\!\!\!\!&+&\!\!\!\!\!\!\ 3342y^iy_jr_{0m}s^m_0\beta^2 -462y^iy_js_{0|0}\beta^2 +648b^iy_js_{0|0}\beta^3 -3b^iy_jr_{00}^2\beta +402b^iy_jr_{00|0}\beta^2
\\
\nonumber  \!\!\!\!\!\!&+&\!\!\!\!\!\!\ 228b^ib_jr_{00}^2\beta^2 -216b^ib_jr_{00|0}\beta^3 +432y^ib_js_{0|0}\beta^3 -351y^ib_jr_{00}^2\beta -78y^ib_jr_{00|0}\beta^2 +1392y^ir_jr_{00}\beta^3\\
\nonumber  \!\!\!\!\!\!&-&\!\!\!\!\!\!\ 432y^ib_jr_{0m}s^m_0\beta^3 +1728y^is_jr_{00}\beta^3 -864\delta_j^ir_{00}^2b^4\beta^2 -2424\delta_j^ir_{00}^2b^2\beta^2 -1008\delta_j^ir_{00}s_0\beta^3
\\
\nonumber  \!\!\!\!\!\!&-&\!\!\!\!\!\!\ 720y^ir_{j0}s_0\beta^3 +600y^ir_{j0}r_{00}\beta^2 -1392y^ir_{j0}r_0\beta^3 +4176y^ir_{00|j}b^2\beta^3 +1152y^ir_{00|j}\beta^3b^4
\\
\nonumber  \!\!\!\!\!\!&-&\!\!\!\!\!\!\ 4176y^ir_{j0|0}b^2\beta^3 -1152y^ir_{j0|0}\beta^3b^4 -5184s_{0|0}^iy_jb^2\beta^3 -1728s_{0|0}^iy_j\beta^3b^4 -480b^ir_{j0}r_{00}\beta^3
\\
\nonumber  \!\!\!\!\!\!&-&\!\!\!\!\!\!\ 648y^is_js_0\beta^3 -1206y^is_jr_{00}\beta^2 -324g^iy_jr_{00}\beta^3 +1536\delta_j^ir_{00}r_0b^3\beta^3 -1836\overline{R}^i_j\beta^4 +39y^iy_jr_{00}^2 \\
\nonumber  \!\!\!\!\!\!&-&\!\!\!\!\!\!\ 607716s^i_0s_{0j}\beta^4 -2056y^iy_jr_{00}r_0b^2\beta -2816y^iy_jr_{00}s_0b^3\beta +1268y^iy_jr_{00}s_0b^2\beta -776y^iy_jr_{00}s_0b\beta \\
\nonumber  \!\!\!\!\!\!&-&\!\!\!\!\!\!\ 776y^iy_jr_{00}r_0b\beta -2816y^iy_jr_{00}r_0b^3\beta -2112y^iy_jr_0s_0b^2\beta^2 +3360y^iy_jr_0s_0b\beta^2 +1536y^iy_jr_0s_0b^3\beta^2 \\
\nonumber  \!\!\!\!\!\!&-&\!\!\!\!\!\!\ 1664y^iy_jr_{00}s_0b^5\beta +1472y^iy_jr_{00}s_0b^4\beta +1392\delta_j^ir_{00}r_0\beta^3 -2088\delta_j^ir_{00|0}b^2\beta^3 -576\delta_j^ir_{00|0}b^4\beta^3 \\
\nonumber  \!\!\!\!\!\!&-&\!\!\!\!\!\!\ 571968s^i_0s_{0j}b^2\beta^4 +432R4r_{j0}\beta^4 +864y^is_{0|j}\beta^4 +432b^ir_{j0|0}\beta^4 -432b^ir_{00|j}\beta^4 -432r_j^ir_{00}\beta^4
\end{eqnarray}
\begin{eqnarray}
\nonumber  \!\!\!\!\!\!&+&\!\!\!\!\!\!\ -864\delta_j^ir_{0m}s^m_0\beta^4 -432\delta_j^is_{0|0}\beta^4 +2772y^ir_{00|j}\beta^3 -2772y^ir_{j0|0}\beta^3 +432y^ir_{k0}s^k_j\beta^4 \\
\nonumber  \!\!\!\!\!\!&-&\!\!\!\!\!\!\ 864y^ir_{jk}s^k_0\beta^4 -432y^is_{j|0}\beta^4 -2808s_{0|0}^iy_j\beta^3 -1005\delta_j^ir_{00}^2\beta^2 -1386\delta_j^ir_{00|0}\beta^3 +492y^iy_jr_{00}^2b^4 \\
\nonumber  \!\!\!\!\!\!&+&\!\!\!\!\!\!\ 324b^is_jr_{00}\beta^3 -264b^iy_jr_{00}r_0\beta^2 +528b^iy_jr_{00|0}b^2\beta^2 +96b^iy_jr_{00|0}b^4\beta^2 +416y^iy_jr_{00|0}b^6\beta \\
\nonumber  \!\!\!\!\!\!&-&\!\!\!\!\!\!\ 5568y^ir_jr_{00}b\beta^3 +768y^ir_jr_{00}\beta^3b^2 +372b^iy_jr_{00}^2b^2\beta -672b^iy_jr_{00}s_0\beta^2 -1664y^iy_jr_{00}r_0b^5\beta \\
\nonumber  \!\!\!\!\!\!&-&\!\!\!\!\!\!\ 832y^iy_jr_{00}r_0b^4\beta +384b^iy_jr_{00}s_0b^2\beta^2 -2112b^iy_jr_{00}s_0b\beta^2 -768b^iy_jr_{00}s_0b^3\beta^2 -96b^iy_jr_{00}r_0b^2\beta^2 \\
\nonumber  \!\!\!\!\!\!&+&\!\!\!\!\!\!\ -2112b^iy_jr_{00}r_0b\beta^2 -768b^iy_jr_{00}r_0b^3\beta^2 -672y^ib_jr_{00}s_0b^2\beta^2 -2304y^ib_jr_{00}s_0b\beta^2 +351y^iy_jr_{00}^2b^2
\\
\nonumber  \!\!\!\!\!\!&+&\!\!\!\!\!\!\ 576y^ib_jr_{00}r_0b^2\beta^2 -2304y^ib_jr_{00}r_0b\beta^2 +144y^ir_{j0}r_{00}b^2\beta^2 -96y^ir_{j0}r_{00}\beta^2b^4 +5568y^ir_{j0}s_0b\beta^3 \\
\nonumber  \!\!\!\!\!\!&+&\!\!\!\!\!\!\ 5568y^ir_{j0}r_0b\beta^3 -576y^ir_{j0}s_0b^2\beta^3 +3072y^ir_{j0}s_0b^3\beta^3 -768y^ir_{j0}r_0b^2\beta^3 +3072y^ir_{j0}r_0b^3\beta^3 \\
\nonumber  \!\!\!\!\!\!&-&\!\!\!\!\!\!\ 96b^ir_{j0}r_{00}\beta^3b^2 -5568y^is_jr_{00}b\beta^3 +864y^is_jr_{00}\beta^3b^2 -3072y^is_jr_{00}b^3\beta^3 +528r^i_0y_jr_{00}b^2\beta^2 \\
\nonumber  \!\!\!\!\!\!&+&\!\!\!\!\!\!\ 96r^i_0y_jr_{00}\beta^2b^4 +648r^i_0y_js_0\beta^3 +402r^i_0y_jr_{00}\beta^2 -216r^i_0b_jr_{00}\beta^3 +729\beta^2s^i_0y_jr_{00} -1728\overline{R}^i_j\beta^4b^2\big]
\\
\nonumber t_6 \!\!\!\!&:=&\!\!\!\!\ -\beta\big[-4536\beta^2s^i_0y_jr_{00}b^2 -2592\beta^2s^i_0y_jr_{00}b^4 +9216y^ir_js_0b\beta^4 +8576\delta_j^ir_{00}s_0b^2\beta^3 \\
\nonumber  \!\!\!\!\!\!&-&\!\!\!\!\!\!\ 11072\delta_j^ir_{00}s_0b\beta^3 -12800\delta_j^ir_{00}s_0b^3\beta^3 -6400\delta_j^ir_{00}r_0b^2\beta^3 -11072\delta_j^ir_{00}r_0b\beta^3 \\
\nonumber  \!\!\!\!\!\!&+&\!\!\!\!\!\!\ 24y^iy_jr_{00}s_0\beta +1068y^iy_jr_{00}r_0\beta -1524y^iy_jr_{00|0}b^2\beta -2208y^iy_jr_{00|0}b^4\beta +2976y^iy_js_0^2b^2\beta^2 \\
\nonumber  \!\!\!\!\!\!&-&\!\!\!\!\!\!\ 7104y^iy_js_0^2b\beta^2 -6144y^iy_js_0^2b^3\beta^2 +7152y^iy_jr_0s_0\beta^2 -1152y^iy_js_ms^m_0b^2\beta^3 -20400y^iy_jr_{0m}s^m_0b^2\beta^2 \\
\nonumber  \!\!\!\!\!\!&-&\!\!\!\!\!\!\ 17664y^iy_jr_{0m}s^m_0b^4\beta^2 +1104y^iy_js_{0|0}b^2\beta^2 +384y^iy_js_{0|0}b^4\beta^2 +4096y^ir_jr_{00}b^5\beta^3 -288b^ib_jr_{00}s_0\beta^3\\
\nonumber  \!\!\!\!\!\!&+&\!\!\!\!\!\!\ 25600y^ir_jr_{00}b^3\beta^3 +1152b^iy_jr_{0m}s^m_0b^2\beta^3 -4608b^iy_js_{0|0}b^2\beta^3 +864b^ib_jr_{00}^2b^2\beta^2  -576b^ib_jr_{00}r_0\beta^3 \\
\nonumber  \!\!\!\!\!\!&+&\!\!\!\!\!\!\ 1152b^ib_jr_{00|0}b^2\beta^3 +512y^ib_jr_{00|0}b^6\beta^2 +576y^ib_jr_{00}^2b^4\beta -1512y^ib_jr_{00}^2b^2\beta +6056y^ib_jr_{00}s_0\beta^2 \\
\nonumber  \!\!\!\!\!\!&-&\!\!\!\!\!\!\ 2320y^ib_jr_{00}r_0\beta^2 +2496y^ib_jr_{00|0}b^2\beta^2 +2880y^ib_jr_{00|0}b^4\beta^2 +9216y^ib_js_0^2b\beta^3 -2304y^ib_jr_0s_0\beta^3 \\
\nonumber  \!\!\!\!\!\!&+&\!\!\!\!\!\!\ 5760y^ib_jr_{0m}s^m_0b^2\beta^3 -2304y^ib_js_{0|0}b^2\beta^3 +2160y^is_jr_{00}b^2\beta^2 +3456y^is_js_0b^2\beta^3 +2880g^iy_jr_{00}\beta^3b^2
\\
\nonumber  \!\!\!\!\!\!&+&\!\!\!\!\!\!\ 1296s^i_ks^k_0y_j\beta^4 -1728b^is_jr_{00}\beta^3b^2 -4608b^ir_jr_{00}b\beta^4 -330y^iy_jr_{00|0}\beta +1992y^iy_js_0^2\beta^2 \\
\nonumber  \!\!\!\!\!\!&-&\!\!\!\!\!\!\ 576y^iy_js_ms^m_0\beta^3 -6276y^iy_jr_{0m}s^m_0\beta^2 +564y^iy_js_{0|0}\beta^2 -4464b^iy_js_{0|0}\beta^3 -108b^iy_jr_{00}^2\beta \\
\nonumber  \!\!\!\!\!\!&-&\!\!\!\!\!\!\ 684b^iy_jr_{00|0}\beta^2 -864b^iy_js_0^2\beta^3 +144b^iy_jr_{0m}s^m_0\beta^3 -2376b^ib_jr_{00}^2\beta^2 +1872b^ib_jr_{00|0}\beta^3 \\
\nonumber  \!\!\!\!\!\!&-&\!\!\!\!\!\!\ 2880y^ib_js_{0|0}\beta^3 +162y^ib_jr_{00}^2\beta +268y^ib_jr_{00|0}\beta^2 +3456r_j^ir_{00}b^2\beta^4 +1152b^ir_jr_{00}\beta^4 \\
\nonumber  \!\!\!\!\!\!&-&\!\!\!\!\!\!\ 5536y^ir_jr_{00}\beta^3 -2304y^ir_js_0\beta^4 +3312y^ib_jr_{0m}s^m_0\beta^3 +1440b^is_jr_{00}\beta^4 -6784y^is_jr_{00}\beta^3
\\
\nonumber  \!\!\!\!\!\!&-&\!\!\!\!\!\!\ 1728y^is_js_0\beta^4 +1152y^is_jr_0\beta^4 +1944s^i_0b_jr_{00}\beta^3 +512\delta_j^ir_{00|0}b^6\beta^3 +3168\delta_j^ir_{00}^2b^4\beta^2 \\
\nonumber  \!\!\!\!\!\!&+&\!\!\!\!\!\!\ 5424\delta_j^ir_{00}^2b^2\beta^2 +2768\delta_j^ir_{00}s_0\beta^3 -4608\delta_j^is_0^2b\beta^4 -2304\delta_j^ir_0s_0\beta^4 +6912\delta_j^ir_{0m}s^m_0b^2\beta^4 \\
\nonumber  \!\!\!\!\!\!&+&\!\!\!\!\!\!\ 3456\delta_j^is_{0|0}b^2\beta^4 +3040y^ir_{j0}s_0\beta^3 -840y^ir_{j0}r_{00}\beta^2 +5536y^ir_{j0}r_0\beta^3 -16608y^ir_{00|j}b^2\beta^3 \\
\nonumber  \!\!\!\!\!\!&-&\!\!\!\!\!\!\ 9600y^ir_{00|j}\beta^3b^4 -1024y^ir_{00|j}b^6\beta^3 +16608y^ir_{j0|0}b^2\beta^3 +9600y^ir_{j0|0}\beta^3b^4 +1024y^ir_{j0|0}b^6\beta^3\\
\nonumber  \!\!\!\!\!\!&-&\!\!\!\!\!\!\ 3456y^ir_{k0}s^k_jb^2\beta^4 +6912y^ir_{kj}s^k_0b^2\beta^4 +3456y^is_{j|0}b^2\beta^4 +23040s_{0|0}^iy_jb^2\beta^3 +18432s_{0|0}^iy_j\beta^3b^4 \\
\nonumber  \!\!\!\!\!\!&+&\!\!\!\!\!\!\ 3072s_{0|0}^iy_jb^6\beta^3 -6912y^is_{0|j}b^2\beta^4 +2592b^ir_{j0}r_{00}\beta^3 -576b^ir_{j0}s_0\beta^4 -1152b^ir_{j0}r_0\beta^4
\\
\nonumber  \!\!\!\!\!\!&-&\!\!\!\!\!\!\ 3456b^ir_{j0|0}b^2\beta^4 +3456b^ir_{00|j}b^2\beta^4 +5616y^is_js_0\beta^3 +2052y^is_jr_{00}\beta^2 +2304g^iy_jr_{00}\beta^3
\\
\nonumber  \!\!\!\!\!\!&-&\!\!\!\!\!\!\ 12800\delta_j^ir_{00}r_0b^3\beta^3 -4608\delta_j^ir_0s_0b\beta^4 -2048\delta_j^ir_{00}s_0b^5\beta^3 +2048\delta_j^ir_{00}s_0b^4\beta^3 -2048\delta_j^ir_{00}r_0b^5\beta^3
\\
\nonumber  \!\!\!\!\!\!&+&\!\!\!\!\!\!\ 8856\overline{R}^i_j\beta^4 +2931336s^i_0s_{0j}\beta^4 +1296s_{j|0}^i\beta^5 -2592s_{j|0}^i\beta^5 +3024y^iy_jr_{00}r_0b^2\beta +4608y^iy_jr_{00}s_0b^3\beta \\
\nonumber  \!\!\!\!\!\!&-&\!\!\!\!\!\!\ 1584y^iy_jr_{00}s_0b^2\beta +1200y^iy_jr_{00}s_0b\beta +1200y^iy_jr_{00}r_0b\beta +4608y^iy_jr_{00}r_0b^3\beta +12096y^iy_jr_0s_0b^2\beta^2\\
\nonumber  \!\!\!\!\!\!&-&\!\!\!\!\!\!\ 7104y^iy_jr_0s_0b\beta^2 -6144y^iy_jr_0s_0b^3\beta^2 +3840y^iy_jr_{00}s_0b^5\beta -3264y^iy_jr_{00}s_0b^4\beta -5536\delta_j^ir_{00}r_0\beta^3\\
\nonumber  \!\!\!\!\!\!&+&\!\!\!\!\!\!\ 8304\delta_j^ir_{00|0}b^2\beta^3 +4800\delta_j^ir_{00|0}b^4\beta^3 +5719680s^i_0s_{0j}b^2\beta^4 +2287872s^i_0s_{0j}\beta^4b^4
\\
\nonumber  \!\!\!\!\!\!&-&\!\!\!\!\!\!\ 4320r^i_0r_{j0}\beta^4 -7776y^is_{0|j}\beta^4 -4320b^ir_{j0|0}\beta^4 +4320b^ir_{00|j}\beta^4 -1296s_{0|0}^ib_j\beta^4 +4320r_j^ir_{00}\beta^4 \\
\nonumber  \!\!\!\!\!\!&+&\!\!\!\!\!\!\ 7776\delta_j^ir_{0m}s^m_0\beta^4 +3888\delta_j^is_{0|0}\beta^4 -7112y^ir_{00|j}\beta^3 +7112y^ir_{j0|0}\beta^3 -3888y^ir_{k0}s^k_j\beta^4
\end{eqnarray}
\begin{eqnarray}
\nonumber  \!\!\!\!\!\!&+&\!\!\!\!\!\!\ 7776y^ir_{kj}s^k_0\beta^4 +3888y^is_{j|0}\beta^4 +7620s_{0|0}^iy_j\beta^3 +1506\delta_j^ir_{00}^2\beta^2 +3556\delta_j^ir_{00|0}\beta^3
\\
\nonumber  \!\!\!\!\!\!&-&\!\!\!\!\!\!\ 3456b^is_jr_{00}\beta^3 -576b^iy_jr_0s_0\beta^3 +9216b^iy_js_0^2b\beta^3 +720b^iy_jr_{00}r_0\beta^2 -1440b^iy_jr_{00|0}b^2\beta^2 \\
\nonumber  \!\!\!\!\!\!&-&\!\!\!\!\!\!\ 576b^iy_jr_{00|0}b^4\beta^2 -960y^iy_jr_{00|0}b^6\beta +22144y^ir_jr_{00}b\beta^3 -6400y^ir_jr_{00}\beta^3b^2 -1024y^ir_jr_{00}\beta^3b^4
\\
\nonumber  \!\!\!\!\!\!&-&\!\!\!\!\!\!\ 648b^iy_jr_{00}^2b^2\beta +3816b^iy_jr_{00}s_0\beta^2 +3840y^iy_jr_{00}r_0b^5\beta +1920y^iy_jr_{00}r_0b^4\beta +3072y^iy_jr_0s_0b^4\beta^2 \\
\nonumber  \!\!\!\!\!\!&+&\!\!\!\!\!\!\ 288b^iy_jr_{00}s_0b^2\beta^2 +5760b^iy_jr_{00}s_0b\beta^2 +4608b^iy_jr_{00}s_0b^3\beta^2 +576b^iy_jr_{00}r_0b^2\beta^2 +5760b^iy_jr_{00}r_0b\beta^2
\\
\nonumber  \!\!\!\!\!\!&+&\!\!\!\!\!\!\ 4608b^iy_jr_{00}r_0b^3\beta^2 +9216b^iy_jr_0s_0b\beta^3 -4608b^ib_jr_{00}s_0b\beta^3 -4608b^ib_jr_{00}r_0b\beta^3 +2336y^ib_jr_{00}s_0b^2\beta^2 \\
\nonumber  \!\!\!\!\!\!&+&\!\!\!\!\!\!\ 7168y^ib_jr_{00}s_0b\beta^2 -512y^ib_jr_{00}s_0b^3\beta^2 -4288y^ib_jr_{00}r_0b^2\beta^2 +7168y^ib_jr_{00}r_0b\beta^2 -512y^ib_jr_{00}r_0b^3\beta^2 \\
\nonumber  \!\!\!\!\!\!&+&\!\!\!\!\!\!\ 9216y^ib_jr_0s_0b\beta^3 -2048y^ib_jr_{00}s_0b^5\beta^2 +2048y^ib_jr_{00}s_0b^4\beta^2 -2048y^ib_jr_{00}r_0b^5\beta^2 \\
\nonumber  \!\!\!\!\!\!&-&\!\!\!\!\!\!\ 1024y^ib_jr_{00}r_0b^4\beta^2 -1024\delta_j^ir_{00}r_0b^4\beta^3 +288y^ir_{j0}r_{00}b^2\beta^2 +768y^ir_{j0}r_{00}\beta^2b^4 -22144y^ir_{j0}s_0b\beta^3 \\
\nonumber  \!\!\!\!\!\!&-&\!\!\!\!\!\!\ 22144y^ir_{j0}r_0b\beta^3 +5440y^ir_{j0}s_0b^2\beta^3 -25600y^ir_{j0}s_0b^3\beta^3 +6400y^ir_{j0}r_0b^2\beta^3 -25600y^ir_{j0}r_0b^3\beta^3 \\
\nonumber  \!\!\!\!\!\!&+&\!\!\!\!\!\!\ 1024y^ir_{j0}s_0\beta^3b^4 -4096y^ir_{j0}s_0b^5\beta^3 +1024y^ir_{j0}r_0b^4\beta^3 -4096y^ir_{j0}r_0b^5\beta^3 +1152b^ir_{j0}r_{00}\beta^3b^2 \\
\nonumber  \!\!\!\!\!\!&+&\!\!\!\!\!\!\ 4608b^ir_{j0}s_0b\beta^4 +4608b^ir_{j0}r_0b\beta^4 -4608b^is_jr_{00}b\beta^4 +22144y^is_jr_{00}b\beta^3 -6880y^is_jr_{00}\beta^3b^2 \\
\nonumber  \!\!\!\!\!\!&-&\!\!\!\!\!\!\ 1024y^is_jr_{00}\beta^3b^4 +25600y^is_jr_{00}b^3\beta^3 +4096y^is_jr_{00}b^5\beta^3 +4608y^is_js_0b\beta^4 -4608y^is_jr_0b\beta^4 \\
\nonumber  \!\!\!\!\!\!&-&\!\!\!\!\!\!\ 3072y^iy_jr_{0m}s^m_0b^6\beta^2 -1440r^i_0y_jr_{00}b^2\beta^2 -576r^i_0y_jr_{00}\beta^2b^4 -4608r^i_0y_js_0b^2\beta^3 +1152r^i_0b_jr_{00}\beta^3b^2 \\
\nonumber  \!\!\!\!\!\!&-&\!\!\!\!\!\!\ 3456r^i_0r_{j0}b^2\beta^4 -4464r^i_0y_js_0\beta^3 -684r^i_0y_jr_{00}\beta^2 +1872r^i_0b_jr_{00}\beta^3 -1944\beta^3s^i_0y_js_0 \\
\nonumber  \!\!\!\!\!\!&-&\!\!\!\!\!\!\ 1620\beta^2s^i_0y_jr_{00} +6912\overline{R}^i_j\beta^4b^4 -27y^iy_jr_{00}^2 -504y^iy_jr_{00}^2b^4 -324y^iy_jr_{00}^2b^2 +17280\overline{R}^i_j\beta^4b^2\big]
\\
\nonumber t_7 \!\!\!\!&:=&\!\!\!\!\  -1728\beta^2s^i_0y_jr_{00}b^2-3888\beta^3s^i_0y_js_0b^2-2160\beta^2s^i_0y_jr_{00}b^4-576\beta^2s^i_0y_jr_{00}b^6  -3072y^ir_js_0b^2\beta^4 \\
\nonumber  \!\!\!\!\!\!&+&\!\!\!\!\!\!\ 17664y^ir_js_0b\beta^4 +12288y^ir_js_0b^3\beta^4 +6656\delta_j^ir_{00}s_0b^2\beta^3 -5792\delta_j^ir_{00}s_0b\beta^3 -10496\delta_j^ir_{00}s_0b^3\beta^3 \\
\nonumber  \!\!\!\!\!\!&-&\!\!\!\!\!\!\ 5248\delta_j^ir_{00}r_0b^2\beta^3 -5792\delta_j^ir_{00}r_0b\beta^3 -1536y^iy_js_ms^m_0b^4\beta^3 +16y^iy_jr_{00}s_0\beta -9120y^iy_jr_{0m}s^m_0b^4\beta^2
\\
\nonumber  \!\!\!\!\!\!&+&\!\!\!\!\!\!\ 148y^iy_jr_{00}r_0\beta-216y^iy_jr_{00|0}b^2\beta-384y^iy_jr_{00|0}b^4\beta +2272y^iy_js_0^2b^2\beta^2 -1600y^iy_js_0^2b\beta^2
\\
\nonumber  \!\!\!\!\!\!&-&\!\!\!\!\!\!\ 1792y^iy_js_0^2b^3\beta^2 +2416y^iy_jr_0s_0\beta^2 -2400y^iy_js_ms^m_0b^2\beta^3 -6984y^iy_jr_{0m}s^m_0b^2\beta^2 +144y^iy_js_{0|0}b^2\beta^2 \\
\nonumber  \!\!\!\!\!\!&-&\!\!\!\!\!\!\ 96y^iy_js_{0|0}b^4\beta^2 +7168y^ir_jr_{00}b^5\beta^3 +20992y^ir_jr_{00}b^3\beta^3+1824b^iy_jr_{0m}s^m_0b^2\beta^3 +1536b^iy_jr_{0m}s^m_0b^4\beta^3 \\
\nonumber  \!\!\!\!\!\!&-&\!\!\!\!\!\!\ 6432b^iy_js_{0|0}b^2\beta^3-2688b^iy_js_{0|0}b^4\beta^3+144b^ib_jr_{00}^2b^2\beta^2-5136b^ib_jr_{00}s_0\beta^3-1056b^ib_jr_{00}r_0\beta^3\\
\nonumber  \!\!\!\!\!\!&+&\!\!\!\!\!\!\ 2112b^ib_jr_{00|0}b^2\beta^3+384b^ib_jr_{00|0}b^4\beta^3+768y^ib_jr_{00|0}b^6\beta^2 -96y^ib_jr_{00}^2b^4\beta -684y^ib_jr_{00}^2b^2\beta
\\
\nonumber  \!\!\!\!\!\!&+&\!\!\!\!\!\!\  1908y^ib_jr_{00}s_0\beta^2-1032y^ib_jr_{00}r_0\beta^2+1152y^ib_jr_{00|0}b^2\beta^2+2016y^ib_jr_{00|0}b^4\beta^2 -1152y^ib_js_0^2b^2\beta^3 \\
\nonumber  \!\!\!\!\!\!&+&\!\!\!\!\!\!\  13440y^ib_js_0^2b\beta^3 +6144y^ib_js_0^2b^3\beta^3 -3648y^ib_jr_0s_0\beta^3 +9120y^ib_jr_{0m}s^m_0b^2\beta^3+5376y^ib_jr_{0m}s^m_0b^4\beta^3 \\
\nonumber  \!\!\!\!\!\!&-&\!\!\!\!\!\!\ 2784y^ib_js_{0|0}b^2\beta^3-768y^ib_js_{0|0}b^4\beta^3 -36y^is_jr_{00}b^2\beta^2 -576y^is_jr_{00}b^6\beta^2-1440y^is_jr_{00}\beta^2b^4
\\
\nonumber  \!\!\!\!\!\!&+&\!\!\!\!\!\!\ 5904y^is_js_0b^2\beta^3+1152y^is_js_0\beta^3b^4 +4128s^iy_jr_{00}\beta^3b^2 +2052s^i_ks^k_0y_j\beta^4 +2112s^iy_jr_{00}\beta^3b^4 \\
\nonumber  \!\!\!\!\!\!&-&\!\!\!\!\!\!\  4032b^is_jr_{00}\beta^3b^2 -576b^is_jr_{00}\beta^3b^4 +1536b^ir_jr_{00}b^2\beta^4 -9984b^ir_jr_{00}b\beta^4 -128y^iy_js_0^2b^4\beta^2 \\
\nonumber  \!\!\!\!\!\!&-&\!\!\!\!\!\!\ 40y^iy_jr_{00|0}\beta+1096y^iy_js_0^2\beta^2 -708y^iy_js_ms^m_0\beta^3- 1628y^iy_jr_{0m}s^m_0\beta^2 +80y^iy_js_{0|0}\beta^2\\
\nonumber  \!\!\!\!\!\!&-&\!\!\!\!\!\!\ 2976b^iy_js_{0|0}\beta^3 -72b^iy_jr_{00}^2\beta-144b^iy_jr_{00|0}\beta^2 +1488b^iy_js_0^2\beta^3-1728b^iy_js_ms^m_0\beta^4
\\
\nonumber  \!\!\!\!\!\!&+&\!\!\!\!\!\!\ 312b^iy_jr_{0m}s^m_0\beta^3 -1836b^ib_jr_{00}^2\beta^2+1608b^ib_jr_{00|0}\beta^3 +864b^ib_jr_{0m}s^m_0\beta^4 +1728b^ib_js_{0|0}\beta^4 \\
\nonumber  \!\!\!\!\!\!&-&\!\!\!\!\!\!\ 1848y^ib_js_{0|0}\beta^3-39y^ib_jr_{00}^2\beta+114y^ib_jr_{00|0}\beta^2 +2016y^ib_js_0^2\beta^3 -432s^ib_jr_{00}\beta^4 +7488r_j^ir_{00}b^2\beta^4 \\
\nonumber  \!\!\!\!\!\!&+&\!\!\!\!\!\!\ 2304r_j^ir_{00}\beta^4b^4 +2496b^ir_jr_{00}\beta^4 -2896y^ir_jr_{00}\beta^3-4416y^ir_js_0\beta^4 -864y^ib_js_ms^m_0\beta^4 \\
\nonumber  \!\!\!\!\!\!&+&\!\!\!\!\!\!\ 2568y^ib_jr_{0m}s^m_0\beta^3 +3312b^is_jr_{00}\beta^4 -3472y^is_jr_{00}\beta^3 -3264y^is_js_0\beta^4 +2208y^is_jr_0\beta^4 \\
\nonumber  \!\!\!\!\!\!&+&\!\!\!\!\!\!\ 3456s^i_ks^k_0y_jb^2\beta^4+2916s^i_0b_jr_{00}\beta^3+2304\delta_j^is_{0|0}\beta^4b^4 +4608\delta_j^ir_{0m}s^m_0b^4\beta^4 +896\delta_j^ir_{00|0}b^6\beta^3 \\
\nonumber  \!\!\!\!\!\!&+&\!\!\!\!\!\!\  768\delta_j^is_0^2b^2\beta^4 -6144\delta_j^is_0^2b^3\beta^4 +1488\delta_j^ir_{00}^2b^4\beta^2 +1752\delta_j^ir_{00}^2b^2\beta^2  +992\delta_j^ir_{00}s_0\beta^3 -8832\delta_j^is_0^2b\beta^4 \\
\nonumber  \!\!\!\!\!\!&-&\!\!\!\!\!\!\ 4416\delta_j^ir_0s_0\beta^4+13248\delta_j^ir_{0m}s^m_0b^2\beta^4+6624\delta_j^is_{0|0}b^2\beta^4+1744y^ir_{j0}s_0\beta^3-156y^ir_{j0}r_{00}\beta^2 \end{eqnarray}
\begin{eqnarray}
\nonumber  \!\!\!\!\!\!&+&\!\!\!\!\!\!\ 2896y^ir_{j0}r_0\beta^3 -8688y^ir_{00|j}b^2\beta^3-7872y^ir_{00|j}\beta^3b^4 -1792y^ir_{00|j}b^6\beta^3 +8688y^ir_{j0|0}b^2\beta^3 \\
\nonumber  \!\!\!\!\!\!&+&\!\!\!\!\!\!\ 7872y^ir_{j0|0}\beta^3b^4+1792y^ir_{j0|0}b^6\beta^3 -6624y^ir_{k0}s^k_jb^2\beta^4 -2304y^ir_{k0}s^k_j\beta^4b^4 +13248y^ir_{jk}s^k_0b^2\beta^4 \\
\nonumber  \!\!\!\!\!\!&+&\!\!\!\!\!\!\ 4608y^ir_{jk}s^k_0\beta^4b^4+6624y^is_{j|0}b^2\beta^4+2304y^is_{j|0}\beta^4b^4+12640s_{0|0}^iy_jb^2\beta^3+16896s_{0|0}^iy_j\beta^3b^4 \\
\nonumber  \!\!\!\!\!\!&+&\!\!\!\!\!\!\ 7168s_{0|0}^iy_jb^6\beta^3+512s_{0|0}^iy_j\beta^3b^8-13248y^is_{0|j}b^2\beta^4 -4608y^is_{0|j}\beta^4b^4+1824b^ir_{j0}r_{00}\beta^3
\\
\nonumber  \!\!\!\!\!\!&-&\!\!\!\!\!\!\ 864b^ir_{j0}s_0\beta^4 -2496b^ir_{j0}r_0\beta^4-7488b^ir_{j0|0}b^2\beta^4-2304b^ir_{j0|0}\beta^4b^4+7488b^ir_{00|j}b^2\beta^4 \\
\nonumber  \!\!\!\!\!\!&+&\!\!\!\!\!\!\ 2304b^ir_{00|j}\beta^4b^4 +4284y^is_js_0\beta^3 +432y^is_jr_{00}\beta^2 +1644s^iy_jr_{00}\beta^3+432s^iy_js_0\beta^4-10496\delta_j^ir_{00}r_0b^3\beta^3 \\
\nonumber  \!\!\!\!\!\!&-&\!\!\!\!\!\!\ 8832\delta_j^ir_0s_0b\beta^4 -3584\delta_j^ir_{00}s_0b^5\beta^3 +3584\delta_j^ir_{00}s_0b^4\beta^3 -3584\delta_j^ir_{00}r_0b^5\beta^3 +5892\overline{R}^i_j\beta^4 \\
\nonumber  \!\!\!\!\!\!&+&\!\!\!\!\!\!\ 1950252s^i_0s_{0j}\beta^4+3024s_{j|0}^i\beta^5-6048s_{j|0}^i\beta^5+520y^iy_jr_{00}r_0b^2\beta +896y^iy_jr_{00}s_0b^3\beta \\
\nonumber  \!\!\!\!\!\!&-&\!\!\!\!\!\!\  524y^iy_jr_{00}s_0b^2\beta +224y^iy_jr_{00}s_0b\beta+224y^iy_jr_{00}r_0b\beta+896y^iy_jr_{00}r_0b^3\beta +6208y^iy_jr_0s_0b^2\beta^2 \\
\nonumber  \!\!\!\!\!\!&-&\!\!\!\!\!\!\ 1600y^iy_jr_0s_0b\beta^2-1792y^iy_jr_0s_0b^3\beta^2 +896y^iy_jr_{00}s_0b^5\beta -1040y^iy_jr_{00}s_0b^4\beta -2896\delta_j^ir_{00}r_0\beta^3 \\
\nonumber  \!\!\!\!\!\!&+&\!\!\!\!\!\!\ 4344\delta_j^ir_{00|0}b^2\beta^3+3936\delta_j^ir_{00|0}b^4\beta^3-864r^i_0s_j\beta^5 +5910336s^i_0s_{0j}b^2\beta^4+4957056s^i_0s_{0j}\beta^4b^4 \\
\nonumber  \!\!\!\!\!\!&+&\!\!\!\!\!\!\ 1016832s^i_0s_{0j}\beta^4b^6-4464r^i_0r_{j0}\beta^4-7200y^is_{0|j}\beta^4-864b^is_{j|0}\beta^5-864y^is_ms^m_0\beta^5+3456s_{j|0}^ib^2\beta^5 \\
\nonumber  \!\!\!\!\!\!&-&\!\!\!\!\!\!\ 6912s_{j|0}^ib^2\beta^5-4464b^ir_{j0|0}\beta^4+4464b^ir_{00|j}\beta^4-2376s_{0|0}^ib_j\beta^4+4464r_j^ir_{00}\beta^4+1728r_j^is_0\beta^5 \\
\nonumber  \!\!\!\!\!\!&+&\!\!\!\!\!\!\ 7200\delta_j^ir_{0m}s^m_0\beta^4+3600\delta_j^is_{0|0}\beta^4 +1728\delta_j^is_ms^m_0\beta^5 -2708y^ir_{00|j}\beta^3 +2708y^ir_{j0|0}\beta^3   -3600y^ir_{k0}s^k_j\beta^4 \\
\nonumber  \!\!\!\!\!\!&+&\!\!\!\!\!\!\ 7200y^ir_{jk}s^k_0\beta^4+3600y^is_{j|0}\beta^4+2960s_{0|0}^iy_j\beta^3-768\delta_j^is_0^2\beta^4 +351\delta_j^ir_{00}^2\beta^2 +1354\delta_j^ir_{00|0}\beta^3
\\
\nonumber  \!\!\!\!\!\!&+&\!\!\!\!\!\!\ 1728b^is_{0|j}\beta^5+864b^ir_{k0}s^k_j\beta^5-1728b^ir_{jk}s^k_0\beta^5 -864s^ir_{j0}\beta^5-3492b^is_jr_{00}\beta^3 -1296b^is_js_0\beta^4 \\
\nonumber  \!\!\!\!\!\!&-&\!\!\!\!\!\!\ 3456s_{0|0}^ib_jb^2\beta^4+9216b^iy_js_0^2b^3\beta^3 -864b^iy_jr_0s_0\beta^3 -1920b^iy_js_0^2b^2\beta^3 +12672b^iy_js_0^2b\beta^3 \\
\nonumber  \!\!\!\!\!\!&+&\!\!\!\!\!\!\ 216b^iy_jr_{00}r_0\beta^2-432b^iy_jr_{00|0}b^2\beta^2-288b^iy_jr_{00|0}b^4\beta^2 -224y^iy_jr_{00|0}b^6\beta +512y^iy_js_0^2b^5\beta^2
\\
\nonumber  \!\!\!\!\!\!&-&\!\!\!\!\!\!\ 6144b^ir_jr_{00}b^3\beta^4 +11584y^ir_jr_{00}b\beta^3 -5248y^ir_jr_{00}\beta^3b^2 -1792y^ir_jr_{00}\beta^3b^4 -180b^iy_jr_{00}^2b^2\beta \\
\nonumber  \!\!\!\!\!\!&+&\!\!\!\!\!\!\ 1404b^iy_jr_{00}s_0\beta^2 +896y^iy_jr_{00}r_0b^5\beta +448y^iy_jr_{00}r_0b^4\beta +512y^iy_jr_0s_0b^5\beta^2 +3328y^iy_jr_0s_0b^4\beta^2 \\
\nonumber  \!\!\!\!\!\!&+&\!\!\!\!\!\!\ 864b^iy_jr_{00}s_0b^2\beta^2 +1728b^iy_jr_{00}s_0b\beta^2 +2304b^iy_jr_{00}s_0b^3\beta^2 +288b^iy_jr_{00}r_0b^2\beta^2+1728b^iy_jr_{00}r_0b\beta^2 \\
\nonumber  \!\!\!\!\!\!&+&\!\!\!\!\!\!\ 2304b^iy_jr_{00}r_0b^3\beta^2-1152b^iy_jr_0s_0b^2\beta^3+12672b^iy_jr_0s_0b\beta^3 +9216b^iy_jr_0s_0b^3\beta^3 +576b^iy_jr_{00}s_0b^4\beta^2 \\
\nonumber  \!\!\!\!\!\!&+&\!\!\!\!\!\!\  960b^ib_jr_{00}s_0b^2\beta^3-8448b^ib_jr_{00}s_0b\beta^3 -3072b^ib_jr_{00}s_0b^3\beta^3 -384b^ib_jr_{00}r_0b^2\beta^3 -8448b^ib_jr_{00}r_0b\beta^3 \\
\nonumber  \!\!\!\!\!\!&-&\!\!\!\!\!\!\  3072b^ib_jr_{00}r_0b^3\beta^3-1440y^ib_jr_{00}s_0b^2\beta^2 +2688y^ib_jr_{00}s_0b\beta^2 -768y^ib_jr_{00}s_0b^3\beta^2 -2976y^ib_jr_{00}r_0b^2\beta^2 \\
\nonumber  \!\!\!\!\!\!&+&\!\!\!\!\!\!\  2688y^ib_jr_{00}r_0b\beta^2 -768y^ib_jr_{00}r_0b^3\beta^2 -3840y^ib_jr_0s_0b^2\beta^3 +13440y^ib_jr_0s_0b\beta^3 +6144y^ib_jr_0s_0b^3\beta^3
\\
\nonumber  \!\!\!\!\!\!&-&\!\!\!\!\!\!\ 3072y^ib_jr_{00}s_0b^5\beta^2 +1728y^ib_jr_{00}s_0b^4\beta^2 -3072y^ib_jr_{00}r_0b^5\beta^2-1536y^ib_jr_{00}r_0b^4\beta^2 -1792\delta_j^ir_{00}r_0b^4\beta^3  \\
\nonumber  \!\!\!\!\!\!&-&\!\!\!\!\!\!\ 3072\delta_j^ir_0s_0b^2\beta^4 -6144\delta_j^ir_0s_0b^3\beta^4 +336y^ir_{j0}r_{00}b^2\beta^2 +576y^ir_{j0}r_{00}\beta^2b^4 -11584y^ir_{j0}s_0b\beta^3 \\
\nonumber  \!\!\!\!\!\!&-&\!\!\!\!\!\!\ 11584y^ir_{j0}r_0b\beta^3 +5152y^ir_{j0}s_0b^2\beta^3-20992y^ir_{j0}s_0b^3\beta^3 +5248y^ir_{j0}r_0b^2\beta^3-20992y^ir_{j0}r_0b^3\beta^3 \\
\nonumber  \!\!\!\!\!\!&+&\!\!\!\!\!\!\ 2176y^ir_{j0}s_0\beta^3b^4-7168y^ir_{j0}s_0b^5\beta^3+1792y^ir_{j0}r_0b^4\beta^3 -7168y^ir_{j0}r_0b^5\beta^3+1344b^ir_{j0}r_{00}\beta^3b^2\\
\nonumber  \!\!\!\!\!\!&-&\!\!\!\!\!\!\ 1152b^ir_{j0}s_0b^2\beta^4 +9984b^ir_{j0}s_0b\beta^4+6144b^ir_{j0}s_0b^3\beta^4 -1536b^ir_{j0}r_0b^2\beta^4 +9984b^ir_{j0}r_0b\beta^4 \\
\nonumber  \!\!\!\!\!\!&+&\!\!\!\!\!\!\ 6144b^ir_{j0}r_0b^3\beta^4+1728b^is_jr_{00}b^2\beta^4 -9984b^is_jr_{00}b\beta^4 -6144b^is_jr_{00}b^3\beta^4 +11584y^is_jr_{00}b\beta^3 \\
\nonumber  \!\!\!\!\!\!&-&\!\!\!\!\!\!\ 5296y^is_jr_{00}\beta^3b^2 -1600y^is_jr_{00}\beta^3b^4 +20992y^is_jr_{00}b^3\beta^3 +7168y^is_jr_{00}b^5\beta^3-1920y^is_js_0b^2\beta^4\\
\nonumber  \!\!\!\!\!\!&+&\!\!\!\!\!\!\ 8832y^is_js_0b\beta^4+6144y^is_js_0b^3\beta^4+1536y^is_jr_0b^2\beta^4 -8832y^is_jr_0b\beta^4-6144y^is_jr_0b^3\beta^4 \\
\nonumber  \!\!\!\!\!\!&-&\!\!\!\!\!\!\ 128y^iy_js_{0|0}\beta^2b^6-3328y^iy_jr_{0m}s^m_0b^6\beta^2 +3888s^i_0b_jr_{00}\beta^3b^2-432r^i_0y_jr_{00}b^2\beta^2 +17856\overline{R}^i_j\beta^4b^2\\
\nonumber  \!\!\!\!\!\!&-&\!\!\!\!\!\!\ 288r^i_0y_jr_{00}\beta^2b^4-6432r^i_0y_js_0b^2\beta^3 -2688r^i_0y_js_0\beta^3b^4 +2112r^i_0b_jr_{00}\beta^3b^2+384r^i_0b_jr_{00}\beta^3b^4\\
\nonumber  \!\!\!\!\!\!&-&\!\!\!\!\!\!\ 7488r^i_0r_{j0}b^2\beta^4-2304r^i_0r_{j0}\beta^4b^4-2976r^i_0y_js_0\beta^3 -144r^i_0y_jr_{00}\beta^2+1608r^i_0b_jr_{00}\beta^3 +1728r^i_0b_js_0\beta^4 \\
\nonumber  \!\!\!\!\!\!&-&\!\!\!\!\!\!\ 1944\beta^3s^i_0y_js_0-396\beta^2s^i_0y_jr_{00}+14976\overline{R}^i_j\beta^4b^4 +3072\overline{R}^i_j\beta^4b^6-3y^iy_jr_{00}^2  -60y^iy_jr_{00}^2b^4-36y^iy_jr_{00}^2b^2
\end{eqnarray}
\begin{eqnarray}
\nonumber t_8 \!\!\!\!&:=&\!\!\!\!\  648s^i_ks^k_j\beta^5 -773216s^i_0s_{0j}\beta^3 -2916s_{j|0}^i\beta^4 +5832s_{j|0}^i\beta^4 -1620s^i_0b_jr_{00}\beta^2 +1944y^is_jr_{00}b^2\beta^2 \\
\nonumber  \!\!\!\!\!\!&+&\!\!\!\!\!\!\ 864y^is_jr_{00}\beta^2b^4 +2720y^is_js_0b^2\beta^3 +512y^is_js_0\beta^3b^4 -3296b^is_jr_{00}\beta^3b^2 -512b^is_jr_{00}\beta^3b^4 \\
\nonumber  \!\!\!\!\!\!&+&\!\!\!\!\!\!\ 2384y^is_js_0\beta^3 +972y^is_jr_{00}\beta^2 -1864192s^i_0s_{0j}b^6\beta^3 -169472s^i_0s_{0j}\beta^3b^8 +2440r^i_0r_{j0}\beta^3 -3456b^is_{0|j}\beta^4 \\
\nonumber  \!\!\!\!\!\!&-&\!\!\!\!\!\!\  3456r_j^is_0\beta^4-3024\delta_j^is_ms^m_0\beta^4+1756y^ir_{k0}s^k_j\beta^3 -2440b^ir_{00|j}\beta^3 +1728s_{0|0}^ib_j\beta^3-1728b^ir_{k0}s^k_j\beta^4 \\
\nonumber  \!\!\!\!\!\!&-&\!\!\!\!\!\!\  306\delta_j^ir_{00|0}\beta^2+1000\delta_j^is_0^2\beta^3 -3512y^ir_{jk}s^k_0\beta^3-1756y^is_{j|0}\beta^3-664s_{0|0}^iy_j\beta^2  -1756\delta_j^is_{0|0}\beta^3
\\
\nonumber  \!\!\!\!\!\!&-&\!\!\!\!\!\!\  48\delta_j^ir_{00}^2\beta+3512y^is_{0|j}\beta^3+1728b^is_{j|0}\beta^4+1512y^is_ms^m_0\beta^4+2y^iy_jr_{00|0}+12b^iy_jr_{00}^2+9y^ib_jr_{00}^2 \\
\nonumber  \!\!\!\!\!\!&-&\!\!\!\!\!\!\ 2440r_j^ir_{00}\beta^3+612y^ir_{00|j}\beta^2-612y^ir_{j0|0}\beta^2+6912s_{j|0}^i\beta^4b^4 +2440b^ir_{j0|0}\beta^3-6912s_{j|0}^ib^2\beta^4\\
\nonumber  \!\!\!\!\!\!&-&\!\!\!\!\!\!\ 3456s_{j|0}^i\beta^4b^4+3456b^ir_{jk}s^k_0\beta^4+13824s_{j|0}^ib^2\beta^4+1728g^ir_{j0}\beta^4-3512\delta_j^ir_{0m}s^m_0\beta^3-3648y^is_{j|0}\beta^3b^4\\
\nonumber  \!\!\!\!\!\!&-&\!\!\!\!\!\!\ 512y^is_{j|0}b^6\beta^3-3680s_{0|0}^iy_jb^2\beta^2+3456b^is_js_0\beta^3 +1692b^is_jr_{00}\beta^2-1476y^is_js_0\beta^2 -36y^is_jr_{00}\beta \\
\nonumber  \!\!\!\!\!\!&-&\!\!\!\!\!\!\ 684b^ib_jr_{00|0}\beta^2 -1584b^ib_jr_{0m}s^m_0\beta^3-1008y^ib_jr_{0m}s^m_0\beta^2+1024y^is_{0|j}b^6\beta^3+368b^ir_{j0}s_0\beta^3 \\
\nonumber  \!\!\!\!\!\!&+&\!\!\!\!\!\!\  564y^ib_js_{0|0}\beta^2-8y^iy_jr_{00}r_0-236y^iy_js_0^2\beta +316y^iy_js_ms^m_0\beta^2+220y^iy_jr_{0m}s^m_0\beta -4y^iy_js_{0|0}\beta\\
\nonumber  \!\!\!\!\!\!&+&\!\!\!\!\!\!\ 12b^iy_jr_{00|0}\beta-7008s_{0|0}^iy_j\beta^2b^4-5120s_{0|0}^iy_j\beta^2b^6-1024s_{0|0}^iy_j\beta^2b^8 +4y^iy_jr_{00}s_0+1728b^is_{j|0}b^2\beta^4\\
\nonumber  \!\!\!\!\!\!&+&\!\!\!\!\!\!\ 5184s_{0|0}^ib_jb^2\beta^3+3456s_{0|0}^ib_j\beta^3b^4 -592g^iy_jr_{00}\beta^2+840\delta_j^ir_{00}r_0\beta^2 -1260\delta_j^ir_{00|0}b^2\beta^2  \\
\nonumber  \!\!\!\!\!\!&-&\!\!\!\!\!\!\ 1584\delta_j^ir_{00|0}b^4\beta^2+24b^iy_jr_{00}^2b^2-1024y^ir_{jk}s^k_0b^6\beta^3 -2144b^ir_jr_{00}\beta^3-708b^ir_{j0}r_{00}\beta^2\\
\nonumber  \!\!\!\!\!\!&+&\!\!\!\!\!\!\ 2144b^ir_{j0}r_0\beta^3 -2736b^ib_js_{0|0}\beta^3+600b^ib_jr_{00}^2\beta+9984y^is_{0|j}b^2\beta^3+7296y^is_{0|j}\beta^3b^4-512r_j^ir_{00}\beta^3b^6\\
\nonumber  \!\!\!\!\!\!&-&\!\!\!\!\!\!\ 6432r_j^ir_{00}\beta^3b^2+24y^iy_jr_{00|0}b^4+12y^iy_jr_{00|0}b^2+576b^is_jr_0\beta^4-9984\delta_j^ir_{0m}s^m_0b^2\beta^3\\
\nonumber  \!\!\!\!\!\!&-&\!\!\!\!\!\!\ 7296\delta_j^ir_{0m}s^m_0b^4\beta^3-4992\delta_j^is_{0|0}b^2\beta^3-3648\delta_j^is_{0|0}b^4\beta^3-576y^ir_{j0}s_0\beta^2-840y^ir_{j0}r_0\beta^2\\
\nonumber  \!\!\!\!\!\!&+&\!\!\!\!\!\!\ 12y^ir_{j0}r_{00}\beta-1152b^ir_js_0\beta^4 +3328y^ir_js_0\beta^3 +840y^ir_jr_{00}\beta^2 +936y^ib_js_ms^m_0\beta^3 +6432b^ir_{j0|0}b^2\beta^3\\
\nonumber  \!\!\!\!\!\!&+&\!\!\!\!\!\!\ 4224b^ir_{j0|0}\beta^3b^4+512b^ir_{j0|0}b^6\beta^3-6432b^ir_{00|j}b^2\beta^3 -4224b^ir_{00|j}\beta^3b^4-512b^ir_{00|j}b^6\beta^3
\\
\nonumber  \!\!\!\!\!\!&-&\!\!\!\!\!\!\ 576g^iy_js_0\beta^3+2520y^ir_{00|j}b^2\beta^2+3168y^ir_{00|j}\beta^2b^4 +1152y^ir_{00|j}\beta^2b^6-2520y^ir_{j0|0}b^2\beta^2\\
\nonumber  \!\!\!\!\!\!&-&\!\!\!\!\!\!\ 3168y^ir_{j0|0}\beta^2b^4-1152y^ir_{j0|0}\beta^2b^6 +4992y^ir_{k0}s^k_jb^2\beta^3 +3648y^ir_{k0}s^k_j\beta^3b^4+512y^ir_{k0}s^k_jb^6\beta^3 \\
\nonumber  \!\!\!\!\!\!&-&\!\!\!\!\!\!\ 9984y^ir_{jk}s^k_0b^2\beta^3-7296y^ir_{jk}s^k_0\beta^3b^4+16y^iy_jr_{00|0}b^6+1728y^is_ms^m_0b^2\beta^4+48y^ib_jr_{00}^2b^4\\
\nonumber  \!\!\!\!\!\!&+&\!\!\!\!\!\!\ 96y^ib_jr_{00}^2b^2-512\delta_j^is_0^2b^4\beta^3 +2048\delta_j^is_0^2b^5\beta^3 -512\delta_j^is_{0|0}b^6\beta^3-1024\delta_j^ir_{0m}s^m_0b^6\beta^3 -576\delta_j^ir_{00|0}b^6\beta^2 \\
\nonumber  \!\!\!\!\!\!&-&\!\!\!\!\!\!\  3456\delta_j^is_ms^m_0b^2\beta^4 -360\delta_j^ir_{00}^2b^4\beta -312\delta_j^ir_{00}^2b^2\beta -228\delta_j^ir_{00}s_0\beta^2 -704\delta_j^is_0^2b^2\beta^3  +6656\delta_j^is_0^2b\beta^3
\\
\nonumber  \!\!\!\!\!\!&+&\!\!\!\!\!\!\ 9728\delta_j^is_0^2b^3\beta^3+3328\delta_j^ir_0s_0\beta^3-24y^ib_jr_{00|0}\beta-3456b^is_{0|j}b^2\beta^4+3456b^ir_{jk}s^k_0b^2\beta^4-1664y^is_jr_0\beta^3\\
\nonumber  \!\!\!\!\!\!&-&\!\!\!\!\!\!\ 4224r_j^ir_{00}\beta^3b^4-3456r_j^is_0b^2\beta^4-1336b^iy_js_0^2\beta^2 +2088b^iy_js_ms^m_0\beta^3-232b^iy_jr_{0m}s^m_0\beta^2 \\
\nonumber  \!\!\!\!\!\!&+&\!\!\!\!\!\!\ 952b^iy_js_{0|0}\beta^2  -1392y^ib_js_0^2\beta^2 +1728g^ir_{j0}b^2\beta^4 +576g^ib_jr_{00}\beta^3-1728b^ir_{k0}s^k_jb^2\beta^4 -4992y^is_{j|0}b^2\beta^3 \\
\nonumber  \!\!\!\!\!\!&-&\!\!\!\!\!\!\ 1328y^iy_jr_0s_0b^2\beta+128y^iy_jr_0s_0b\beta +128y^iy_jr_0s_0b^3\beta -256y^iy_jr_0s_0b^5\beta -1088y^iy_jr_0s_0b^4\beta \\
\nonumber  \!\!\!\!\!\!&-&\!\!\!\!\!\!\  24b^iy_jr_{00}s_0b^2\beta -192b^iy_jr_{00}s_0b\beta -384b^iy_jr_{00}s_0b^3\beta -48b^iy_jr_{00}r_0b^2\beta+768y^ib_jr_{00}r_0b^4\beta \\
\nonumber  \!\!\!\!\!\!&+&\!\!\!\!\!\!\ 1536y^ib_jr_0s_0b^4\beta^2-480y^ib_jr_{00}s_0b\beta +384y^ib_jr_{00}s_0b^3\beta +912y^ib_jr_{00}r_0b^2\beta-480y^ib_jr_{00}r_0b\beta
\\
\nonumber  \!\!\!\!\!\!&+&\!\!\!\!\!\!\ 384y^ib_jr_{00}r_0b^3\beta+1704y^ib_jr_{00}s_0b^2\beta +6912b^ib_jr_0s_0b\beta^3 +4800y^ib_jr_0s_0b^2\beta^2-7104y^ib_jr_0s_0b\beta^2\\
\nonumber  \!\!\!\!\!\!&-&\!\!\!\!\!\!\ 6144y^ib_jr_0s_0b^3\beta^2+1536y^ib_jr_{00}s_0b^5\beta +96y^ib_jr_{00}s_0b^4\beta +1536y^ib_jr_{00}r_0b^5\beta+4608b^ib_jr_{00}r_0b^3\beta^2\\
\nonumber  \!\!\!\!\!\!&+&\!\!\!\!\!\!\ 576b^ib_jr_{00}r_0b^2\beta^2+5760b^ib_jr_{00}r_0b\beta^2 +5760b^ib_jr_{00}s_0b\beta^2+4608b^ib_jr_{00}s_0b^3\beta^2-2048b^iy_jr_0s_0b^5\beta^2\\
\nonumber  \!\!\!\!\!\!&+&\!\!\!\!\!\!\ 512b^iy_jr_0s_0b^4\beta^2+1728b^ib_jr_{00}s_0b^2\beta^2-9728b^iy_jr_0s_0b^3\beta^2 -288b^iy_jr_{00}s_0b^4\beta-192b^iy_jr_{00}r_0b\beta\\
\nonumber  \!\!\!\!\!\!&+&\!\!\!\!\!\!\ -384b^iy_jr_{00}r_0b^3\beta +1280b^iy_jr_0s_0b^2\beta^2 -6080b^iy_jr_0s_0b\beta^2-3456s^i_ks^k_0y_j\beta^3b^4 -1944s^i_0b_js_0\beta^3\\
\nonumber  \!\!\!\!\!\!&-&\!\!\!\!\!\!\ 4320s^i_ks^k_0y_jb^2\beta^3+2592\beta^2s^i_0y_js_0b^4+432\beta s^i_0y_jr_{00}b^4 +288\beta s^i_0y_jr_{00}b^6 -3032b^is_jr_{00}\beta^3\\
\nonumber  \!\!\!\!\!\!&-&\!\!\!\!\!\!\ 864b^is_js_0\beta^4 -4536s^i_0b_jr_{00}b^2\beta^2-2592s^i_0b_jr_{00}\beta^2b^4 +48r^i_0y_jr_{00}b^2\beta+3216r^i_0y_js_0b^2\beta^2\\
\nonumber  \!\!\!\!\!\!&+&\!\!\!\!\!\!\ 48r^i_0y_jr_{00}\beta b^4+2880r^i_0y_js_0\beta^2b^4 +512r^i_0y_js_0b^6\beta^2 -1440r^i_0b_jr_{00}b^2\beta^2-576r^i_0b_jr_{00}\beta^2b^4
\end{eqnarray}
\begin{eqnarray}
\nonumber  \!\!\!\!\!\!&-&\!\!\!\!\!\!\ 2880r^i_0b_js_0b^2\beta^3+2592\beta^2s^i_0y_js_0b^2+216\beta s^i_0y_jr_{00}b^2 +216y^is_jr_{00}b^2\beta -2880y^is_js_0b^2\beta^2
\\
\nonumber  \!\!\!\!\!\!&+&\!\!\!\!\!\!\ 864y^is_jr_{00}\beta b^4+1152b^is_jr_{00}\beta^2b^4 +1728b^is_js_0b^2\beta^3 +3312b^is_jr_{00}b^2\beta^2 -576b^ib_jr_{0m}s^m_0b^2\beta^3 \\
\nonumber  \!\!\!\!\!\!&-&\!\!\!\!\!\!\  2208g^iy_jr_{00}b^2\beta^2 -512g^iy_jr_{00}b^6\beta^2 +512b^ir_{j0}s_0\beta^3b^4 -2048b^ir_{j0}s_0b^5\beta^3-2496\delta_j^ir_{00}s_0b^4\beta^2
\\
\nonumber  \!\!\!\!\!\!&+&\!\!\!\!\!\!\ 2304\delta_j^ir_{00}r_0b^5\beta^2+576y^is_jr_{00}\beta b^6+576y^ib_js_ms^m_0b^2\beta^3-576b^ib_jr_{00|0}b^4\beta^2 +4864\delta_j^ir_0s_0b^2\beta^3\\
\nonumber  \!\!\!\!\!\!&+&\!\!\!\!\!\!\ 6656\delta_j^ir_0s_0b\beta^3 +9728\delta_j^ir_0s_0b^3\beta^3 +2304\delta_j^ir_{00}s_0b^5\beta^2 +2048\delta_j^ir_0s_0b^5\beta^3 +48y^iy_js_{0|0}b^4\beta
\\
\nonumber  \!\!\!\!\!\!&+&\!\!\!\!\!\!\ 512b^iy_js_{0|0}\beta^2b^6-13312y^ir_js_0b\beta^3 -256b^iy_js_0^2b^2\beta^2-6080b^iy_js_0^2b\beta^2 -9728b^iy_js_0^2b^3\beta^2\\
\nonumber  \!\!\!\!\!\!&-&\!\!\!\!\!\!\ 1200b^iy_jr_{0m}s^m_0b^2\beta^2-1728b^iy_jr_{0m}s^m_0b^4\beta^2+3216b^iy_js_{0|0}b^2\beta^2+2880b^iy_js_{0|0}b^4\beta^2+120b^ib_jr_{00}^2b^2\beta \\
\nonumber  \!\!\!\!\!\!&+&\!\!\!\!\!\!\  6480b^ib_jr_{00}s_0\beta^2 +512b^ir_{j0}r_0b^4\beta^3 -2048b^ir_{j0}r_0b^5\beta^3+8576b^is_jr_{00}b\beta^3 +11264b^is_jr_{00}b^3\beta^3\\
\nonumber  \!\!\!\!\!\!&+&\!\!\!\!\!\!\  2048b^is_jr_{00}b^5\beta^3 +4864y^ir_js_0b^2\beta^3 +960y^ib_js_0^2b^2\beta^2 -7104y^ib_js_0^2b\beta^2+720b^ib_jr_{00}r_0\beta^2\\
\nonumber  \!\!\!\!\!\!&-&\!\!\!\!\!\!\ 1440b^ib_jr_{00|0}b^2\beta^2 -2304g^iy_jr_{00}\beta^2b^4 +2048y^is_jr_0b^5\beta^3 +64y^iy_js_{0|0}\beta b^6+1088y^iy_jr_{0m}s^m_0b^6\beta\\
\nonumber  \!\!\!\!\!\!&+&\!\!\!\!\!\!\  1104y^ib_js_{0|0}b^2\beta^2+384y^ib_js_{0|0}b^4\beta^2 -2048y^is_js_0b^5\beta^3 -512y^is_jr_0b^4\beta^3+228y^ib_jr_{00}r_0\beta\\
\nonumber  \!\!\!\!\!\!&-&\!\!\!\!\!\!\ 264y^ib_jr_{00|0}b^2\beta -11264b^ir_{j0}r_0b^3\beta^3 -24b^iy_jr_{00}r_0\beta +1152y^iy_jr_{0m}s^m_0b^2\beta +1968y^iy_jr_{0m}s^m_0b^4\beta \\
\nonumber  \!\!\!\!\!\!&+&\!\!\!\!\!\!\  184y^iy_jr_{00}s_0b^4-64y^iy_jr_{00}s_0b^5 -32y^iy_jr_{00}r_0b^4 -64y^iy_jr_{00}r_0b^5 +1464y^iy_js_ms^m_0b^2\beta^2 \\
\nonumber  \!\!\!\!\!\!&+&\!\!\!\!\!\!\  1920y^iy_js_ms^m_0b^4\beta^2 -256y^iy_js_0^2b^5\beta-32y^iy_js_0^2b^4\beta -624y^ib_jr_{00|0}b^4\beta-8576b^ir_{j0}s_0b\beta^3 \\
\nonumber  \!\!\!\!\!\!&-&\!\!\!\!\!\!\  768b^ir_{j0}r_{00}b^2\beta^2 -1728y^ir_{j0}s_0\beta^2b^4+4608y^ir_{j0}s_0b^5\beta^2 -1152y^ir_{j0}r_0b^4\beta^2+4608y^ir_{j0}r_0b^5\beta^2\\
\nonumber  \!\!\!\!\!\!&+&\!\!\!\!\!\!\ +2112\delta_j^ir_{00}r_0b^2\beta^2+1680\delta_j^ir_{00}r_0b\beta^2 -488y^iy_js_0^2b^2\beta+128y^iy_js_0^2b\beta +128y^iy_js_0^2b^3\beta -392y^iy_jr_0s_0\beta \\
\nonumber  \!\!\!\!\!\!&-&\!\!\!\!\!\!\  4608y^is_jr_{00}b^5\beta^2-8576b^ir_{j0}r_0b\beta^3 +1856b^ir_{j0}s_0b^2\beta^3 -11264b^ir_{j0}s_0b^3\beta^3 +2816b^ir_{j0}r_0b^2\beta^3 \\
\nonumber  \!\!\!\!\!\!&+&\!\!\!\!\!\!\ 4224\delta_j^ir_{00}r_0b^3\beta^2 +512y^iy_js_ms^m_0b^6\beta^2 +8448y^ir_{j0}s_0b^3\beta^2 +8448y^ir_{j0}r_0b^3\beta^2 -19456y^ir_js_0b^3\beta^3\\
\nonumber  \!\!\!\!\!\!&+&\!\!\!\!\!\!\  1024y^ir_js_0\beta^3b^4-4096y^ir_js_0b^5\beta^3 -2784\delta_j^ir_{00}s_0b^2\beta^2 +1680\delta_j^ir_{00}s_0b\beta^2+4224\delta_j^ir_{00}s_0b^3\beta^2\\
\nonumber  \!\!\!\!\!\!&-&\!\!\!\!\!\!\ 3360y^is_jr_{00}b\beta^2-2880b^ib_js_{0|0}b^2\beta^3 -1536y^ib_jr_{0m}s^m_0b^6\beta^2 -384y^ib_jr_{00|0}b^6\beta -5424y^ib_jr_{0m}s^m_0b^2\beta^2
\\
\nonumber  \!\!\!\!\!\!&-&\!\!\!\!\!\!\ 6720y^ib_jr_{0m}s^m_0b^4\beta^2-16y^iy_jr_{00}s_0b +100y^iy_jr_{00}s_0b^2 -64y^iy_jr_{00}s_0b^3-16y^iy_jr_{00}r_0b
\\
\nonumber  \!\!\!\!\!\!&-&\!\!\!\!\!\!\ 32y^iy_jr_{00}r_0b^2-64y^iy_jr_{00}r_0b^3 -6656y^is_js_0b\beta^3+6656y^is_jr_0b\beta^3 -9728y^is_js_0b^3\beta^3 -12864\overline{R}^i_j\beta^3b^4\\
\nonumber  \!\!\!\!\!\!&-&\!\!\!\!\!\!\ 2432y^is_jr_0b^2\beta^3 +9728y^is_jr_0b^3\beta^3-144y^ib_jr_{00}s_0\beta -6144y^ib_js_0^2b^3\beta^2 +2160y^ib_jr_0s_0\beta^2 \\
\nonumber  \!\!\!\!\!\!&+&\!\!\!\!\!\!\ 48b^iy_jr_{00|0}b^2\beta+48b^iy_jr_{00|0}b^4\beta  -4608y^ir_jr_{00}b^5\beta^2 +512b^iy_jr_0s_0\beta^2+2880b^iy_js_ms^m_0b^2\beta^3 \\
\nonumber  \!\!\!\!\!\!&+&\!\!\!\!\!\!\ 1152g^ib_jr_{00}\beta^3b^2+8576b^ir_jr_{00}b\beta^3 -1152g^iy_js_0b^2\beta^3 -8448y^is_jr_{00}b^3\beta^2 +6912b^ib_js_0^2b\beta^3 \\
\nonumber  \!\!\!\!\!\!&+&\!\!\!\!\!\!\  2304b^is_js_0b\beta^4 -2304b^is_jr_0b\beta^4-2112y^ir_{j0}r_0\beta^2b^2 +1152\delta_j^ir_{00}r_0b^4\beta^2+1024\delta_j^ir_0s_0b^4\beta^3\\
\nonumber  \!\!\!\!\!\!&+&\!\!\!\!\!\!\ 3360y^ir_{j0}s_0b\beta^2 +3360y^ir_{j0}r_0b\beta^2-108y^ir_{j0}r_{00}b^2\beta- 2448y^ir_{j0}s_0b^2\beta^2-192y^ir_{j0}r_{00}\beta b^4 \\
\nonumber  \!\!\!\!\!\!&-&\!\!\!\!\!\!\  2048b^iy_js_0^2b^5\beta^2 +512b^iy_js_0^2b^4\beta^2 +2048b^ir_jr_{00}b^5\beta^3 +4608b^ir_js_0b\beta^4 +2112y^ir_jr_{00}b^2\beta^2 \\
\nonumber  \!\!\!\!\!\!&-&\!\!\!\!\!\!\  3360y^ir_jr_{00}b\beta^2 -8448y^ir_jr_{00}b^3\beta^2+1152y^ir_jr_{00}\beta^2b^4 -48b^iy_jr_{00}s_0\beta-512b^iy_jr_{0m}s^m_0b^6\beta^2\\
\nonumber  \!\!\!\!\!\!&-&\!\!\!\!\!\!\ 2816b^ir_jr_{00}\beta^3b^2 -512b^ir_jr_{00}\beta^3b^4 +11264b^ir_jr_{00}b^3\beta^3 +1728r^i_0s_jb^2\beta^4 +952r^i_0y_js_0\beta^2 \\
\nonumber  \!\!\!\!\!\!&+&\!\!\!\!\!\!\  12r^i_0y_jr_{00}\beta-2736r^i_0b_js_0\beta^3-684r^i_0b_jr_{00}\beta^2+36\beta s^i_0y_jr_{00} +648\beta^2s^i_0y_js_0 +6432r^i_0r_{j0}b^2\beta^3 \\
\nonumber  \!\!\!\!\!\!&+&\!\!\!\!\!\!\ 4224r^i_0r_{j0}\beta^3b^4+512r^i_0r_{j0}b^6\beta^3 -1152y^is_js_0\beta^2b^4 -1296s^i_ks^k_0y_j\beta^3 +1728r^i_0s_j\beta^4 +648s^i_ks^k_0b_j\beta^4 \\
\nonumber  \!\!\!\!\!\!&-&\!\!\!\!\!\!\  3230560s^i_0s_{0j}b^2\beta^3-4257984s^i_0s_{0j}\beta^3b^4 -2336\overline{R}^i_j\beta^3 -9760\overline{R}^i_j\beta^3b^2   -5632\overline{R}^i_j\beta^3b^6 -512\overline{R}^i_j\beta^3b^8
\\
\nonumber t_9 \!\!\!\!&:=&\!\!\!\!\ -432s^is_j\beta^4 -1188s^i_ks^k_j\beta^4 +182712s^i_0s_{0j}\beta^2 +3\delta_j^ir_{00}^2 +1488s_{j|0}^i\beta^3 -2976s_{0|j}^i\beta^3 +396s^i_0b_jr_{00}\beta \\
\nonumber  \!\!\!\!\!\!&+&\!\!\!\!\!\!\ 254208s^i_0s_{0j}\beta^2b^8 -744r^i_0r_{j0}\beta^2 +744r_j^ir_{00}\beta^2 +2088\delta_j^is_ms^m_0\beta^3 +5472s_{j|0}^ib^2\beta^3-476y^ir_{k0}s^k_j\beta^2\\
\nonumber  \!\!\!\!\!\!&+&\!\!\!\!\!\!\ 744b^ir_{00|j}\beta^2-624s_{0|0}^ib_j\beta^2+432b^is_ms^m_0\beta^4 +1368b^ir_{k0}s^k_j\beta^3+38\delta_j^ir_{00|0}\beta +952y^ir_{jk}s^k_0\beta^2\\
\nonumber  \!\!\!\!\!\!&+&\!\!\!\!\!\!\ 476y^is_{j|0}\beta^2+80s_{0|0}^iy_j\beta+36\delta_j^ir_{00}^2b^4 +24\delta_j^ir_{00}^2b^2 -952y^is_{0|j}\beta^2 +2736b^is_{0|j}\beta^3 -1044y^is_ms^m_0\beta^3\\
\nonumber  \!\!\!\!\!\!&-&\!\!\!\!\!\!\ 10944s_{0|j}^ib^2\beta^3-12y^iy_jr_{0m}s^m_0 +12y^iy_js_0^2-72b^ib_jr_{00}^2 +2y^ib_jr_{00|0}-1368b^is_{j|0}\beta^3+2736r_j^is_0\beta^3
\end{eqnarray}
\begin{eqnarray}
\nonumber  \!\!\!\!\!\!&-&\!\!\!\!\!\!\ 76y^ir_{00|j}\beta+76y^ir_{j0|0}\beta-3072s_{0|j}^ib^6\beta^3-744b^ir_{j0|0}\beta^2 -452\delta_j^is_0^2\beta^2 +5760s_{j|0}^i\beta^3b^4+1536s_{j|0}^ib^6\beta^3
\\
\nonumber  \!\!\!\!\!\!&-&\!\!\!\!\!\!\  2736b^ir_{jk}s^k_0\beta^3-1368s^ir_{j0}\beta^3-11520s_{0|j}^i\beta^3b^4+952\delta_j^ir_{0m}s^m_0\beta^2+476\delta_j^is_{0|0}\beta^2+228y^is_js_0\beta\\
\nonumber  \!\!\!\!\!\!&+&\!\!\!\!\!\!\ 144b^ib_jr_{00|0}\beta-3120b^ib_js_0^2\beta^2 +1296b^ib_js_ms^m_0\beta^3 +1056b^ib_jr_{0m}s^m_0\beta^2+1608b^ib_js_{0|0}\beta^2 \\
\nonumber  \!\!\!\!\!\!&-&\!\!\!\!\!\!\ 144y^is_jr_{00}b^6+200y^ib_jr_{0m}s^m_0\beta-1280y^is_{0|j}\beta^2b^6+24b^ir_{j0}s_0\beta^2-80y^ib_js_{0|0}\beta-2736b^ir_{j0|0}b^2\beta^2 \\
\nonumber  \!\!\!\!\!\!&+&\!\!\!\!\!\!\  24y^iy_jr_0s_0-20y^ib_jr_{00}s_0-60y^iy_js_ms^m_0\beta+1408s_{0|0}^iy_j\beta b^6+512s_{0|0}^iy_j\beta b^8+24y^ir_{j0}r_{00}b^4 \\
\nonumber  \!\!\!\!\!\!&-&\!\!\!\!\!\!\ 2880b^is_{j|0}b^2\beta^3-1152b^is_{j|0}\beta^3b^4+5760b^is_{0|j}b^2\beta^3+1632b^is_js_0\beta^3+1380b^is_jr_{00}\beta^2-2880s_{0|0}^ib_jb^2\beta^2 \\
\nonumber  \!\!\!\!\!\!&-&\!\!\!\!\!\!\ 4032s_{0|0}^ib_j\beta^2b^4-1536s_{0|0}^ib_j\beta^2b^6-36b^ib_jr_{00}^2b^2 +312s^iy_js_0\beta^2 +108s^iy_jr_{00}\beta-128\delta_j^ir_{00}r_0\beta \\
\nonumber  \!\!\!\!\!\!&+&\!\!\!\!\!\!\ 192\delta_j^ir_{00|0}b^2\beta+312\delta_j^ir_{00|0}b^4\beta +352\delta_j^is_0^2b^2\beta^2 +12y^ir_{j0}r_{00}b^2+1848y^is_{j|0}b^2\beta^2 +2304r_j^is_0\beta^3b^4 \\
\nonumber  \!\!\!\!\!\!&+&\!\!\!\!\!\!\ 1920b^ir_js_0\beta^3-20y^ib_jr_{00}r_0 -912b^ir_{j0}r_0\beta^2+144b^ir_{j0}r_{00}\beta -3696y^is_{0|j}b^2\beta^2-4224y^is_{0|j}\beta^2b^4 \\
\nonumber  \!\!\!\!\!\!&-&\!\!\!\!\!\!\  276s^ib_jr_{00}\beta^2 +2736r_j^ir_{00}b^2\beta^2+768r_j^ir_{00}b^6\beta^2 +2880r_j^ir_{00}\beta^2b^4-24y^iy_js_0^2b^2-144y^iy_jr_{0m}s^m_0b^4 \\
\nonumber  \!\!\!\!\!\!&-&\!\!\!\!\!\!\ 72y^iy_jr_{0m}s^m_0b^2-36b^iy_jr_{00}s_0 -960b^is_jr_0\beta^3-2448y^is_ms^m_0b^2\beta^3-1232\delta_j^ir_0s_0\beta^2+4896\delta_j^is_ms^m_0b^2\beta^3 \\
\nonumber  \!\!\!\!\!\!&+&\!\!\!\!\!\!\ 3696\delta_j^ir_{0m}s^m_0b^2\beta^2 +4224\delta_j^ir_{0m}s^m_0b^4\beta^2 +1848\delta_j^is_{0|0}b^2\beta^2+2112\delta_j^is_{0|0}b^4\beta^2+104y^ir_{j0}s_0\beta \\
\nonumber  \!\!\!\!\!\!&+&\!\!\!\!\!\!\  128y^ir_{j0}r_0\beta+912b^ir_jr_{00}\beta^2 -1232y^ir_js_0\beta^2-128y^ir_jr_{00}\beta +72y^ib_jr_{00|0}b^4+24y^ib_jr_{00|0}b^2\\
\nonumber  \!\!\!\!\!\!&-&\!\!\!\!\!\!\ 5760b^ir_{jk}s^k_0b^2\beta^3-2880b^ir_{j0|0}\beta^2b^4-768b^ir_{j0|0}\beta^2b^6+2736b^ir_{00|j}b^2\beta^2+2880b^ir_{00|j}\beta^2b^4 \\
\nonumber  \!\!\!\!\!\!&+&\!\!\!\!\!\!\ 768b^ir_{00|j}\beta^2b^6-832y^is_js_0\beta^2-140y^is_jr_{00}\beta-384y^ir_{00|j}b^2\beta-624y^ir_{00|j}\beta b^4-320y^ir_{00|j}\beta b^6\\
\nonumber  \!\!\!\!\!\!&+&\!\!\!\!\!\!\ 384y^ir_{j0|0}b^2\beta+624y^ir_{j0|0}\beta b^4+320y^ir_{j0|0}\beta b^6-1848y^ir_{k0}s^k_jb^2\beta^2-2112y^ir_{k0}s^k_j\beta^2b^4\\
\nonumber  \!\!\!\!\!\!&-&\!\!\!\!\!\!\ 640y^ir_{k0}s^k_j\beta^2b^6+3696y^ir_{jk}s^k_0b^2\beta^2 +4224y^ir_{jk}s^k_0\beta^2b^4 +1280y^ir_{jk}s^k_0\beta^2b^6 -96y^iy_js_0^2b^4 \\
\nonumber  \!\!\!\!\!\!&-&\!\!\!\!\!\!\  96y^iy_jr_{0m}s^m_0b^6-1152y^is_ms^m_0\beta^3b^4+640\delta_j^is_{0|0}\beta^2b^6+1280\delta_j^ir_{0m}s^m_0b^6\beta^2+160\delta_j^ir_{00|0}b^6\beta \\
\nonumber  \!\!\!\!\!\!&-&\!\!\!\!\!\!\  2560\delta_j^is_0^2b^5\beta^2+640\delta_j^is_0^2b^4\beta^2 +2304\delta_j^is_ms^m_0b^4\beta^3 +40\delta_j^ir_{00}s_0\beta-2464\delta_j^is_0^2b\beta^2 -5632\delta_j^is_0^2b^3\beta^2 \\
\nonumber  \!\!\!\!\!\!&+&\!\!\!\!\!\!\ 188y^ib_js_0^2\beta+64y^ib_jr_{00|0}b^6+2304b^is_{0|j}\beta^3b^4 -2304b^ir_{jk}s^k_0\beta^3b^4 +616y^is_jr_0\beta^2 -2880s^ir_{j0}b^2\beta^3
\\
\nonumber  \!\!\!\!\!\!&+&\!\!\!\!\!\!\ 5760r_j^is_0b^2\beta^3 +144b^iy_js_0^2\beta -888b^iy_js_ms^m_0\beta^2 +72b^iy_jr_{0m}s^m_0\beta-144b^iy_js_{0|0}\beta -336y^ib_js_ms^m_0\beta^2 \\
\nonumber  \!\!\!\!\!\!&-&\!\!\!\!\!\!\  1152s^ir_{j0}\beta^3b^4 +2880b^ir_{k0}s^k_jb^2\beta^3+1152b^ir_{k0}s^k_j\beta^3b^4 +2112y^is_{j|0}\beta^2b^4+640y^is_{j|0}\beta^2b^6\\
\nonumber  \!\!\!\!\!\!&+&\!\!\!\!\!\!\ 544s_{0|0}^iy_jb^2\beta+1344s_{0|0}^iy_j\beta b^4+552\overline{R}^i_j\beta^2 +2976\overline{R}^i_j\beta^2b^2+5472\overline{R}^i_j\beta^2b^4+3840\overline{R}^i_j\beta^2b^6\\
\nonumber  \!\!\!\!\!\!&+&\!\!\!\!\!\!\ 768\overline{R}^i_j\beta^2b^8+3552b^ib_js_{0|0}b^2\beta^2+1536b^ib_js_{0|0}b^4\beta^2+128y^ib_js_{0|0}\beta b^6 +3072b^iy_jr_0s_0b^3\beta \\
\nonumber  \!\!\!\!\!\!&-&\!\!\!\!\!\!\ 512y^ib_jr_0s_0b^5\beta   -1408y^ib_jr_0s_0b^4\beta -576b^ib_jr_{00}s_0b^4\beta +384b^ib_jr_0s_0b^2\beta^2 -8832b^ib_jr_0s_0b\beta^2 \\
\nonumber  \!\!\!\!\!\!&-&\!\!\!\!\!\!\  6144b^ib_jr_0s_0b^3\beta^2 -1984y^ib_jr_0s_0b^2\beta +1600y^ib_jr_0s_0b\beta +1792y^ib_jr_0s_0b^3\beta-2304b^ib_jr_{00}r_0b^3\beta \\
\nonumber  \!\!\!\!\!\!&-&\!\!\!\!\!\!\ 288b^ib_jr_{00}r_0b^2\beta -1728b^ib_jr_{00}r_0b\beta -1728b^ib_jr_{00}s_0b\beta -2304b^ib_jr_{00}s_0b^3\beta -384b^iy_jr_0s_0b^4\beta \\
\nonumber  \!\!\!\!\!\!&-&\!\!\!\!\!\!\  2160b^ib_jr_{00}s_0b^2\beta+1536b^iy_jr_0s_0b^5\beta-480b^iy_jr_0s_0b^2\beta +1152b^iy_jr_0s_0b\beta+1536s^i_ks^k_0y_j\beta^2b^6\\
\nonumber  \!\!\!\!\!\!&+&\!\!\!\!\!\!\ 1944s^i_0b_js_0\beta^2+2016s^i_ks^k_0y_jb^2\beta^2 +3168s^i_ks^k_0y_j\beta^2b^4-864\beta s^i_0y_js_0b^4-576\beta s^i_0y_js_0b^6\\
\nonumber  \!\!\!\!\!\!&-&\!\!\!\!\!\!\ 432\beta s^i_0y_js_0b^2  1728s^i_0b_jr_{00}b^2\beta +3888s^i_0b_js_0b^2\beta^2 +2160s^i_0b_jr_{00}\beta b^4+576s^i_0b_jr_{00}\beta b^6 \\
\nonumber  \!\!\!\!\!\!&-&\!\!\!\!\!\!\ 672r^i_0y_js_0b^2\beta -960r^i_0y_js_0\beta b^4 -384r^i_0y_js_0b^6\beta +432r^i_0b_jr_{00}b^2\beta +3552r^i_0b_js_0b^2\beta^2 +288r^i_0b_jr_{00}\beta b^4 \\
\nonumber  \!\!\!\!\!\!&+&\!\!\!\!\!\!\ 1536r^i_0b_js_0\beta^2b^4  +1408y^ib_jr_{0m}s^m_0b^6\beta +2784y^ib_jr_{0m}s^m_0b^4\beta -144y^ib_js_{0|0}b^2\beta +96y^iy_jr_0s_0b^2 \\
\nonumber  \!\!\!\!\!\!&-&\!\!\!\!\!\!\ 332y^is_jr_{00}b^2\beta -1312y^is_js_0b^2\beta^2 -176y^is_jr_{00}\beta b^4  +1408y^is_jr_0\beta^2b^2 +5632y^is_js_0b^3\beta^2\\
\nonumber  \!\!\!\!\!\!&-&\!\!\!\!\!\!\ 5632y^is_jr_0b^3\beta^2 +1792y^ib_js_0^2b^3\beta -568y^ib_jr_0s_0\beta-144b^iy_jr_0s_0\beta-2544b^iy_js_ms^m_0b^2\beta^2\\
\nonumber  \!\!\!\!\!\!&-&\!\!\!\!\!\!\  1536b^iy_js_ms^m_0b^4\beta^2 +768b^is_jr_{00}\beta^2b^4 +960b^is_js_0b^2\beta^3 -1248s^ib_jr_{00}b^2\beta^2 -960s^ib_jr_{00}\beta^2b^4 \\
 \nonumber  \!\!\!\!\!\!&+&\!\!\!\!\!\!\ 1920b^ir_jr_{00}b^2\beta^2 -3648b^ir_jr_{00}b\beta^2 -7680b^ir_jr_{00}b^3\beta^2 +960s^iy_js_0\beta^2b^4+2352b^is_jr_{00}b^2\beta^2\\
 \nonumber  \!\!\!\!\!\!&-&\!\!\!\!\!\!\ 6144b^ib_js_0^2b^3\beta^2 -96b^ib_jr_0s_0\beta^2 -3072b^is_js_0b^3\beta^3 -768b^is_jr_0b^2\beta^3+3072b^is_jr_0b^3\beta^3 -3840b^is_js_0b\beta^3 \\
 \nonumber  \!\!\!\!\!\!&+&\!\!\!\!\!\!\  3840b^is_jr_0b\beta^3-320\delta_j^ir_{00}r_0b^4\beta -2560\delta_j^ir_0s_0b^5\beta^2 -1280\delta_j^ir_0s_0b^4\beta^2+584y^ir_{j0}s_0b^2\beta-512y^ir_{j0}s_0b\beta
 \end{eqnarray}
\begin{eqnarray}
 \nonumber  \!\!\!\!\!\!&-&\!\!\!\!\!\!\ 1664y^ir_{j0}s_0b^3\beta+416y^ir_{j0}r_0b^2\beta -512y^ir_{j0}r_0b\beta -1664y^ir_{j0}r_0b^3\beta +608y^ir_{j0}s_0\beta b^4+192b^iy_js_0^2b^4\beta
\\
 \nonumber  \!\!\!\!\!\!&-&\!\!\!\!\!\!\  7680b^ir_js_0b\beta^3 +1536b^ir_js_0b^2\beta^3 -6144b^ir_js_0b^3\beta^3 +512y^ir_jr_{00}b\beta +1664y^ir_jr_{00}b^3\beta +4928y^ir_js_0b\beta^2 \\
 \nonumber  \!\!\!\!\!\!&-&\!\!\!\!\!\!\  416y^ir_jr_{00}b^2\beta-2816y^ir_js_0b^2\beta^2-320y^ir_jr_{00}\beta b^4+1280y^ir_jr_{00}b^5\beta-384b^iy_js_{0|0}\beta b^6
\\
 \nonumber  \!\!\!\!\!\!&+&\!\!\!\!\!\!\  384b^iy_jr_{0m}s^m_0b^6\beta +1536b^iy_js_0^2b^5\beta +768b^ir_jr_{00}\beta^2b^4 -2880r^i_0s_jb^2\beta^3-1152r^i_0s_j\beta^3b^4\\
 \nonumber  \!\!\!\!\!\!&-&\!\!\!\!\!\!\ 144r^i_0y_js_0\beta+1608r^i_0b_js_0\beta^2+144r^i_0b_jr_{00}\beta-72\beta s^i_0y_js_0 -1728s^i_ks^k_0b_jb^2\beta^3-2736r^i_0r_{j0}b^2\beta^2 \\
 \nonumber  \!\!\!\!\!\!&-&\!\!\!\!\!\!\ 2880r^i_0r_{j0}\beta^2b^4 -768r^i_0r_{j0}\beta^2b^6 +432y^is_js_0b^2\beta -144y^is_js_0\beta b^4-576b^is_js_0\beta^2b^4-1152b^is_jr_{00}b^2\beta\\
 \nonumber  \!\!\!\!\!\!&-&\!\!\!\!\!\!\ 4464b^is_js_0b^2\beta^2-720b^is_jr_{00}\beta b^4 +1056b^ib_jr_{0m}s^m_0b^2\beta^2 -448y^is_js_0\beta^2b^4-192y^is_js_0b^6\beta\\
 \nonumber  \!\!\!\!\!\!&+&\!\!\!\!\!\!\ +528s^iy_jr_{00}b^2\beta -768b^ir_{j0}s_0\beta^2b^4 +3072b^ir_{j0}s_0b^5\beta^2 -768b^ir_{j0}r_0b^4\beta^2-224y^ib_jr_{00}s_0b^4 \\
 \nonumber  \!\!\!\!\!\!&-&\!\!\!\!\!\!\ 256y^ib_jr_{00}s_0b^5 -128y^ib_jr_{00}r_0b^4 -256y^ib_jr_{00}r_0b^5 +32y^ib_jr_{00}s_0b 368y^ib_jr_{00}s_0b^2 -64y^ib_jr_{00}s_0b^3 \\
 \nonumber  \!\!\!\!\!\!&+&\!\!\!\!\!\!\ 32y^ib_jr_{00}r_0b -104y^ib_jr_{00}r_0b^2 -64y^ib_jr_{00}r_0b^3 -640\delta_j^ir_{00}r_0b^5\beta +1280y^is_jr_{00}b^5\beta -48y^ib_js_ms^m_0b^2\beta^2
\\
 \nonumber  \!\!\!\!\!\!&+&\!\!\!\!\!\!\  384y^ib_js_ms^m_0b^4\beta^2 +1440y^ib_jr_{0m}s^m_0b^2\beta +960b^ib_js_0^2b^2\beta^2 -8832b^ib_js_0^2b\beta^2 -2816\delta_j^ir_0s_0b^2\beta^2 \\
 \nonumber  \!\!\!\!\!\!&-&\!\!\!\!\!\!\  2464\delta_j^ir_0s_0b\beta^2 -5632\delta_j^ir_0s_0b^3\beta^2 -640\delta_j^ir_{00}s_0b^5\beta +832\delta_j^ir_{00}s_0b^4\beta +11264y^ir_js_0b^3\beta^2  +96b^iy_js_0^2b^2\beta\\
 \nonumber  \!\!\!\!\!\!&+&\!\!\!\!\!\!\  1152b^iy_js_0^2b\beta +3072b^iy_js_0^2b^3\beta +384b^iy_jr_{0m}s^m_0b^2\beta +672b^iy_jr_{0m}s^m_0b^4\beta-672b^iy_js_{0|0}b^2\beta \\
 \nonumber  \!\!\!\!\!\!&-&\!\!\!\!\!\!\ 960b^iy_js_{0|0}b^4\beta -2808b^ib_jr_{00}s_0\beta+3072b^ir_{j0}r_0b^5\beta^2 -3648b^is_jr_{00}b\beta^2-7680b^is_jr_{00}b^3\beta^2\\
 \nonumber  \!\!\!\!\!\!&-&\!\!\!\!\!\!\ 3072b^is_jr_{00}b^5\beta^2 -1312y^ib_js_0^2b^2\beta +1600y^ib_js_0^2b\beta -216b^ib_jr_{00}r_0\beta+432b^ib_jr_{00|0}b^2\beta+288b^ib_jr_{00|0}b^4\beta \\
 \nonumber  \!\!\!\!\!\!&+&\!\!\!\!\!\!\  1104s^iy_js_0b^2\beta^2+816s^iy_jr_{00}\beta b^4+384s^iy_jr_{00}\beta b^6+96y^ib_js_{0|0}b^4\beta+2560y^is_js_0b^5\beta^2 +640y^is_jr_0b^4\beta^2 \\
 \nonumber  \!\!\!\!\!\!&-&\!\!\!\!\!\!\ 2560y^is_jr_0b^5\beta^2 +96y^iy_jr_0s_0b^4-336y^iy_js_ms^m_0b^2\beta -624y^iy_js_ms^m_0b^4\beta-512y^ib_js_0^2b^5\beta -64y^ib_js_0^2b^4\beta\\
 \nonumber  \!\!\!\!\!\!&-&\!\!\!\!\!\!\ 1920b^ir_{j0}r_0\beta^2b^2 +7680b^ir_{j0}s_0b^3\beta^2+3648b^ir_{j0}s_0b\beta^2 +3648b^ir_{j0}r_0b\beta^2+216b^ir_{j0}r_{00}b^2\beta \\
 \nonumber  \!\!\!\!\!\!&-&\!\!\!\!\!\!\ 1056b^ir_{j0}s_0b^2\beta^2-416\delta_j^ir_{00}r_0b^2\beta -256\delta_j^ir_{00}r_0b\beta +7680b^ir_{j0}r_0b^3\beta^2 -3072b^ir_jr_{00}b^5\beta^2
\\
 \nonumber  \!\!\!\!\!\!&-&\!\!\!\!\!\!\ 832\delta_j^ir_{00}r_0b^3\beta -384y^iy_js_ms^m_0b^6\beta-1280y^ir_{j0}s_0b^5\beta +320y^ir_{j0}r_0b^4\beta-1280y^ir_{j0}r_0b^5\beta-72b^iy_jr_{00}s_0b^2 \\
 \nonumber  \!\!\!\!\!\!&-&\!\!\!\!\!\!\  1280y^ir_js_0\beta^2b^4+5120y^ir_js_0b^5\beta^2+640\delta_j^ir_{00}s_0b^2\beta -256\delta_j^ir_{00}s_0b\beta -832\delta_j^ir_{00}s_0b^3\beta+512y^is_jr_{00}b\beta
\\
\nonumber  \!\!\!\!\!\!&+&\!\!\!\!\!\!\  1664y^is_jr_{00}b^3\beta +2464y^is_js_0b\beta^2 -2464y^is_jr_0b\beta^2 -384b^ib_jr_{0m}s^m_0b^4\beta^2 +408s^i_ks^k_0y_j\beta^2 \\
\nonumber  \!\!\!\!\!\!&-&\!\!\!\!\!\!\  1368r^i_0s_j\beta^3-864s^i_ks^k_0b_j\beta^3-1728s^i_ks^k_jb^2\beta^4 +985056s^i_0s_{0j}b^2\beta^2 +1811232s^i_0s_{0j}\beta^2b^4 \\
\nonumber  \!\!\!\!\!\!&+&\!\!\!\!\!\!\ 1271040s^i_0s_{0j}\beta^2b^6 -144y^is_jr_{00}b^4 -36y^is_jr_{00}b^2-2736b^is_js_0\beta^2-396b^is_jr_{00}\beta
\\
\nonumber t_{10} \!\!\!\!&:=&\!\!\!\!\ 648s^is_j\beta^3+864s^i_ks^k_j\beta^3-23832s^i_0s_{0j}\beta-2\delta_j^ir_{00|0}-424s_{j|0}^i\beta^2+848s_{j|0}^i\beta^2+4y^ir_{00|j}-4s_{0|0}^iy_j\\
\nonumber  \!\!\!\!\!\!&-&\!\!\!\!\!\!\ 4y^ir_{j0|0}-36s^i_0b_jr_{00} +120r^i_0r_{j0}\beta-12y^is_js_0 +36b^is_jr_{00} -712\delta_j^is_ms^m_0\beta^2+68y^ir_{k0}s^k_j\beta\\
\nonumber  \!\!\!\!\!\!&-&\!\!\!\!\!\!\ 120b^ir_{00|j}\beta+112s_{0|0}^ib_j\beta-648b^is_ms^m_0\beta^3 -536b^ir_{k0}s^k_j\beta^2 -136y^ir_{jk}s^k_0\beta-68y^is_{j|0}\beta\\
\nonumber  \!\!\!\!\!\!&-&\!\!\!\!\!\!\ 24\delta_j^ir_{00|0}b^4-12\delta_j^ir_{00|0}b^2+136y^is_{0|j}\beta+356y^is_ms^m_0\beta^2-2144s_{j|0}^ib^2\beta^2-64s_{0|0}^iy_jb^8\\
\nonumber  \!\!\!\!\!\!&-&\!\!\!\!\!\!\ 128s_{0|0}^iy_jb^6-4\delta_j^ir_{00}s_0+8\delta_j^ir_{00}r_0-16\delta_j^ir_{00|0}b^6-96s_{0|0}^iy_jb^4-32s_{0|0}^iy_jb^2-8s^iy_jr_{00}+8y^ir_jr_{00}\\
\nonumber  \!\!\!\!\!\!&-&\!\!\!\!\!\!\ 8y^ir_{j0}s_0-8y^ir_{j0}r_0+48y^ir_{00|j}b^4+24y^ir_{00|j}b^2+32y^ir_{00|j}b^6-48y^ir_{j0|0}b^4-24y^ir_{j0|0}b^2\\
\nonumber  \!\!\!\!\!\!&-&\!\!\!\!\!\!\ 32y^ir_{j0|0}b^6+7104s_{j|0}^i\beta^2b^4-12b^ir_{j0}r_{00} +8y^is_jr_{00} +4y^iy_js_ms^m_0-8b^iy_jr_{0m}s^m_0+8b^iy_js_{0|0}
\\
\nonumber  \!\!\!\!\!\!&+&\!\!\!\!\!\!\ 40b^iy_js_0^2-12b^ib_jr_{00|0}-16y^ib_jr_{0m}s^m_0+4y^ib_js_{0|0}+20y^ib_js_0^2+536b^is_{j|0}\beta^2-1072r_j^is_0\beta^2\\
\nonumber  \!\!\!\!\!\!&-&\!\!\!\!\!\!\ 120r_j^ir_{00}\beta+4096s_{j|0}^i\beta^2b^6+512s_{j|0}^i\beta^2b^8+120b^ir_{j0|0}\beta+80\delta_j^is_0^2\beta
-1072b^is_{0|j}\beta^2\\
\nonumber  \!\!\!\!\!\!&-&\!\!\!\!\!\!\ 3552s_{j|0}^i\beta^2b^4-2048s_{j|0}^i\beta^2b^6-256s_{j|0}^i\beta^2b^8+4288s_{j|0}^ib^2\beta^2+1072b^ir_{jk}s^k_0\beta^2+536s^ir_{j0}\beta^2
\\
\nonumber  \!\!\!\!\!\!&-&\!\!\!\!\!\!\ 136\delta_j^ir_{0m}s^m_0\beta-68\delta_j^is_{0|0}\beta+144y^is_js_0b^4+864b^is_js_0\beta +144b^is_jr_{00}b^2+144b^is_jr_{00}b^4-576b^ir_{00|j}b^2\beta\\
\nonumber  \!\!\!\!\!\!&-&\!\!\!\!\!\!\ 864b^ir_{00|j}\beta b^4-384b^ir_{00|j}\beta b^6+136y^is_js_0\beta-80s^iy_js_0\beta +336y^ir_{k0}s^k_jb^2\beta+528y^ir_{k0}s^k_j\beta b^4\\
\nonumber  \!\!\!\!\!\!&+&\!\!\!\!\!\!\ 256y^ir_{k0}s^k_j\beta b^6-672y^ir_{jk}s^k_0b^2\beta-1056y^ir_{jk}s^k_0\beta b^4-512y^ir_{jk}s^k_0\beta b^6-32y^ir_jr_{00}b +32y^ir_jr_{00}b^2 \\
\nonumber  \!\!\!\!\!\!&-&\!\!\!\!\!\!\  128y^ir_jr_{00}b^3 +32y^ir_jr_{00}b^4-128y^ir_jr_{00}b^5+32y^iy_js_ms^m_0b^6 +416y^ib_js_0^2b^2+1248y^is_ms^m_0\beta^2b^4
 \end{eqnarray}
\begin{eqnarray}
\nonumber  \!\!\!\!\!\!&+&\!\!\!\!\!\!\ 256y^is_ms^m_0\beta^2b^6 -256\delta_j^is_{0|0}\beta b^6 -512\delta_j^ir_{0m}s^m_0b^6\beta +1024\delta_j^is_0^2b^5\beta  -352\delta_j^is_0^2b^4\beta -512\delta_j^is_ms^m_0b^6\beta^2 \\
\nonumber  \!\!\!\!\!\!&+&\!\!\!\!\!\!\  448\delta_j^is_0^2b\beta +1408\delta_j^is_0^2b^3\beta +320y^ib_js_0^2b^4+56y^ib_jr_0s_0-512b^is_{0|j}\beta^2b^6-864b^is_ms^m_0b^2\beta^3
\\
\nonumber  \!\!\!\!\!\!&+&\!\!\!\!\!\!\ 1536s^ir_{j0}\beta^2b^4 +512b^ir_{jk}s^k_0\beta^2b^6 +1776s^ir_{j0}b^2\beta^2-1776b^ir_{k0}s^k_jb^2\beta^2-32y^is_jr_{00}b-112y^is_jr_0\beta\\
\nonumber  \!\!\!\!\!\!&-&\!\!\!\!\!\!\ 128y^ib_js_0^2b^3-3072r_j^is_0\beta^2b^4 +152b^iy_js_ms^m_0\beta +40y^ib_js_ms^m_0\beta +256s^ir_{j0}\beta^2b^6-72s^ib_js_0\beta^2\\
\nonumber  \!\!\!\!\!\!&-&\!\!\!\!\!\!\ 1536b^ir_{k0}s^k_j\beta^2b^4-256y^is_{j|0}\beta b^6 +1024b^ib_jr_0s_0b^5\beta -256b^ib_jr_0s_0b^4\beta +3712b^ib_jr_0s_0b\beta
\\
\nonumber  \!\!\!\!\!\!&+&\!\!\!\!\!\!\ 5632b^ib_jr_0s_0b^3\beta -256b^ib_jr_0s_0b^2\beta -256s^i_ks^k_0y_j\beta b^8-648s^i_0b_js_0\beta-416s^i_ks^k_0y_jb^2\beta-960s^i_ks^k_0y_j\beta b^4\\
\nonumber  \!\!\!\!\!\!&-&\!\!\!\!\!\!\ 896s^i_ks^k_0y_j\beta b^6-72\overline{R}^i_j\beta-1152\overline{R}^i_j\beta b^4-1152\overline{R}^i_j\beta b^6-480\overline{R}^i_j\beta b^2-384\overline{R}^i_j\beta b^8-2592s^i_0b_js_0b^2\beta\\
\nonumber  \!\!\!\!\!\!&-&\!\!\!\!\!\!\ 2592s^i_0b_js_0\beta b^4-1440r^i_0b_js_0b^2\beta-1344r^i_0b_js_0\beta b^4-256r^i_0b_js_0b^6\beta+408b^ib_jr_{00}s_0+24b^ib_jr_{00}r_0\\
\nonumber  \!\!\!\!\!\!&-&\!\!\!\!\!\!\ 48b^ib_jr_{00|0}b^4-48b^ib_jr_{00|0}b^2-64y^ib_js_{0|0}b^6-160\delta_j^is_0^2b^2\beta -256b^iy_js_0^2b^5+64b^iy_js_{0|0}b^6+16b^iy_jr_0s_0 \\
\nonumber  \!\!\!\!\!\!&+&\!\!\!\!\!\!\ 64b^iy_js_0^2b^4-64b^iy_jr_{0m}s^m_0b^6-64b^iy_js_0^2b +112b^iy_js_0^2b^2 -256b^iy_js_0^2b^3-96b^iy_jr_{0m}s^m_0b^4\\
\nonumber  \!\!\!\!\!\!&-&\!\!\!\!\!\!\ 48b^iy_jr_{0m}s^m_0b^2+48b^iy_js_{0|0}b^2+96b^iy_js_{0|0}b^4+32y^ir_{j0}s_0b-56y^ir_{j0}s_0b^2 +128y^ir_{j0}s_0b^3\\
\nonumber  \!\!\!\!\!\!&+&\!\!\!\!\!\!\ 32y^ir_{j0}r_0b-32y^ir_{j0}r_0b^2+128y^ir_{j0}r_0b^3-192b^ir_jr_{00}\beta-336y^is_{j|0}b^2\beta-528y^is_{j|0}\beta b^4-512r_j^is_0b^6\beta^2 \\
\nonumber  \!\!\!\!\!\!&-&\!\!\!\!\!\!\ 1184b^ir_js_0\beta^2+192b^ir_{j0}r_0\beta -96s^iy_jr_{00}b^4-64s^iy_jr_{00}b^6 +1056y^is_{0|j}\beta b^4+512y^is_{0|j}\beta b^6
\\
\nonumber  \!\!\!\!\!\!&-&\!\!\!\!\!\!\  576r_j^ir_{00}b^2\beta-3552r_j^is_0b^2\beta^2-864r_j^ir_{00}\beta b^4-384r_j^ir_{00}\beta b^6 +24y^iy_js_ms^m_0b^2 +48y^iy_js_ms^m_0b^4 \\
\nonumber  \!\!\!\!\!\!&+&\!\!\!\!\!\!\  592b^is_jr_0\beta^2+1272y^is_ms^m_0b^2\beta^2+224\delta_j^ir_0s_0\beta -2544\delta_j^is_ms^m_0b^2\beta^2-2496\delta_j^is_ms^m_0b^4\beta^2\\
\nonumber  \!\!\!\!\!\!&-&\!\!\!\!\!\!\ 672\delta_j^ir_{0m}s^m_0b^2\beta-1056\delta_j^ir_{0m}s^m_0b^4\beta-336\delta_j^is_{0|0}b^2\beta-528\delta_j^is_{0|0}b^4\beta +224y^ir_js_0\beta +56s^ib_jr_{00}\beta \\
\nonumber  \!\!\!\!\!\!&-&\!\!\!\!\!\!\  128y^ib_js_0^2b+3072b^ir_{jk}s^k_0\beta^2b^4 -256b^ir_{k0}s^k_j\beta^2b^6 +3552b^ir_{jk}s^k_0b^2\beta^2+384b^ir_{j0|0}\beta b^6-216s^i_0b_jr_{00}b^2 \\
\nonumber  \!\!\!\!\!\!&-&\!\!\!\!\!\!\  432s^i_0b_jr_{00}b^4 -288s^i_0b_jr_{00}b^6 +1536r^i_0s_j\beta^2b^4 +256r^i_0s_j\beta^2b^6-416r^i_0b_js_0\beta+1728s^i_ks^k_0b_jb^2\beta^2\\
\nonumber  \!\!\!\!\!\!&+&\!\!\!\!\!\!\ 1728s^i_ks^k_0b_j\beta^2b^4 +48r^i_0y_js_0b^2+96r^i_0y_js_0b^4 +64r^i_0y_js_0b^6-48r^i_0b_jr_{00}b^2 -48r^i_0b_jr_{00}b^4\\
\nonumber  \!\!\!\!\!\!&+&\!\!\!\!\!\!\ 576r^i_0r_{j0}b^2\beta+864r^i_0r_{j0}\beta b^4+384r^i_0r_{j0}\beta b^6+1776r^i_0s_jb^2\beta^2+2592b^is_js_0b^2\beta+1728b^is_js_0\beta b^4
\\
\nonumber  \!\!\!\!\!\!&-&\!\!\!\!\!\!\  576b^ib_jr_{0m}s^m_0b^2\beta +192b^ib_jr_{0m}s^m_0b^4\beta -384s^iy_js_0b^2\beta +256y^ib_jr_0s_0b^5 +320y^ib_jr_0s_0b^4 \\
\nonumber  \!\!\!\!\!\!&+&\!\!\!\!\!\!\  288b^ib_jr_{00}s_0b^4+192b^ib_jr_{00}s_0b +528b^ib_jr_{00}s_0b^2 +384b^ib_jr_{00}s_0b^3 -192y^ib_js_ms^m_0b^2\beta -128y^ib_jr_0s_0b \\
\nonumber  \!\!\!\!\!\!&-&\!\!\!\!\!\!\  672y^ib_js_ms^m_0b^4\beta +1664b^ib_js_0^2b^2\beta +3712b^ib_js_0^2b\beta+704\delta_j^ir_0s_0b^2\beta+448\delta_j^ir_0s_0b\beta  +1408\delta_j^ir_0s_0b^3\beta \\
\nonumber  \!\!\!\!\!\!&-&\!\!\!\!\!\!\  2048y^ir_js_0b^5\beta-256b^ib_js_{0|0}\beta b^6+256b^ib_jr_{0m}s^m_0b^6\beta +1024b^ib_js_0^2b^5\beta-256b^ib_js_0^2b^4\beta+768b^is_jr_{00}b\beta\\
\nonumber  \!\!\!\!\!\!&+&\!\!\!\!\!\!\ 192b^ib_jr_{00}r_0b +48b^ib_jr_{00}r_0b^2 +384b^ib_jr_{00}r_0b^3 +2304b^is_jr_{00}b^3\beta +2368b^is_js_0b\beta^2 -2368b^is_jr_0b\beta^2
\\
\nonumber  \!\!\!\!\!\!&+&\!\!\!\!\!\!\ 1536b^is_jr_{00}b^5\beta+1024b^is_jr_0\beta^2b^2-576s^iy_js_0\beta b^4+272y^ib_jr_0s_0b^2-128y^ib_jr_0s_0b^3+1408y^is_jr_0b^3\beta\\
\nonumber  \!\!\!\!\!\!&-&\!\!\!\!\!\!\ 256y^ib_js_ms^m_0b^6\beta -1536b^ir_{j0}s_0b^5\beta +240b^ir_{j0}s_0b^2\beta -768b^ir_{j0}s_0b\beta-2304b^ir_{j0}s_0b^3\beta\\
\nonumber  \!\!\!\!\!\!&+&\!\!\!\!\!\!\ 576b^ir_{j0}r_0b^2\beta  -768b^ir_{j0}r_0b\beta -2304b^ir_{j0}r_0b^3\beta +384b^ir_{j0}s_0\beta b^4+384b^ir_{j0}r_0b^4\beta -1536b^ir_{j0}r_0b^5\beta \\
\nonumber  \!\!\!\!\!\!&-&\!\!\!\!\!\!\ 2016b^ib_js_ms^m_0b^2\beta^2 +64b^iy_jr_0s_0b^4-448y^is_js_0b\beta -1408y^is_js_0b^3\beta-352y^is_jr_0b^2\beta+448y^is_jr_0b\beta\\
\nonumber  \!\!\!\!\!\!&-&\!\!\!\!\!\!\ 1440b^ib_js_{0|0}b^2\beta-1344b^ib_js_{0|0}b^4\beta-256b^iy_jr_0s_0b^5 -64b^iy_jr_0s_0b+64b^iy_jr_0s_0b^2-256b^iy_jr_0s_0b^3
\\
\nonumber  \!\!\!\!\!\!&+&\!\!\!\!\!\!\ 232y^is_js_0b^2\beta+64y^is_js_0\beta b^4-1024y^is_js_0b^5\beta -256y^is_jr_0b^4\beta+1024y^is_jr_0b^5\beta +672b^iy_js_ms^m_0b^2\beta\\
\nonumber  \!\!\!\!\!\!&+&\!\!\!\!\!\!\ 864b^iy_js_ms^m_0b^4\beta -256b^is_js_0\beta^2b^4 +432s^ib_jr_{00}b^2\beta +288s^ib_js_0b^2\beta^2 +768s^ib_jr_{00}\beta b^4+256s^ib_jr_{00}\beta b^6 \\
\nonumber  \!\!\!\!\!\!&+&\!\!\!\!\!\!\  768b^ir_jr_{00}b\beta +2304b^ir_jr_{00}b^3\beta +4736b^ir_js_0b\beta^2 -256s^iy_js_0b^6\beta-744b^is_jr_{00}b^2\beta -112\delta_j^ir_{00}s_0b^4
\\
\nonumber  \!\!\!\!\!\!&+&\!\!\!\!\!\!\ 5632b^ib_js_0^2b^3\beta +80b^ib_jr_0s_0\beta -4096b^is_jr_0b^3\beta^2 +1024b^is_js_0b^5\beta^2+256b^is_jr_0b^4\beta^2-1024b^is_jr_0b^5\beta^2\\
\nonumber  \!\!\!\!\!\!&+&\!\!\!\!\!\!\ 4096b^is_js_0b^3\beta^2+1024\delta_j^ir_0s_0b^5\beta +512\delta_j^ir_0s_0b^4\beta +256b^iy_js_ms^m_0b^6\beta+8192b^ir_js_0b^3\beta^2-512b^ir_js_0\beta^2b^4\\
\nonumber  \!\!\!\!\!\!&+&\!\!\!\!\!\!\ 2048b^ir_js_0b^5\beta^2+704y^ir_js_0b^2\beta-896y^ir_js_0b\beta -2816y^ir_js_0b^3\beta+512y^ir_js_0\beta b^4-576b^ir_jr_{00}b^2\beta\\
\nonumber  \!\!\!\!\!\!&-&\!\!\!\!\!\!\ 2048b^ir_js_0b^2\beta^2-384b^ir_jr_{00}\beta b^4+1536b^ir_jr_{00}b^5\beta-1408b^is_js_0b^2\beta^2-384b^is_jr_{00}\beta b^4-64s^i_ks^k_0y_j\beta\\
\nonumber  \!\!\!\!\!\!&+&\!\!\!\!\!\!\ 536r^i_0s_j\beta^2 +864s^is_jb^2\beta^3+432s^i_ks^k_0b_j\beta^2 +2592s^i_ks^k_jb^2\beta^3+1728s^i_ks^k_j\beta^3b^4 -158880s^i_0s_{0j}b^2\beta
\end{eqnarray}
\begin{eqnarray}
\nonumber  \!\!\!\!\!\!&-&\!\!\!\!\!\!\ 381312s^i_0s_{0j}\beta b^4-381312s^i_0s_{0j}\beta b^6-127104s^i_0s_{0j}\beta b^8+8r^i_0y_js_0-12r^i_0b_jr_{00}-24b^ir_{j0}r_{00}b^2\\
\nonumber  \!\!\!\!\!\!&+&\!\!\!\!\!\!\ 2480b^ib_js_0^2\beta -1440b^ib_js_ms^m_0\beta^2 -304b^ib_jr_{0m}s^m_0\beta -416b^ib_js_{0|0}\beta+192y^is_js_0b^6 +8y^is_jr_{00}b^4\\
\nonumber  \!\!\!\!\!\!&-&\!\!\!\!\!\!\ 48s^iy_jr_{00}b^2-48b^ir_{j0}s_0\beta -384y^ib_jr_{0m}s^m_0b^4-144y^ib_jr_{0m}s^m_0b^2-48y^ib_js_{0|0}b^4+20y^is_jr_{00}b^2
\\
\nonumber  \!\!\!\!\!\!&+&\!\!\!\!\!\!\ 16\delta_j^ir_{00}s_0b-64\delta_j^ir_{00}s_0b^2 +64\delta_j^ir_{00}s_0b^3 +16\delta_j^ir_{00}r_0b+32\delta_j^ir_{00}r_0b^2 +256s_{0|0}^ib_j\beta b^8+1776b^is_{j|0}b^2\beta^2\\
\nonumber  \!\!\!\!\!\!&+&\!\!\!\!\!\!\ 64\delta_j^ir_{00}r_0b^3+576b^ir_{j0|0}b^2\beta+864b^ir_{j0|0}\beta b^4+256y^ib_js_0^2b^5-320y^ib_jr_{0m}s^m_0b^6+672y^is_{0|j}b^2\beta\\
\nonumber  \!\!\!\!\!\!&-&\!\!\!\!\!\!\  80y^ir_{j0}s_0b^4+128y^ir_{j0}s_0b^5-32y^ir_{j0}r_0b^4 +128y^ir_{j0}r_0b^5 -128y^is_jr_{00}b^3-128y^is_jr_{00}b^5\\
\nonumber  \!\!\!\!\!\!&+&\!\!\!\!\!\!\ 64\delta_j^ir_{00}s_0b^5+32\delta_j^ir_{00}r_0b^4+64\delta_j^ir_{00}r_0b^5 +1536b^is_{j|0}\beta^2b^4+256b^is_{j|0}\beta^2b^6 +1536s_{0|0}^ib_j\beta b^4 \\
\nonumber  \!\!\!\!\!\!&-&\!\!\!\!\!\!\ 3552b^is_{0|j}b^2\beta^2-3072b^is_{0|j}\beta^2b^4 -1144b^is_js_0\beta^2-312b^is_jr_{00}\beta +704s_{0|0}^ib_jb^2\beta +1280s_{0|0}^ib_j\beta b^6
\\
\nonumber t_{11} \!\!\!\!&:=&\!\!\!\!\ 4s^i_ks^k_0y_j+31776s^i_0s_{0j}b^4 +10592s^i_0s_{0j}b^2-360s^is_j\beta^2 -312s^i_ks^k_j\beta^2+21184s^i_0s_{0j}b^8+42368s^i_0s_{0j}b^6\\
\nonumber  \!\!\!\!\!\!&-&\!\!\!\!\!\!\ 8r^i_0r_{j0}-4y^ir_{k0}s^k_j-8s_{0|0}^ib_j+8\delta_j^ir_{0m}s^m_0+4\delta_j^is_{0|0}-4\delta_j^is_0^2-8y^is_{0|j}-8b^ir_{j0|0}+8b^ir_{00|j}+64s_{0|j}^i\beta\\
\nonumber  \!\!\!\!\!\!&-&\!\!\!\!\!\!\ 128s_{j0}^i\beta+4y^is_{0|j}+8y^ir_{kj}s^k_0+8r_j^ir_{00} +96s^i_ks^k_0y_jb^4 +32s^i_ks^k_0y_jb^2+64s^i_ks^k_0y_jb^8+128s^i_ks^k_0y_jb^6\\
\nonumber  \!\!\!\!\!\!&+&\!\!\!\!\!\!\ 72s^i_0b_js_0+32\overline{R}^i_jb^2 +64\overline{R}^i_jb^8+96\overline{R}^i_jb^4 +128\overline{R}^i_jb^6-96b^is_js_0 +120\delta_j^is_ms^m_0\beta +360b^is_ms^m_0\beta^2 \\
\nonumber  \!\!\!\!\!\!&+&\!\!\!\!\!\!\  104b^ir_{k0}s^k_j\beta+64\delta_j^ir_{0m}s^m_0b^6-32\delta_j^is_0^2b+208b^is_{0|j}\beta -60y^is_ms^m_0\beta+960s_{0|j}^i\beta b^4-128\delta_j^is_0^2b^5\\
\nonumber  \!\!\!\!\!\!&+&\!\!\!\!\!\!\ 32\delta_j^is_{0|0}b^6-16\delta_j^ir_0s_0+80\delta_j^is_0^2b^4 -24y^ir_{k0}s^k_jb^2-32y^ir_{k0}s^k_jb^6+96y^ir_{kj}s^k_0b^4 +48y^ir_{kj}s^k_0b^2\\
\nonumber  \!\!\!\!\!\!&+&\!\!\!\!\!\!\ 64y^ir_{kj}s^k_0b^6+48y^is_{0|j}b^4+24y^is_{0|j}b^2+32y^is_{0|j}b^6+64b^ir_{00|j}b^6-8y^is_js_0+8s^iy_js_0+28b^is_jr_{00}\\
\nonumber  \!\!\!\!\!\!&-&\!\!\!\!\!\!\ 192s_{0|0}^ib_jb^4 -64s_{0|0}^ib_jb^2-128s_{0|0}^ib_jb^8-256s_{0|0}^ib_jb^6 -4s^ib_jr_{00}+48r_j^ir_{00}b^2+96r_j^ir_{00}b^4\\
\nonumber  \!\!\!\!\!\!&+&\!\!\!\!\!\!\ 64r_j^ir_{00}b^6+16b^ir_jr_{00}-16y^ir_js_0 -128\delta_j^is_0^2b^3 +96\delta_j^ir_{0m}s^m_0b^4+48\delta_j^ir_{0m}s^m_0b^2  +24\delta_j^is_{0|0}b^2\\
\nonumber  \!\!\!\!\!\!&+&\!\!\!\!\!\!\ 48\delta_j^is_{0|0}b^4-48y^ir_{k0}s^k_jb^4-1920s_{j0}^i\beta b^4-1792s_{j0}^i\beta b^6-512s_{j0}^i\beta b^8-96y^is_{0|j}b^4-48y^is_{0|j}b^2\\
\nonumber  \!\!\!\!\!\!&-&\!\!\!\!\!\!\ 64y^is_{0|j}b^6+8b^ir_{j0}s_0-16b^ir_{j0}r_0 +8y^is_jr_0 -8b^iy_js_ms^m_0+32b^ib_jr_{0m}s^m_0+40b^ib_js_{0|0}-496b^ib_js_0^2
\\
\nonumber  \!\!\!\!\!\!&-&\!\!\!\!\!\!\ 96b^ir_{j0|0}b^4-48b^ir_{j0|0}b^2-64b^ir_{j0|0}b^6+96b^ir_{00|j}b^4+48b^ir_{00|j}b^2-104b^is_{0|j}\beta+208r_j^is_0\beta\\
\nonumber  \!\!\!\!\!\!&+&\!\!\!\!\!\!\ 416s_{0|j}^ib^2\beta+896s_{0|j}^i\beta b^6+256s_{0|j}^i\beta b^8-832s_{j0}^ib^2\beta -208b^ir_{kj}s^k_0\beta-104s^ir_{j0}\beta+32\delta_j^is_0^2b^2
\\
\nonumber  \!\!\!\!\!\!&-&\!\!\!\!\!\!\ 432b^is_js_0b^2-576b^is_js_0b^4+64b^ir_{j0}r_0b -64b^ir_{j0}r_0b^2 +256b^ir_{j0}r_0b^3-64b^ir_{j0}s_0b^4+256b^ir_{j0}s_0b^5
\\
\nonumber  \!\!\!\!\!\!&-&\!\!\!\!\!\!\ 64b^ir_{j0}r_0b^4+256b^ir_{j0}r_0b^5+64b^ir_jr_{00}b^2 +528b^ib_js_ms^m_0\beta-64b^ir_jr_{00}b +48s^iy_js_0b^2+96s^iy_js_0b^4
\\
\nonumber  \!\!\!\!\!\!&+&\!\!\!\!\!\!\  64s^iy_js_0b^6+16y^is_js_0b^4-32\delta_j^ir_0s_0b -64\delta_j^ir_0s_0b^2+64b^ir_{j0}s_0b-16b^ir_{j0}s_0b^2 +256b^ir_{j0}s_0b^3\\
\nonumber  \!\!\!\!\!\!&+&\!\!\!\!\!\!\ 128y^is_js_0b^5+32y^is_jr_0b^4-128y^is_jr_0b^5 -128\delta_j^ir_0s_0b^5 -64\delta_j^ir_0s_0b^4-256b^is_{0|j}\beta b^6+960b^is_{0|j}b^2\beta\\
\nonumber  \!\!\!\!\!\!&+&\!\!\!\!\!\!\ 1344b^is_{0|j}\beta b^4+352b^is_js_0\beta-672b^is_{0|j}\beta b^4-512b^ib_js_0^2b^5 +128b^ib_js_{0|0}b^6-16b^ib_jr_0s_0 \\
\nonumber  \!\!\!\!\!\!&-&\!\!\!\!\!\!\ 448b^ib_js_0^2b^4  -128b^ib_jr_{0m}s^m_0b^6-512b^ib_js_0^2b -784b^ib_js_0^2b^2-1280b^ib_js_0^2b^3+96b^ib_jr_{0m}s^m_0b^2
\\
\nonumber  \!\!\!\!\!\!&+&\!\!\!\!\!\!\ 192b^ib_js_{0|0}b^2+288b^ib_js_{0|0}b^4-64b^iy_js_ms^m_0b^6-48b^iy_js_ms^m_0b^2-96b^iy_js_ms^m_0b^4-128\delta_j^ir_0s_0b^3\\
\nonumber  \!\!\!\!\!\!&+&\!\!\!\!\!\!\ 512r_j^is_0b^6\beta+320b^ir_js_0\beta +88b^is_jr_{00}b^2-192b^is_js_0b^6+64b^is_jr_{00}b^4+960r_j^is_0b^2\beta+1344r_j^is_0\beta b^4 \\
\nonumber  \!\!\!\!\!\!&-&\!\!\!\!\!\!\  160b^is_jr_0\beta-288y^is_ms^m_0b^2\beta-432y^is_ms^m_0\beta b^4-256b^ir_jr_{00}b^3+576\delta_j^is_ms^m_0b^2\beta+864\delta_j^is_ms^m_0b^4\beta\\
\nonumber  \!\!\!\!\!\!&-&\!\!\!\!\!\!\ 960b^ir_{kj}s^k_0b^2\beta+64b^ir_jr_{00}b^4-256b^ir_jr_{00}b^5 +64y^ir_js_0b-64y^ir_js_0b^2+256y^ir_js_0b^3-64y^ir_js_0b^4\\
\nonumber  \!\!\!\!\!\!&+&\!\!\!\!\!\!\ 256y^ir_js_0b^5+192y^ib_js_ms^m_0b^4+48y^ib_js_ms^m_0b^2 -64b^is_jr_{00}b+4\overline{R}^i_j-576s^i_ks^k_0b_jb^2\beta+1324s^i_0s_{0j}\\
\nonumber  \!\!\!\!\!\!&+&\!\!\!\!\!\!\ 432s^i_0b_js_0b^2 +864s^i_0b_js_0b^4+576s^i_0b_js_0b^6-256r^i_0s_j\beta b^6-1152s^i_ks^k_0b_j\beta b^4-768s^i_ks^k_0b_j\beta b^6\\
\nonumber  \!\!\!\!\!\!&+&\!\!\!\!\!\!\ 192r^i_0b_js_0b^2 +288r^i_0b_js_0b^4 +128r^i_0b_js_0b^6 -480r^i_0s_jb^2\beta -672r^i_0s_j\beta b^4+960b^ib_js_ms^m_0b^4\beta \\
\nonumber  \!\!\!\!\!\!&-&\!\!\!\!\!\!\ 512b^ib_jr_0s_0b^5  +128b^ib_jr_0s_0b^4-512b^ib_jr_0s_0b +32b^ib_jr_0s_0b^2-640b^is_js_0b\beta-1792b^is_js_0b^3\beta \\
\nonumber  \!\!\!\!\!\!&-&\!\!\!\!\!\!\ 448b^is_jr_0b^2\beta+640b^is_jr_0b\beta +1792b^is_jr_0b^3\beta-1280b^ib_jr_0s_0b^3 +1536b^ib_js_ms^m_0b^2\beta-96s^ib_js_0b^2\beta\\
\nonumber  \!\!\!\!\!\!&-&\!\!\!\!\!\!\ 384s^ib_js_0\beta b^4  +896b^ir_js_0b^2\beta -1280b^ir_js_0b\beta-3584b^ir_js_0b^3\beta +688b^is_js_0b^2\beta+1024b^is_jr_0b^5\beta
\end{eqnarray}
\begin{eqnarray}
\nonumber  \!\!\!\!\!\!&-&\!\!\!\!\!\!\ 1024b^is_js_0b^5\beta-256b^is_jr_0b^4\beta-2048b^ir_js_0b^5\beta+512b^ir_js_0\beta b^4+256b^is_js_0\beta b^4-104r^i_0s_j\beta \\
\nonumber  \!\!\!\!\!\!&-&\!\!\!\!\!\!\ 1008s^is_jb^2\beta^2 -576s^is_j\beta^2b^4 -96s^i_ks^k_0b_j\beta-1440s^i_ks^k_jb^2\beta^2 -2016s^i_ks^k_j\beta^2b^4 -768s^i_ks^k_j\beta^2b^6\\
\nonumber  \!\!\!\!\!\!&-&\!\!\!\!\!\!\ 96r^i_0r_{j0}b^4-48r^i_0r_{j0}b^2-64r^i_0r_{j0}b^6 -256b^is_jr_{00}b^3-256b^is_jr_{00}b^5 -192y^is_ms^m_0\beta b^6
\\
\nonumber  \!\!\!\!\!\!&+&\!\!\!\!\!\!\ 384\delta_j^is_ms^m_0b^6\beta +576b^is_ms^m_0\beta^2b^4 +192y^ib_js_ms^m_0b^6+512b^is_{0|j}\beta b^6+1008b^is_ms^m_0b^2\beta^2 \\
\nonumber  \!\!\!\!\!\!&-&\!\!\!\!\!\!\ 256s^ir_{j0}\beta b^6 -512b^ir_{kj}s^k_0\beta b^6-480s^ir_{j0}b^2\beta +480b^ir_{k0}s^k_jb^2\beta -144s^ib_jr_{00}b^4 +40r^i_0b_js_0
\\
\nonumber  \!\!\!\!\!\!&-&\!\!\!\!\!\!\ 128y^is_jr_0b^3-48s^ib_jr_{00}b^2-672s^ir_{j0}\beta b^4-128s^ib_jr_{00}b^6+32y^is_js_0b-8y^is_js_0b^2+128y^is_js_0b^3\\
\nonumber  \!\!\!\!\!\!&-&\!\!\!\!\!\!\ 32y^is_jr_0b+32y^is_jr_0b^2-480b^is_{0|j}b^2\beta+48s^ib_js_0\beta -1344b^ir_{kj}s^k_0\beta b^4+672b^ir_{k0}s^k_j\beta b^4+256b^ir_{k0}s^k_j\beta b^6
\\
\nonumber t_{12} \!\!\!\!&:=&\!\!\!\!\ -(4(2b^2+1))(-2s^i_ks^k_0b_j-22s^is_j\beta -14s^i_ks^k_j\beta -2r^i_0s_j-2b^is_{0|j}+4b^is_{0|j}+2b^ir_{k0}s^k_j+2\delta_j^is_ms^m_0\\
\nonumber  \!\!\!\!\!\!&+&\!\!\!\!\!\!\ 12s_{0|j}^ib^4+6s_{0|j}^ib^2+8s_{0|j}^ib^6-24s_{j0}^ib^4-y^is_ms^m_0-2s^ir_{j0} -4b^ir_{kj}s^k_0 +4r_j^is_0-16s_{j0}^ib^6  -12s_{j0}^ib^2 \\
\nonumber  \!\!\!\!\!\!&+&\!\!\!\!\!\!\  8b^ir_{k0}s^k_jb^4 +8b^ir_{k0}s^k_jb^2-16b^ir_{kj}s^k_0b^4 -16b^ir_{kj}s^k_0b^2 -8s^ir_{j0}b^4-8s^ir_{j0}b^2-8b^is_{0|j}b^2 +22b^is_ms^m_0\beta \\
\nonumber  \!\!\!\!\!\!&+&\!\!\!\!\!\!\ 10b^is_js_0-8b^is_{0|j}b^4+2s^ib_js_0+16r_j^is_0b^2+16r_j^is_0b^4+8b^ir_js_0 +8\delta_j^is_ms^m_0b^2+8\delta_j^is_ms^m_0b^4 \\
\nonumber  \!\!\!\!\!\!&+&\!\!\!\!\!\!\ 16b^is_{0|j}b^4-4b^is_jr_0-4y^is_ms^m_0b^4 -4y^is_ms^m_0b^2 +16b^ib_js_ms^m_0+16b^is_{0|j}b^2-32b^ir_js_0b \\
\nonumber  \!\!\!\!\!\!&-&\!\!\!\!\!\!\  32b^is_js_0b^3+40b^ib_js_ms^m_0b^2 +16b^ib_js_ms^m_0b^4+8b^is_js_0b^2 +16b^is_jr_0b-8b^is_jr_0b^2+32b^is_jr_0b^3 \\
\nonumber  \!\!\!\!\!\!&-&\!\!\!\!\!\!\ 16b^ir_js_0b^2+52b^is_ms^m_0b^2\beta+16b^is_ms^m_0\beta b^4-16s^ib_js_0b^4 -64b^ir_js_0b^3-2s_{j0}^i+s_{0|j}^i \\
\nonumber  \!\!\!\!\!\!&-&\!\!\!\!\!\!\ 4s^ib_js_0b^2-8r^i_0s_jb^4-52s^is_jb^2\beta-16s^is_j\beta b^4-60s^i_ks^k_jb^2\beta -72s^i_ks^k_j\beta b^4 -16b^is_js_0b\\
\nonumber  \!\!\!\!\!\!&-&\!\!\!\!\!\!\ 16s^i_ks^k_j\beta b^6-8r^i_0s_jb^2-24s^i_ks^k_0b_jb^4-12s^i_ks^k_0b_jb^2-16s^i_ks^k_0b_jb^6)\\
\nonumber t_{13} \!\!\!\!&:=&\!\!\!\!\ 4(2b^2+1)^3(-2s^i_ks^k_jb^2-s^i_ks^k_j-2s^is_j+2b^is_ms^m_j).
\end{eqnarray}
\section{Appendix 2: Coefficients in (\ref{a8})}
\begin{eqnarray}
\nonumber d_0 \!\!\!\!&:=&\!\!\!\!\ -288r_{00}^2(8n-11)\beta^7,\\
\nonumber d_1 \!\!\!\!&:=&\!\!\!\!\ 48\beta^6(-54r_{00|0}\beta +80r_{00}^2nb^2-163r_{00}^2 -104r_{00}^2b^2+124r_{00}^2n+36r_{00|0}n\beta),
\\
\nonumber d_2 \!\!\!\!&:=&\!\!\!\!\ -24\beta^5(108{\bf \overline{Ric}}\beta^2+384r_{00}s_0b\beta+96r_{00}s_0n\beta+198r_{00|0}n\beta-192r_{00}s_0nb\beta
-192r_{00}r_0nb\beta\\
\nonumber  \!\!\!\!\!\!&-&\!\!\!\!\!\!\ 276r_{00|0}\beta-120r_{00}s_0\beta +384r_{00}r_0b\beta-96r_{00}r_0n\beta-192r_{00|0}\beta b^2+144r_{00}r_0\beta+144r_{00|0}n\beta b^2 \\
 \nonumber  \!\!\!\!\!\!&+&\!\!\!\!\!\!\ 64r_{00}^2nb^4 -64r_{00}^2b^4 -512r_{00}^2b^2 +392r_{00}^2nb^2 -360r_{00}^2+273r_{00}^2n)
\\
\nonumber d_3 \!\!\!\!&:=&\!\!\!\!\ 12\beta^4(784r_{00}r_0\beta+462r_{00|0}n\beta-224r_{00|0}b^4\beta -112r_{00}s_0\beta -800r_{00|0}\beta b^2+288r_{00}^2nb^4\\
 \nonumber  \!\!\!\!\!\!&+&\!\!\!\!\!\!\  808r_{00}^2nb^2+288r_{0m}s^m_0n\beta^2 +144s_{0|0}n\beta^2+144r_{00}r_m^m\beta^2+144r_{00|m}b^m\beta^2-144r_{0m|0}b^m\beta^2 \\
 \nonumber  \!\!\!\!\!\!&+&\!\!\!\!\!\!\  336r_{00}s_0n\beta+1632r_{00}s_0b\beta-512r_{00}s_0\beta b^2+768r_{00}s_0\beta b^3-464r_{00}r_0n\beta +1632r_{00}r_0b\beta  \\
  \nonumber  \!\!\!\!\!\!&-&\!\!\!\!\!\!\ 512r_{00}s_0nb^3\beta+384r_{00}s_0nb^2\beta-928r_{00}r_0nb\beta-512r_{00}r_0nb^3\beta -256r_{00}r_0nb^2\beta+320r_{00}r_0\beta b^2 \\
 \nonumber  \!\!\!\!\!\!&+&\!\!\!\!\!\!\ 696r_{00|0}n\beta b^2+192r_{00|0}n\beta b^4-468r_{00}^2-256r_{00}^2b^4 -1184r_{00}^2b^2+335r_{00}^2n-596r_{00|0}\beta +768r_{00}r_0\beta b^3 \\
  \nonumber  \!\!\!\!\!\!&-&\!\!\!\!\!\!\ 144r_{0m}r^m_0\beta^2
 -288r_{0m}s^m_0\beta^2 -360s_{0|0}\beta^2+612{\bf \overline{Ric}}\beta^2+576{\bf \overline{Ric}}\beta^2b^2 -928r_{00}s_0nb\beta)
\\
\nonumber d_4 \!\!\!\!&:=&\!\!\!\!\ -2\beta^3(5360r_{00}r_0\beta+1778r_{00|0}n\beta-2112r_{00|0}b^4\beta-256r_{00|0}b^6\beta+2144r_{00}s_0\beta -3936r_{00|0}\beta b^2\\
 \nonumber  \!\!\!\!\!\!&+&\!\!\!\!\!\!\  6912s_0^2b\beta^2+576rr_{00}\beta^2+2304r_0^2b\beta^2 -288r_0s_0\beta^2 +1584r_{00}^2nb^4+2712r_{00}^2nb^2\\
  \nonumber  \!\!\!\!\!\!&-&\!\!\!\!\!\!\  1728r_{0m}r^m_0b^2\beta^2  +3888r_{0m}s^m_0n\beta^2-2880r_{0m}s^m_0\beta^2b^2 +1944s_{0|0}n\beta^2 -4032s_{0|0}\beta^2b^2 \\
   \nonumber  \!\!\!\!\!\!&+&\!\!\!\!\!\!\ 2160r_{00}r_m^m\beta^2 +2160r_{00|m}b^m\beta^2
 -2160r_{0m|0}b^m\beta^2 +1384r_{00}s_0n\beta -2304s_0^2nb\beta^2\\
 \nonumber  \!\!\!\!\!\!&+&\!\!\!\!\!\!\  8416r_{00}s_0b\beta -4288r_{00}s_0\beta b^2+1024r_{00}s_0\beta b^5+6400r_{00}s_0\beta b^3-2768r_{00}r_0n\beta \\
 \nonumber  \!\!\!\!\!\!&+&\!\!\!\!\!\!\  8416r_{00}r_0b\beta+1024r_{00}s_0nb^4\beta -5536r_{00}s_0nb\beta -6400r_{00}s_0nb^3\beta -1024r_{00}s_0nb^5\beta \\
  \nonumber  \!\!\!\!\!\!&+&\!\!\!\!\!\!\ 4288r_{00}s_0nb^2\beta-512r_{00}r_0nb^4\beta -5536r_{00}r_0nb\beta -6400r_{00}r_0nb^3\beta -1024r_{00}r_0nb^5\beta
\end{eqnarray}
\begin{eqnarray}
\nonumber  \!\!\!\!\!\!&-&\!\!\!\!\!\!\ 3200r_{00}r_0nb^2\beta-2304r_0s_0nb\beta^2+4352r_{00}r_0\beta b^2+512r_{00}r_0\beta b^4 +1024r_{00}r_0\beta b^5 \\
\nonumber  \!\!\!\!\!\!&+&\!\!\!\!\!\!\ 6400r_{00}r_0\beta b^3+4152r_{00|0}n\beta b^2+2400r_{00|0}n\beta b^4+256r_{00|0}n\beta b^6-1024r_{00}s_0b^4\beta  \\
\nonumber  \!\!\!\!\!\!&-&\!\!\!\!\!\!\ 2304rr_{00}b\beta^2-1152r_0s_0n\beta^2 +8064r_0s_0b\beta^2 +3456r_{0m}s^m_0n\beta^2b^2 -1149r_{00}^2 \\
\nonumber  \!\!\!\!\!\!&-&\!\!\!\!\!\!\ 1440r_{00}^2b^4-4944r_{00}^2b^2 +753r_{00}^2n-2120r_{00|0}\beta-144s_0^2\beta^2 -576r_0^2\beta^2-2160r_{0m}r^m_0\beta^2
\\
\nonumber  \!\!\!\!\!\!&-&\!\!\!\!\!\!\ 3168r_{0m}s^m_0\beta^2-4824s_{0|0}\beta^2 -1296s_{0|m}^m\beta^3 +4428{\bf \overline{Ric}}\beta^2 +3456{\bf \overline{Ric}}\beta^2b^4
\\
\nonumber  \!\!\!\!\!\!&+&\!\!\!\!\!\!\ 8640{\bf \overline{Ric}}\beta^2b^2 +1728s_{0|0}n\beta^2b^2 +1728r_{00}r_m^mb^2\beta^2 +1728r_{00|m}b^mb^2\beta^2-1728r_{0m|0}b^mb^2\beta^2)
\\
\nonumber d_5 \!\!\!\!&:=&\!\!\!\!\ -\beta^2(-6544r_{00}r_0\beta-1354r_{00|0}n\beta+2112r_{00|0}b^4\beta+512r_{00|0}b^6\beta -5488r_{00}s_0\beta +3168r_{00|0}\beta b^2 \\
\nonumber  \!\!\!\!\!\!&+&\!\!\!\!\!\!\ 1632s_0^2n\beta^2+1920s_0^2\beta^2b^2 -15360s_0^2b^3\beta^2 -21504s_0^2b\beta^2-2496rr_{00}\beta^2 +2304r_{0m|0}b^mb^4\beta^2\\
\nonumber  \!\!\!\!\!\!&+&\!\!\!\!\!\!\ 1536r_0^2b^2\beta^2-9984r_0^2b\beta^2-6144r_0^2b^3\beta^2 -2400r_0s_0\beta^2 -1488r_{00}^2nb^4-1752r_{00}^2nb^2\\
\nonumber  \!\!\!\!\!\!&+&\!\!\!\!\!\!\ 7488r_{0m}r^m_0b^2\beta^2+2304r_{0m}r^m_0b^4\beta^2 -7200r_{0m}s^m_0n\beta^2 +7104r_{0m}s^m_0\beta^2b^2 \\
\nonumber  \!\!\!\!\!\!&+&\!\!\!\!\!\!\ 3072r_{0m}s^m_0\beta^2b^4-3600s_{0|0}n\beta^2+14784s_{0|0}\beta^2b^2+4992s_{0|0}\beta^2b^4-1728s_ms^m_0n\beta^3\\
\nonumber  \!\!\!\!\!\!&-&\!\!\!\!\!\!\ 4464r_{00}r_m^m\beta^2-4464r_{00|m}b^m\beta^2+4464r_{0m|0}b^m\beta^2 -648s_ms_0^m\beta^2-1728s_0r_m^m\beta^3\\
\nonumber  \!\!\!\!\!\!&-&\!\!\!\!\!\!\ 1728s_{0|m}b^m\beta^3+864s_{m|0}b^m\beta^3+6912s_{0|m}^m\beta^3b^2 -992r_{00}s_0n\beta-7520r_{00}s_0b\beta \\
\nonumber  \!\!\!\!\!\!&+&\!\!\!\!\!\!\ 5312r_{00}s_0\beta b^2-512r_{00}s_0\beta b^5-4352r_{00}s_0\beta b^3+2896r_{00}r_0n\beta -7520r_{00}r_0b\beta -3584r_{00}s_0nb^4\beta \\
\nonumber  \!\!\!\!\!\!&+&\!\!\!\!\!\!\ 5792r_{00}s_0nb\beta+10496r_{00}s_0nb^3\beta+3584r_{00}s_0nb^5\beta -6656r_{00}s_0nb^2\beta+1792r_{00}r_0nb^4\beta \\
\nonumber  \!\!\!\!\!\!&+&\!\!\!\!\!\!\ 5792r_{00}r_0nb\beta+10496r_{00}r_0nb^3\beta +3584r_{00}r_0nb^5\beta +5248r_{00}r_0nb^2\beta +8832r_0s_0nb\beta^2
\\
\nonumber  \!\!\!\!\!\!&+&\!\!\!\!\!\!\ 6144r_0s_0nb^3\beta^2+3072r_0s_0n\beta^2b^2-7936r_{00}r_0\beta b^2-1792r_{00}r_0\beta b^4-512r_{00}r_0\beta b^5\\
\nonumber  \!\!\!\!\!\!&-&\!\!\!\!\!\!\ 4352r_{00}r_0\beta b^3-4344r_{00|0}n\beta b^2-3936r_{00|0}n\beta b^4-896r_{00|0}n\beta b^6+2048r_{00}s_0b^4\beta+8832s_0^2nb\beta^2\\
\nonumber  \!\!\!\!\!\!&+&\!\!\!\!\!\!\ 6144s_0^2nb^3\beta^2-768s_0^2n\beta^2b^2 -1536rr_{00}b^2\beta^2 +9984rr_{00}b\beta^2 +6144rr_{00}b^3\beta^2 \\
\nonumber  \!\!\!\!\!\!&+&\!\!\!\!\!\!\  4416r_0s_0n\beta^2+1536r_0s_0\beta^2b^2-18432r_0s_0b^3\beta^2 -28800r_0s_0b\beta^2 -13248r_{0m}s^m_0n\beta^2b^2 \\
\nonumber  \!\!\!\!\!\!&-&\!\!\!\!\!\!\  4608r_{0m}s^m_0n\beta^2b^4 +567r_{00}^2+1632r_{00}^2b^4+4200r_{00}^2b^2 -351r_{00}^2n+1498r_{00|0}\beta  -5520s_0^2\beta^2
\\
\nonumber  \!\!\!\!\!\!&+&\!\!\!\!\!\!\  2496r_0^2\beta^2+4464r_{0m}r^m_0\beta^2+4512r_{0m}s^m_0\beta^2+8952s_{0|0}\beta^2+3456s_ms^m_0\beta^3 \\
\nonumber  \!\!\!\!\!\!&-&\!\!\!\!\!\!\ 2592s_mr^m_0\beta^3+1728r_ms^m_0\beta^3+6048s_{0|m}^m\beta^3 -5892{\bf \overline{Ric}}\beta^2-14976{\bf \overline{Ric}}\beta^2b^4
\\
\nonumber  \!\!\!\!\!\!&-&\!\!\!\!\!\!\ 17856{\bf \overline{Ric}}\beta^2b^2-3072{\bf \overline{Ric}}\beta^2b^6-6624s_{0|0}n\beta^2b^2 -2304s_{0|0}n\beta^2b^4 +7488r_{0m|0}b^mb^2\beta^2\\
\nonumber  \!\!\!\!\!\!&-&\!\!\!\!\!\!\ 7488r_{00}r_m^mb^2\beta^2-2304r_{00}r_m^mb^4\beta^2-7488r_{00|m}b^mb^2\beta^2 -2304r_{00|m}b^mb^4\beta^2)
\\
\nonumber d_6 \!\!\!\!&:=&\!\!\!\!\ 2\beta(-1128r_{00}r_0\beta-153r_{00|0}n\beta+96r_{00|0}b^4\beta -1386r_{00}s_0\beta +312r_{00|0}\beta b^2+1148s_0^2n\beta^2 \\
\nonumber  \!\!\!\!\!\!&-&\!\!\!\!\!\!\  2976s_0^2\beta^2b^2 +384s_0^2\beta^2b^4 -2048s_0^2b^5\beta^2-6272s_0^2b^3\beta^2 -6368s_0^2b\beta^2-1072rr_{00}\beta^2 \\
\nonumber  \!\!\!\!\!\!&+&\!\!\!\!\!\!\  1408r_0^2b^2\beta^2+256r_0^2b^4\beta^2 -4288r_0^2b\beta^2-5632r_0^2b^3\beta^2 -1024r_0^2b^5\beta^2-2208r_0s_0\beta^2\\
\nonumber  \!\!\!\!\!\!&-&\!\!\!\!\!\!\ 180r_{00}^2nb^4-156r_{00}^2nb^2 +3216r_{0m}r^m_0b^2\beta^2 +2112r_{0m}r^m_0b^4\beta^2 +256r_{0m}r^m_0b^6\beta^2-1728s_0r_m^m\beta^3b^2 \\
\nonumber  \!\!\!\!\!\!&-&\!\!\!\!\!\!\ 1756r_{0m}s^m_0n\beta^2+1008r_{0m}s^m_0\beta^2b^2 +768r_{0m}s^m_0\beta^2b^4 +256r_{0m}s^m_0b^6\beta^2 -878s_{0|0}n\beta^2\\
\nonumber  \!\!\!\!\!\!&+&\!\!\!\!\!\!\ 5328s_{0|0}\beta^2b^2+3552s_{0|0}\beta^2b^4+512s_{0|0}b^6\beta^2-576rs_0\beta^3 -1512s_ms^m_0n\beta^3 +3168s_ms^m_0\beta^3b^2 \\
\nonumber  \!\!\!\!\!\!&+&\!\!\!\!\!\!\ 2592s_mr^m_0\beta^3b^2+1728r_ms^m_0\beta^3b^2-1220r_{00}r_m^m\beta^2-1220r_{00|m}b^m\beta^2+1220r_{0m|0}b^m\beta^2 \\
\nonumber  \!\!\!\!\!\!&-&\!\!\!\!\!\!\ 432s_ms_0^m\beta^2-1728s_0r_m^m\beta^3-1728s_{0|m}b^m\beta^3+864s_{m|0}b^m\beta^3 +6912s_{0|m}^m\beta^3b^2 +864s_{m|0}b^m\beta^3b^2\\
\nonumber  \!\!\!\!\!\!&+&\!\!\!\!\!\!\ 3456s_{0|m}^m\beta^3b^4-114r_{00}s_0n\beta-936r_{00}s_0b\beta+1200r_{00}s_0\beta b^2 +1152r_{00}s_0\beta b^5 -1728s_{0|m}b^m\beta^3b^2\\
\nonumber  \!\!\!\!\!\!&+&\!\!\!\!\!\!\ 576r_{00}s_0\beta b^3+420r_{00}r_0n\beta-936r_{00}r_0b\beta -1248r_{00}s_0nb^4\beta +840r_{00}s_0nb\beta +2112r_{00}s_0nb^3\beta \\
\nonumber  \!\!\!\!\!\!&+&\!\!\!\!\!\!\ 1152r_{00}s_0nb^5\beta-1392r_{00}s_0nb^2\beta +576r_{00}r_0nb^4\beta +840r_{00}r_0nb\beta +2112r_{00}r_0nb^3\beta \\
\nonumber  \!\!\!\!\!\!&+&\!\!\!\!\!\!\ 1152r_{00}r_0nb^5\beta+1056r_{00}r_0nb^2\beta +3328r_0s_0nb\beta^2 +4864r_0s_0nb^3\beta^2 +1024r_0s_0nb^5\beta^2 \\
\nonumber  \!\!\!\!\!\!&+&\!\!\!\!\!\!\ 2432r_0s_0n\beta^2b^2+512r_0s_0n\beta^2b^4-1824r_{00}r_0\beta b^2-576r_{00}r_0\beta b^4+1152r_{00}r_0\beta b^5 -864s_ms_0^mb^2\beta^2\\
\nonumber  \!\!\!\!\!\!&+&\!\!\!\!\!\!\ 576r_{00}r_0\beta b^3-630r_{00|0}n\beta b^2-792r_{00|0}n\beta b^4-288r_{00|0}n\beta b^6+384r_{00}s_0b^4\beta+3328s_0^2nb\beta^2\\
\nonumber  \!\!\!\!\!\!&+&\!\!\!\!\!\!\ 4864s_0^2nb^3\beta^2+1024s_0^2nb^5\beta^2 +1376s_0^2n\beta^2b^2 -256s_0^2n\beta^2b^4 -1408rr_{00}b^2\beta^2 +256r_{0m|0}b^mb^6\beta^2
 \end{eqnarray}
\begin{eqnarray}
\nonumber  \!\!\!\!\!\!&-&\!\!\!\!\!\!\ 256rr_{00}b^4\beta^2+4288rr_{00}b\beta^2 +5632rr_{00}b^3\beta^2 +1024rr_{00}b^5\beta^2 +1664r_0s_0n\beta^2 +2112r_{0m|0}b^mb^4\beta^2 \\
\nonumber  \!\!\!\!\!\!&-&\!\!\!\!\!\!\ 2208r_0s_0\beta^2b^2 +384r_0s_0\beta^2b^4 -2560r_0s_0b^5\beta^2 -11776r_0s_0b^3\beta^2-10144r_0s_0b\beta^2\\
\nonumber  \!\!\!\!\!\!&-&\!\!\!\!\!\!\ 4992r_{0m}s^m_0n\beta^2b^2 -3648r_{0m}s^m_0n\beta^2b^4 +36r_{00}^2+264r_{00}^2b^4 +492r_{00}^2b^2 +3216r_{0m|0}b^mb^2\beta^2\\
\nonumber  \!\!\!\!\!\!&-&\!\!\!\!\!\!\ 24r_{00}^2n+159r_{00|0}\beta -4932s_0^2\beta^2 +1072r_0^2\beta^2 +1220r_{0m}r^m_0\beta^2 +776r_{0m}s^m_0\beta^2  -256r_{00|m}b^mb^6\beta^2 \\
\nonumber  \!\!\!\!\!\!&+&\!\!\!\!\!\!\ 2218s_{0|0}\beta^2+2880s_ms^m_0\beta^3+2592s_mr^m_0\beta^3 +1728r_ms^m_0\beta^3+2916s_{0|m}^m\beta^3 -2112r_{00}r_m^mb^4\beta^2
\\
\nonumber  \!\!\!\!\!\!&-&\!\!\!\!\!\!\  1168{\bf \overline{Ric}}\beta^2+324s_m^is_i^m\beta^4 -6432{\bf \overline{Ric}}\beta^2b^4 -256{\bf \overline{Ric}}\beta^2b^8 -4880{\bf \overline{Ric}}\beta^2b^2 -256r_{00}r_m^mb^6\beta^2\\
\nonumber  \!\!\!\!\!\!&-&\!\!\!\!\!\!\  2816{\bf \overline{Ric}}\beta^2b^6-512r_{0m}s^m_0n\beta^2b^6-2496s_{0|0}n\beta^2b^2-1824s_{0|0}n\beta^2b^4 -3216r_{00|m}b^mb^2\beta^2\\
\nonumber  \!\!\!\!\!\!&-&\!\!\!\!\!\!\ 256s_{0|0}n\beta^2b^6+2304rs_0b\beta^3-1728s_ms^m_0n\beta^3b^2 -3216r_{00}r_m^mb^2\beta^2 -2112r_{00|m}b^mb^4\beta^2)
\\
\nonumber d_7 \!\!\!\!&:=&\!\!\!\!\ (416r_{00}r_0\beta+38r_{00|0}n\beta+120r_{00|0}b^4\beta+128r_{00|0}b^6\beta+680r_{00}s_0\beta-48r_{00|0}\beta b^2 -1172s_0^2n\beta^2 \\
\nonumber  \!\!\!\!\!\!&+&\!\!\!\!\!\!\  9280s_0^2\beta^2b^2 +6016s_0^2\beta^2b^4-2048s_0^2b^5\beta^2 -128s_0^2b^3\beta^2 +3616s_0^2b\beta^2 +912rr_{00}\beta^2 +24r_{00}^2nb^2\\
 \nonumber  \!\!\!\!\!\!&-&\!\!\!\!\!\!\ 1920r_0^2b^2\beta^2 -768r_0^2b^4\beta^2+3648r_0^2b\beta^2 +7680r_0^2b^3\beta^2 +3072r_0^2b^5\beta^2 -2736r_{0m}r^m_0b^2\beta^2\\
 \nonumber  \!\!\!\!\!\!&+&\!\!\!\!\!\!\ 2624r_0s_0\beta^2 36r_{00}^2nb^4-2880r_{0m}r^m_0b^4\beta^2 -768r_{0m}r^m_0b^6\beta^2 +952r_{0m}s^m_0n\beta^2-3840s_ms^m_0\beta^3b^4 \\
 \nonumber  \!\!\!\!\!\!&+&\!\!\!\!\!\!\ 624r_{0m}s^m_0\beta^2b^2 +1536r_{0m}s^m_0\beta^2b^4 +256r_{0m}s^m_0b^6\beta^2 +476s_{0|0}n\beta^2-3792s_{0|0}\beta^2b^2 \\
  \nonumber  \!\!\!\!\!\!&+&\!\!\!\!\!\!\ -3552s_{0|0}\beta^2b^4-1024s_{0|0}b^6\beta^2+1920rs_0\beta^3+2088s_ms^m_0n\beta^3 -7872s_ms^m_0\beta^3b^2 \\
\nonumber  \!\!\!\!\!\!&-&\!\!\!\!\!\!\ 8640s_mr^m_0\beta^3b^2-3456s_mr^m_0\beta^3b^4 -5760r_ms^m_0\beta^3b^2 -2304r_ms^m_0\beta^3b^4 +744r_{00}r_m^m\beta^2
\\
\nonumber  \!\!\!\!\!\!&+&\!\!\!\!\!\!\ 744r_{00|m}b^m\beta^2-744r_{0m|0}b^m\beta^2 +432s_ms_0^m\beta^2 +2736s_0r_m^m\beta^3 +2736s_{0|m}b^m\beta^3 -1664r_{00}r_0\beta b^5\\
\nonumber  \!\!\!\!\!\!&-&\!\!\!\!\!\!\ 1368s_{m|0}b^m\beta^3-10944s_{0|m}^m\beta^3b^2-11520s_{0|m}^m\beta^3b^4-3072s_{0|m}^m\beta^3b^6 -5632r_0s_0nb^3\beta^2\\
\nonumber  \!\!\!\!\!\!&+&\!\!\!\!\!\!\ 40r_{00}s_0n\beta+256r_{00}s_0b\beta-784r_{00}s_0\beta b^2-1664r_{00}s_0\beta b^5-896r_{00}s_0\beta b^3-128r_{00}r_0n\beta +256r_{00}r_0b\beta
\\
\nonumber  \!\!\!\!\!\!&+&\!\!\!\!\!\!\  832r_{00}s_0nb^4\beta-256r_{00}s_0nb\beta-832r_{00}s_0nb^3\beta -640r_{00}s_0nb^5\beta+640r_{00}s_0nb^2\beta -320r_{00}r_0nb^4\beta \\
\nonumber  \!\!\!\!\!\!&-&\!\!\!\!\!\!\  256r_{00}r_0nb\beta -832r_{00}r_0nb^3\beta-640r_{00}r_0nb^5\beta -416r_{00}r_0nb^2\beta -2464r_0s_0nb\beta^2 -896r_{00}r_0\beta b^3 \\
\nonumber  \!\!\!\!\!\!&-&\!\!\!\!\!\!\ 2560r_0s_0nb^5\beta^2-2816r_0s_0n\beta^2b^2-1280r_0s_0n\beta^2b^4 +848r_{00}r_0\beta b^2+320r_{00}r_0\beta b^4 -1728s_m^is_i^m\beta^4b^2
\\
\nonumber  \!\!\!\!\!\!&+&\!\!\!\!\!\!\ 192r_{00|0}n\beta b^2+312r_{00|0}n\beta b^4+160r_{00|0}n\beta b^6-256r_{00}s_0b^4\beta -2464s_0^2nb\beta^2 -5632s_0^2nb^3\beta^2
\\
\nonumber  \!\!\!\!\!\!&-&\!\!\!\!\!\!\  2560s_0^2nb^5\beta^2-3392s_0^2n\beta^2b^2-3968s_0^2n\beta^2b^4 +1920rr_{00}b^2\beta^2+768rr_{00}b^4\beta^2  -3648rr_{00}b\beta^2 \\
\nonumber  \!\!\!\!\!\!&-&\!\!\!\!\!\!\  7680rr_{00}b^3\beta^2-3072rr_{00}b^5\beta^2 -1232r_0s_0n\beta^2 +6368r_0s_0\beta^2b^2 +3968r_0s_0\beta^2b^4+3584r_0s_0b^5\beta^2 \\
\nonumber  \!\!\!\!\!\!&+&\!\!\!\!\!\!\  9728r_0s_0b^3\beta^2+7136r_0s_0b\beta^2 +3696r_{0m}s^m_0n\beta^2b^2 +4224r_{0m}s^m_0n\beta^2b^4-3r_{00}^2-72r_{00}^2b^4 \\
\nonumber  \!\!\!\!\!\!&-&\!\!\!\!\!\!\ 96r_{00}^2b^2+3r_{00}^2n-38r_{00|0}\beta+6340s_0^2\beta^2-912r_0^2\beta^2 -744r_{0m}r^m_0\beta^2-256r_{0m}s^m_0\beta^2  -1244s_{0|0}\beta^2 \\
\nonumber  \!\!\!\!\!\!&-&\!\!\!\!\!\!\  3840s_ms^m_0\beta^3-4104s_mr^m_0\beta^3-2736r_ms^m_0\beta^3-2976s_{0|m}^m\beta^3+552{\bf \overline{Ric}}\beta^2
-864s_ms^m\beta^4 \\
\nonumber  \!\!\!\!\!\!&-&\!\!\!\!\!\!\ 1188s_m^is_i^m\beta^4 +5472{\bf \overline{Ric}}\beta^2b^4+768{\bf \overline{Ric}}\beta^2b^8 +2976{\bf \overline{Ric}}\beta^2b^2+3840{\bf \overline{Ric}}\beta^2b^6 -6144rs_0\beta^3b^3\\
\nonumber  \!\!\!\!\!\!&+&\!\!\!\!\!\!\ 1280r_{0m}s^m_0n\beta^2b^6+1848s_{0|0}n\beta^2b^2+2112s_{0|0}n\beta^2b^4 +640s_{0|0}n\beta^2b^6+1536rs_0\beta^3b^2 \\
\nonumber  \!\!\!\!\!\!&-&\!\!\!\!\!\!\  7680rs_0b\beta^3+4896s_ms^m_0n\beta^3b^2+2304s_ms^m_0n\beta^3b^4+2736r_{00}r_m^mb^2\beta^2 +2880r_{00}r_m^mb^4\beta^2\\
\nonumber  \!\!\!\!\!\!&+&\!\!\!\!\!\!\ 768r_{00}r_m^mb^6\beta^2+2736r_{00|m}b^mb^2\beta^2 +2880r_{00|m}b^mb^4\beta^2 +768r_{00|m}b^mb^6\beta^2-2736r_{0m|0}b^mb^2\beta^2 \\
\nonumber  \!\!\!\!\!\!&-&\!\!\!\!\!\!\ 2880r_{0m|0}b^mb^4\beta^2-768r_{0m|0}b^mb^6\beta^2+1728s_ms_0^mb^4\beta^2 +1728s_ms_0^mb^2\beta^2+5760s_0r_m^m\beta^3b^2\\
\nonumber  \!\!\!\!\!\!&+&\!\!\!\!\!\!\ 2304s_0r_m^m\beta^3b^4+5760s_{0|m}b^m\beta^3b^2+2304s_{0|m}b^m\beta^3b^4-2880s_{m|0}b^m\beta^3b^2
-1152s_{m|0}b^m\beta^3b^4
\\
\nonumber d_8 \!\!\!\!&:=&\!\!\!\!\ 64r_{00}s_0b^4+112r_{00}s_0b^2-4r_{00}s_0n+320r_{00}s_0b^5 +128r_{00}s_0b^3-16r_{00}s_0b-32r_{00}r_0b^4\\
\nonumber  \!\!\!\!\!\!&-&\!\!\!\!\!\!\  80r_{00}r_0b^2 +8r_{00}r_0n+320r_{00}r_0b^5+128r_{00}r_0b^3-16r_{00}r_0b-12r_{00|0}nb^2+4352s_0^2\beta b^5\\
\nonumber  \!\!\!\!\!\!&-&\!\!\!\!\!\!\ 24r_{00|0}nb^4-16r_{00|0}nb^6+256s_0^2n\beta-3904s_0^2\beta b^2-5440s_0^2\beta b^4-512s_0^2b\beta+2048s_0^2\beta b^3 \\
\nonumber  \!\!\!\!\!\!&-&\!\!\!\!\!\!\ 2816s_0^2\beta b^6+1024s_0^2\beta b^7-192rr_{00}\beta +576r_0^2b^2\beta+384r_0^2b^4\beta-768r_0^2b\beta -2304r_0^2b^3\beta\\
\nonumber  \!\!\!\!\!\!&-&\!\!\!\!\!\!\ 1536r_0^2b^5\beta -672r_0s_0\beta +576r_{0m}r^m_0b^2\beta +864r_{0m}r^m_0b^4\beta+384r_{0m}r^m_0b^6\beta -384r_{0m}s^m_0\beta b^2\\
\nonumber  \!\!\!\!\!\!&-&\!\!\!\!\!\!\ 1152r_{0m}s^m_0\beta b^4-136r_{0m}s^m_0n\beta-640r_{0m}s^m_0\beta b^6+672s_{0|0}\beta b^2+720s_{0|0}\beta b^4-68s_{0|0}n\beta
 \end{eqnarray}
\begin{eqnarray}
\nonumber  \!\!\!\!\!\!&+&\!\!\!\!\!\!\ 256s_{0|0}\beta b^6-1184rs_0\beta^2 -712s_ms^m_0n\beta^2+3504s_ms^m_0\beta^2b^2 +3072s_ms^m_0\beta^2b^4+768s_ms^m_0b^6\beta^2\\
\nonumber  \!\!\!\!\!\!&+&\!\!\!\!\!\!\ 5328s_mr^m_0b^2\beta^2+4608s_mr^m_0\beta^2b^4 +768s_mr^m_0\beta^2b^6 +3552r_ms^m_0b^2\beta^2 +3072r_ms^m_0\beta^2b^4\\
\nonumber  \!\!\!\!\!\!&+&\!\!\!\!\!\!\ 512r_ms^m_0\beta^2b^6-120r_{00}r_m^m\beta-120r_{00|m}b^m\beta+120r_{0m|0}b^m\beta-96s_ms_0^m\beta-1072s_0r_m^m\beta^2\\
\nonumber  \!\!\!\!\!\!&-&\!\!\!\!\!\!\ 1072s_{0|m}b^m\beta^2+536s_{m|0}b^m\beta^2+4288s_{0|m}^m\beta^2b^2+7104s_{0|m}^m\beta^2b^4+4096s_{0|m}^mb^6\beta^2 +448r_0s_0nb\beta\\
\nonumber  \!\!\!\!\!\!&+&\!\!\!\!\!\!\ 512s_{0|m}^m\beta^2b^8+1296s_ms^m\beta^3 +864s_m^is_i^m\beta^3 -112r_{00}s_0nb^4 +16r_{00}s_0nb +64r_{00}s_0nb^3
\\
\nonumber  \!\!\!\!\!\!&+&\!\!\!\!\!\!\ 64r_{00}s_0nb^5-64r_{00}s_0nb^2+32r_{00}r_0nb^4 +16r_{00}r_0nb+64r_{00}r_0nb^3 +64r_{00}r_0nb^5  +448s_0^2nb\beta\\
\nonumber  \!\!\!\!\!\!&+&\!\!\!\!\!\!\ 1408r_0s_0nb^3\beta+1024r_0s_0nb^5\beta+704r_0s_0n\beta b^2+512r_0s_0n\beta b^4+32r_{00}r_0nb^2+1408s_0^2nb^3\beta \\
\nonumber  \!\!\!\!\!\!&+&\!\!\!\!\!\!\ 1024s_0^2nb^5\beta+1184s_0^2n\beta b^2+2912s_0^2n\beta b^4+2560s_0^2n\beta b^6-576rr_{00}b^2\beta-384rr_{00}b^4\beta\\
\nonumber  \!\!\!\!\!\!&+&\!\!\!\!\!\!\  768rr_{00}b\beta+2304rr_{00}b^3\beta +1536rr_{00}b^5\beta+224r_0s_0n\beta -2528r_0s_0\beta b^2-3584r_0s_0\beta b^4\\
\nonumber  \!\!\!\!\!\!&-&\!\!\!\!\!\!\ 1280r_0s_0b\beta -1280r_0s_0\beta b^3+256r_0s_0\beta b^5-1280r_0s_0\beta b^6-672r_{0m}s^m_0n\beta b^2 -68r_{00}s_0\\
\nonumber  \!\!\!\!\!\!&-&\!\!\!\!\!\!\ 1056r_{0m}s^m_0n\beta b^4-512r_{0m}s^m_0nb^6\beta+2r_{00|0}-32r_{00}r_0-2r_{00|0}n-24r_{00|0}b^4-32r_{00|0}b^6
\\
\nonumber  \!\!\!\!\!\!&-&\!\!\!\!\!\!\  1808s_0^2\beta+192r_0^2\beta+120r_{0m}r^m_0\beta+16r_{0m}s^m_0\beta+188s_{0|0}\beta +1296s_ms^m_0\beta^2 +2592s_m^is_i^m\beta^3b^2\\
\nonumber  \!\!\!\!\!\!&+&\!\!\!\!\!\!\ 1608s_mr^m_0\beta^2+1072r_ms^m_0\beta^2+848s_{0|m}^m\beta^2-72{\bf \overline{Ric}}\beta -1152{\bf \overline{Ric}}\beta b^4-384{\bf \overline{Ric}}\beta b^8\\
\nonumber  \!\!\!\!\!\!&-&\!\!\!\!\!\!\ 480{\bf \overline{Ric}}\beta b^2-1152{\bf \overline{Ric}}\beta b^6-336s_{0|0}n\beta b^2-528s_{0|0}n\beta b^4-256s_{0|0}nb^6\beta +2048rs_0b^5\beta^2 \\
\nonumber  \!\!\!\!\!\!&-&\!\!\!\!\!\!\  2048rs_0b^2\beta^2+4736rs_0b\beta^2+8192rs_0b^3\beta^2 -2544s_ms^m_0n\beta^2b^2-2496s_ms^m_0n\beta^2b^4\\
\nonumber  \!\!\!\!\!\!&-&\!\!\!\!\!\!\ 512s_ms^m_0nb^6\beta^2 -512rs_0b^4\beta^2-576r_{00}r_m^mb^2\beta -864r_{00}r_m^mb^4\beta-384r_{00}r_m^mb^6\beta
\\
\nonumber  \!\!\!\!\!\!&-&\!\!\!\!\!\!\ 576r_{00|m}b^mb^2\beta-864r_{00|m}b^mb^4\beta-384r_{00|m}b^mb^6\beta+576r_{0m|0}b^mb^2\beta+864r_{0m|0}b^mb^4\beta \\
\nonumber  \!\!\!\!\!\!&+&\!\!\!\!\!\!\ 384r_{0m|0}b^mb^6\beta+1728s_m^is_i^m\beta^3b^4-576s_ms_0^m\beta b^2-1152s_ms_0^m\beta b^4-768s_ms_0^m\beta b^6
\\
\nonumber  \!\!\!\!\!\!&-&\!\!\!\!\!\!\  3552s_0r_m^mb^2\beta^2-3072s_0r_m^m\beta^2b^4-512s_0r_m^m\beta^2b^6-3552s_{0|m}b^mb^2\beta^2-3072s_{0|m}b^m\beta^2b^4 \\
\nonumber  \!\!\!\!\!\!&-&\!\!\!\!\!\!\ 512s_{0|m}b^m\beta^2b^6+1776s_{m|0}b^mb^2\beta^2+1536s_{m|0}b^m\beta^2b^4+256s_{m|0}b^m\beta^2b^6+1728s_ms^m\beta^3b^2
\\
\nonumber d_9 \!\!\!\!&:=&\!\!\!\!\   4{\bf \overline{Ric}}-128s_0^2nb^2-496s_0^2nb^4-896s_0^2nb^6-32s_0^2nb -128s_0^2nb^3-128s_0^2nb^5-512s_0^2nb^8 \\
\nonumber  \!\!\!\!\!\!&+&\!\!\!\!\!\!\  64rr_{00}b^2+64rr_{00}b^4-64rr_{00}b-256rr_{00}b^3-256rr_{00}b^5-16r_0s_0n+320r_0s_0b^2+704r_0s_0b^4\\
\nonumber  \!\!\!\!\!\!&+&\!\!\!\!\!\!\ 640r_0s_0b^6-384r_0s_0b^5+96r_0s_0b +48r_{0m}s^m_0nb^2+96r_{0m}s^m_0nb^4 +64r_{0m}s^m_0nb^6
\\
\nonumber  \!\!\!\!\!\!&+&\!\!\!\!\!\!\ 24s_{0|0}nb^2+48s_{0|0}nb^4+32s_{0|0}nb^6+320rs_0\beta+120s_ms^m_0n\beta -672s_ms^m_0\beta b^2-576s_ms^m_0\beta b^4 \\
\nonumber  \!\!\!\!\!\!&-&\!\!\!\!\!\!\ 256s_ms^m_0\beta b^6-1440s_mr^m_0\beta b^2-2016s_mr^m_0\beta b^4-768s_mr^m_0\beta b^6-960r_ms^m_0\beta b^2 +64s_ms_0^mb^2\\
\nonumber  \!\!\!\!\!\!&-&\!\!\!\!\!\!\ 1344r_ms^m_0\beta b^4-512r_ms^m_0\beta b^6+48r_{00}r_m^mb^2 +96r_{00}r_m^mb^4+64r_{00}r_m^mb^6 +48r_{00|m}b^mb^2\\
\nonumber  \!\!\!\!\!\!&+&\!\!\!\!\!\!\ 96r_{00|m}b^mb^4+64r_{00|m}b^mb^6-48r_{0m|0}b^mb^2-96r_{0m|0}b^mb^4-64r_{0m|0}b^mb^6 +192s_ms_0^mb^4 \\
\nonumber  \!\!\!\!\!\!&+&\!\!\!\!\!\!\  256s_ms_0^mb^6+128s_ms_0^mb^8+208s_0r_m^m\beta+208s_{0|m}b^m\beta-104s_{m|0}b^m\beta-832s_{0|m}^m\beta b^2\\
\nonumber  \!\!\!\!\!\!&-&\!\!\!\!\!\!\ 1920s_{0|m}^m\beta b^4-1792s_{0|m}^m\beta b^6-512s_{0|m}^m\beta b^8-720s_ms^m\beta^2-312s_m^is_i^m\beta^2 -64r_0s_0nb^2 \\
\nonumber  \!\!\!\!\!\!&-&\!\!\!\!\!\!\  64r_0s_0nb^4-32r_0s_0nb-128r_0s_0nb^3-128r_0s_0nb^5-8r_{0m}r^m_0+16rr_{00}+8r_{00}r_m^m+8r_{00|m}b^m\\
\nonumber  \!\!\!\!\!\!&-&\!\!\!\!\!\!\ 8r_{0m|0}b^m+196s_0^2-16r_0^2-12s_{0|0}-20s_0^2n+512s_0^2b^2 +1120s_0^2b^4+1152s_0^2b^6 -384s_0^2b^3 \\
\nonumber  \!\!\!\!\!\!&-&\!\!\!\!\!\!\  512s_0^2b^7-1152s_0^2b^5+32s_0^2b+512s_0^2b^8-64r_0^2b^2-64r_0^2b^4+64r_0^2b+256r_0^2b^3 +256r_0^2b^5 \\
\nonumber  \!\!\!\!\!\!&+&\!\!\!\!\!\!\ 64r_0s_0-48r_{0m}r^m_0b^2-96r_{0m}r^m_0b^4-64r_{0m}r^m_0b^6+48r_{0m}s^m_0b^2  +192r_{0m}s^m_0b^4+192r_{0m}s^m_0b^6\\
\nonumber  \!\!\!\!\!\!&+&\!\!\!\!\!\!\ 8r_{0m}s^m_0n-48s_{0|0}b^2-48s_{0|0}b^4+4s_{0|0}n-224s_ms^m_0\beta-312s_mr^m_0\beta-208r_ms^m_0\beta \\
\nonumber  \!\!\!\!\!\!&+&\!\!\!\!\!\!\ 8s_ms_0^m-128s_{0|m}^m\beta+96{\bf \overline{Ric}}b^4+64{\bf \overline{Ric}}b^8+32{\bf \overline{Ric}}b^2+128{\bf \overline{Ric}}b^6-1280rs_0b\beta+896rs_0\beta b^2\\
\nonumber  \!\!\!\!\!\!&+&\!\!\!\!\!\!\ 512rs_0\beta b^4-2048rs_0\beta b^5-3584rs_0\beta b^3+576s_ms^m_0n\beta b^2+864s_ms^m_0n\beta b^4+384s_ms^m_0n\beta b^6\\
\nonumber  \!\!\!\!\!\!&-&\!\!\!\!\!\!\ 768s_m^is_i^m\beta^2b^6+960s_0r_m^m\beta b^2+1344s_0r_m^m\beta b^4+512s_0r_m^m\beta b^6+960s_{0|m}b^m\beta b^2\\
\nonumber  \!\!\!\!\!\!&+&\!\!\!\!\!\!\ 1344s_{0|m}b^m\beta b^4+512s_{0|m}b^m\beta b^6-480s_{m|0}b^m\beta b^2 -672s_{m|0}b^m\beta b^4 -256s_{m|0}b^m\beta b^6 \\
\nonumber  \!\!\!\!\!\!&-&\!\!\!\!\!\!\  2016s_ms^mb^2\beta^2-1152s_ms^m\beta^2b^4-1440s_m^is_i^mb^2\beta^2-2016s_m^is_i^m\beta^2b^4
\end{eqnarray}
\begin{eqnarray}
\nonumber d_{10} \!\!\!\!&:=&\!\!\!\!\ 8(1+2b^2)^2(4s_m^is_i^m\beta b^4+4s_{0|m}^mb^4-2s_ms^m_0nb^2 +6s_mr^m_0b^2-4s_0r_m^mb^2+8s_ms^m\beta b^2\\
\nonumber  \!\!\!\!\!\!&-&\!\!\!\!\!\!\ 2s_ms^m_0b^2+16s_m^is_i^m\beta b^2+4r_ms^m_0b^2+2s_{m|0}b^mb^2-4s_{0|m}b^mb^2+4s_{0|m}^mb^2+16rs_0b+7s_m^is_i^m\beta\\
\nonumber  \!\!\!\!\!\!&+&\!\!\!\!\!\!\ 3s_mr^m_0+2r_ms^m_0+s_{m|0}b^m-4rs_0-2s_0r_m^m-2s_{0|m}b^m+22s_ms^m\beta+2s_ms^m_0+s_{0|m}^m-s_ms^m_0n)
\\
\nonumber d_{11} \!\!\!\!&:=&\!\!\!\!\ -4(1+2b^2)^3(2s_m^is_i^mb^2+s_m^is_i^m+4s_ms^m)
\end{eqnarray}

\section{Appendix 3: Coefficients in (\ref{a1})}
\begin{eqnarray*}
d'_0 \!\!\!\!&:=&\!\!\!\!\ 288(8n-11)\beta^7r_{00}^2,\\
d'_2 \!\!\!\!&:=&\!\!\!\!\  -24\beta^5(108{\bf \overline{Ric}}\beta^2-384r_{00}s_0b\beta-96r_{00}s_0n\beta-198r_{00|0}n\beta+192r_{00}s_0nb\beta+192r_{00}r_0nb\beta \\
\nonumber  \!\!\!\!\!\!&+&\!\!\!\!\!\!\ 276r_{00|0}\beta+120r_{00}s_0\beta -384r_{00}r_0b\beta+96r_{00}r_0n\beta+192r_{00|0}\beta b^2-144r_{00}r_0\beta-144r_{00|0}n\beta b^2 \\
 \nonumber  \!\!\!\!\!\!&-&\!\!\!\!\!\!\ 64r_{00}^2nb^4 +64r_{00}^2b^4 +512r_{00}^2b^2 -392r_{00}^2nb^2
 +360r_{00}^2-273r_{00}^2n)
\\
\nonumber d'_4 \!\!\!\!&:=&\!\!\!\!\ -2\beta^3(-5360r_{00}r_0\beta-1778r_{00|0}n\beta+2112r_{00|0}b^4\beta+256r_{00|0}b^6\beta-2144r_{00}s_0\beta +3936r_{00|0}\beta b^2\\
 \nonumber  \!\!\!\!\!\!&-&\!\!\!\!\!\!\  6912s_0^2b\beta^2-576rr_{00}\beta^2-2304r_0^2b\beta^2 +288r_0s_0\beta^2 -1584r_{00}^2nb^4-2712r_{00}^2nb^2\\
  \nonumber  \!\!\!\!\!\!&+&\!\!\!\!\!\!\  1728r_{0m}r^m_0b^2\beta^2  -3888r_{0m}s^m_0n\beta^2+2880r_{0m}s^m_0\beta^2b^2 -1944s_{0|0}n\beta^2 +4032s_{0|0}\beta^2b^2 \\
   \nonumber  \!\!\!\!\!\!&-&\!\!\!\!\!\!\ 2160r_{00}r_m^m\beta^2 -2160r_{00|m}b^m\beta^2
 +2160r_{0m|0}b^m\beta^2 -1384r_{00}s_0n\beta +2304s_0^2nb\beta^2\\
 \nonumber  \!\!\!\!\!\!&-&\!\!\!\!\!\!\  8416r_{00}s_0b\beta+4288r_{00}s_0\beta b^2-1024r_{00}s_0\beta b^5-6400r_{00}s_0\beta b^3+2768r_{00}r_0n\beta \\
 \nonumber  \!\!\!\!\!\!&-&\!\!\!\!\!\!\  8416r_{00}r_0b\beta-1024r_{00}s_0nb^4\beta +5536r_{00}s_0nb\beta +6400r_{00}s_0nb^3\beta +1024r_{00}s_0nb^5\beta \\
  \nonumber  \!\!\!\!\!\!&-&\!\!\!\!\!\!\ 4288r_{00}s_0nb^2\beta+512r_{00}r_0nb^4\beta+5536r_{00}r_0nb\beta +6400r_{00}r_0nb^3\beta +1024r_{00}r_0nb^5\beta \\
   \nonumber  \!\!\!\!\!\!&+&\!\!\!\!\!\!\ 3200r_{00}r_0nb^2\beta+2304r_0s_0nb\beta^2-4352r_{00}r_0\beta b^2-512r_{00}r_0\beta b^4 -1024r_{00}r_0\beta b^5
\\
\nonumber  \!\!\!\!\!\!&-&\!\!\!\!\!\!\ 6400r_{00}r_0\beta b^3-4152r_{00|0}n\beta b^2-2400r_{00|0}n\beta b^4-256r_{00|0}n\beta b^6+1024r_{00}s_0b^4\beta  \\
\nonumber  \!\!\!\!\!\!&+&\!\!\!\!\!\!\ 2304rr_{00}b\beta^2+1152r_0s_0n\beta^2 -8064r_0s_0b\beta^2 -3456r_{0m}s^m_0n\beta^2b^2 +1149r_{00}^2 \\
\nonumber  \!\!\!\!\!\!&+&\!\!\!\!\!\!\ 1440r_{00}^2b^4+4944r_{00}^2b^2 -753r_{00}^2n+2120r_{00|0}\beta+144s_0^2\beta^2 +576r_0^2\beta^2+2160r_{0m}r^m_0\beta^2
\\
\nonumber  \!\!\!\!\!\!&+&\!\!\!\!\!\!\ 3168r_{0m}s^m_0\beta^2+4824s_{0|0}\beta^2 +1296s_{0|m}^m\beta^3 +4428(\textbf{Ric}-{\bf \overline{Ric}})\beta^2 +3456(\textbf{Ric}-{\bf \overline{Ric}})\beta^2b^4 \\
\nonumber  \!\!\!\!\!\!&+&\!\!\!\!\!\!\ 8640(\textbf{Ric}-{\bf \overline{Ric}})\beta^2b^2 -1728s_{0|0}n\beta^2b^2 -1728r_{00}r_m^mb^2\beta^2 -1728r_{00|m}b^mb^2\beta^2-1728r_{0m|0}b^mb^2\beta^2)
\\
\nonumber d'_6 \!\!\!\!&:=&\!\!\!\!\ 2\beta(1128r_{00}r_0\beta+153r_{00|0}n\beta-96r_{00|0}b^4\beta +1386r_{00}s_0\beta -312r_{00|0}\beta b^2-1148s_0^2n\beta^2 \\
\nonumber  \!\!\!\!\!\!&+&\!\!\!\!\!\!\  2976s_0^2\beta^2b^2 -384s_0^2\beta^2b^4 +2048s_0^2b^5\beta^2+6272s_0^2b^3\beta^2 +6368s_0^2b\beta^2+1072rr_{00}\beta^2 \\
\nonumber  \!\!\!\!\!\!&-&\!\!\!\!\!\!\  1408r_0^2b^2\beta^2-256r_0^2b^4\beta^2 +4288r_0^2b\beta^2+5632r_0^2b^3\beta^2 +1024r_0^2b^5\beta^2+2208r_0s_0\beta^2\\
\nonumber  \!\!\!\!\!\!&+&\!\!\!\!\!\!\ 180r_{00}^2nb^4+156r_{00}^2nb^2 -3216r_{0m}r^m_0b^2\beta^2 -2112r_{0m}r^m_0b^4\beta^2 -256r_{0m}r^m_0b^6\beta^2+1728s_0r_m^m\beta^3b^2 \\
\nonumber  \!\!\!\!\!\!&+&\!\!\!\!\!\!\ 1756r_{0m}s^m_0n\beta^2-1008r_{0m}s^m_0\beta^2b^2 -768r_{0m}s^m_0\beta^2b^4 -256r_{0m}s^m_0b^6\beta^2 +878s_{0|0}n\beta^2
\\
\nonumber  \!\!\!\!\!\!&-&\!\!\!\!\!\!\ 5328s_{0|0}\beta^2b^2-3552s_{0|0}\beta^2b^4-512s_{0|0}b^6\beta^2+576rs_0\beta^3 +1512s_ms^m_0n\beta^3 -3168s_ms^m_0\beta^3b^2 \\
\nonumber  \!\!\!\!\!\!&-&\!\!\!\!\!\!\ 2592s_mr^m_0\beta^3b^2-1728r_ms^m_0\beta^3b^2+1220r_{00}r_m^m\beta^2+1220r_{00|m}b^m\beta^2-1220r_{0m|0}b^m\beta^2 \\
\nonumber  \!\!\!\!\!\!&+&\!\!\!\!\!\!\ 432s_ms_0^m\beta^2+1728s_0r_m^m\beta^3+1728s_{0|m}b^m\beta^3-864s_{m|0}b^m\beta^3 -6912s_{0|m}^m\beta^3b^2 -864s_{m|0}b^m\beta^3b^2\\
\nonumber  \!\!\!\!\!\!&-&\!\!\!\!\!\!\ 3456s_{0|m}^m\beta^3b^4+114r_{00}s_0n\beta+936r_{00}s_0b\beta-1200r_{00}s_0\beta b^2 -1152r_{00}s_0\beta b^5 +1728s_{0|m}b^m\beta^3b^2\\
\nonumber  \!\!\!\!\!\!&-&\!\!\!\!\!\!\ 576r_{00}s_0\beta b^3-420r_{00}r_0n\beta+936r_{00}r_0b\beta +1248r_{00}s_0nb^4\beta -840r_{00}s_0nb\beta -2112r_{00}s_0nb^3\beta \\
\nonumber  \!\!\!\!\!\!&-&\!\!\!\!\!\!\ 1152r_{00}s_0nb^5\beta+1392r_{00}s_0nb^2\beta -576r_{00}r_0nb^4\beta -840r_{00}r_0nb\beta -2112r_{00}r_0nb^3\beta \\
\nonumber  \!\!\!\!\!\!&-&\!\!\!\!\!\!\ 1152r_{00}r_0nb^5\beta-1056r_{00}r_0nb^2\beta -3328r_0s_0nb\beta^2 -4864r_0s_0nb^3\beta^2 -1024r_0s_0nb^5\beta^2 \\
\nonumber  \!\!\!\!\!\!&-&\!\!\!\!\!\!\ 2432r_0s_0n\beta^2b^2-+512r_0s_0n\beta^2b^4+1824r_{00}r_0\beta b^2+576r_{00}r_0\beta b^4-1152r_{00}r_0\beta b^5+864s_ms_0^mb^2\beta^2
\end{eqnarray*}
\begin{eqnarray*}
\nonumber  \!\!\!\!\!\!&-&\!\!\!\!\!\!\ 576r_{00}r_0\beta b^3+630r_{00|0}n\beta b^2+792r_{00|0}n\beta b^4+288r_{00|0}n\beta b^6-384r_{00}s_0b^4\beta-3328s_0^2nb\beta^2\\
\nonumber  \!\!\!\!\!\!&-&\!\!\!\!\!\!\ 4864s_0^2nb^3\beta^2-1024s_0^2nb^5\beta^2 -1376s_0^2n\beta^2b^2 +256s_0^2n\beta^2b^4 +1408rr_{00}b^2\beta^2 -256r_{0m|0}b^mb^6\beta^2
\\
\nonumber  \!\!\!\!\!\!&+&\!\!\!\!\!\!\ 256rr_{00}b^4\beta^2-4288rr_{00}b\beta^2 -5632rr_{00}b^3\beta^2 -1024rr_{00}b^5\beta^2 +1664r_0s_0n\beta^2 -2112r_{0m|0}b^mb^4\beta^2 \\
\nonumber  \!\!\!\!\!\!&+&\!\!\!\!\!\!\ 2208r_0s_0\beta^2b^2 -384r_0s_0\beta^2b^4 +2560r_0s_0b^5\beta^2 +11776r_0s_0b^3\beta^2+10144r_0s_0b\beta^2\\
\nonumber  \!\!\!\!\!\!&+&\!\!\!\!\!\!\ 4992r_{0m}s^m_0n\beta^2b^2 +3648r_{0m}s^m_0n\beta^2b^4 -36r_{00}^2-264r_{00}^2b^4 -492r_{00}^2b^2 -3216r_{0m|0}b^mb^2\beta^2\\
\nonumber  \!\!\!\!\!\!&+&\!\!\!\!\!\!\ 24r_{00}^2n-159r_{00|0}\beta +4932s_0^2\beta^2 -1072r_0^2\beta^2 -1220r_{0m}r^m_0\beta^2 -776r_{0m}s^m_0\beta^2  +256r_{00|m}b^mb^6\beta^2 \\
\nonumber  \!\!\!\!\!\!&-&\!\!\!\!\!\!\ 2218s_{0|0}\beta^2-2880s_ms^m_0\beta^3-2592s_mr^m_0\beta^3 -1728r_ms^m_0\beta^3-2916s_{0|m}^m\beta^3 +2112r_{00}r_m^mb^4\beta^2
\\
\nonumber  \!\!\!\!\!\!&-&\!\!\!\!\!\!\  (\textbf{Ric}-{\bf \overline{Ric}})\{1168\beta^2 +6432\beta^2b^4 +256\beta^2b^8 +4880\beta^2b^2+2816\beta^2b^6\} +256r_{00}r_m^mb^6\beta^2
\\
\nonumber  \!\!\!\!\!\!&+&\!\!\!\!\!\!\  512r_{0m}s^m_0n\beta^2b^6+2496s_{0|0}n\beta^2b^2+1824s_{0|0}n\beta^2b^4 +3216r_{00|m}b^mb^2\beta^2-324s_m^is_i^m\beta^4\\
\nonumber  \!\!\!\!\!\!&+&\!\!\!\!\!\!\ 256s_{0|0}n\beta^2b^6-2304rs_0b\beta^3+1728s_ms^m_0n\beta^3b^2 +3216r_{00}r_m^mb^2\beta^2 +2112r_{00|m}b^mb^4\beta^2 )
\\
\nonumber d'_8 \!\!\!\!&:=&\!\!\!\!\ -64r_{00}s_0b^4-112r_{00}s_0b^2+4r_{00}s_0n-320r_{00}s_0b^5 -128r_{00}s_0b^3+16r_{00}s_0b+32r_{00}r_0b^4\\
\nonumber  \!\!\!\!\!\!&+&\!\!\!\!\!\!\  80r_{00}r_0b^2 -8r_{00}r_0n-320r_{00}r_0b^5-128r_{00}r_0b^3+16r_{00}r_0b+12r_{00|0}nb^2-4352s_0^2\beta b^5\\
\nonumber  \!\!\!\!\!\!&+&\!\!\!\!\!\!\ 24r_{00|0}nb^4+16r_{00|0}nb^6-256s_0^2n\beta+3904s_0^2\beta b^2+5440s_0^2\beta b^4+512s_0^2b\beta-2048s_0^2\beta b^3 \\
\nonumber  \!\!\!\!\!\!&+&\!\!\!\!\!\!\ 2816s_0^2\beta b^6-1024s_0^2\beta b^7+192rr_{00}\beta -576r_0^2b^2\beta-384r_0^2b^4\beta+768r_0^2b\beta +2304r_0^2b^3\beta\\
\nonumber  \!\!\!\!\!\!&+&\!\!\!\!\!\!\ 1536r_0^2b^5\beta +672r_0s_0\beta -576r_{0m}r^m_0b^2\beta -864r_{0m}r^m_0b^4\beta-384r_{0m}r^m_0b^6\beta +384r_{0m}s^m_0\beta b^2\\
\nonumber  \!\!\!\!\!\!&+&\!\!\!\!\!\!\ 1152r_{0m}s^m_0\beta b^4+136r_{0m}s^m_0n\beta+640r_{0m}s^m_0\beta b^6-672s_{0|0}\beta b^2-720s_{0|0}\beta b^4+68s_{0|0}n\beta\\
\nonumber  \!\!\!\!\!\!&-&\!\!\!\!\!\!\ 256s_{0|0}\beta b^6+1184rs_0\beta^2 +712s_ms^m_0n\beta^2-3504s_ms^m_0\beta^2b^2 -3072s_ms^m_0\beta^2b^4-768s_ms^m_0b^6\beta^2\\
\nonumber  \!\!\!\!\!\!&-&\!\!\!\!\!\!\ 5328s_mr^m_0b^2\beta^2-4608s_mr^m_0\beta^2b^4 -768s_mr^m_0\beta^2b^6 -3552r_ms^m_0b^2\beta^2 -3072r_ms^m_0\beta^2b^4\\
\nonumber  \!\!\!\!\!\!&-&\!\!\!\!\!\!\ 512r_ms^m_0\beta^2b^6+120r_{00}r_m^m\beta+120r_{00|m}b^m\beta-120r_{0m|0}b^m\beta+96s_ms_0^m\beta+1072s_0r_m^m\beta^2\\
\nonumber  \!\!\!\!\!\!&+&\!\!\!\!\!\!\ 1072s_{0|m}b^m\beta^2-536s_{m|0}b^m\beta^2-4288s_{0|m}^m\beta^2b^2-7104s_{0|m}^m\beta^2b^4-4096s_{0|m}^mb^6\beta^2 -448r_0s_0nb\beta
\\
\nonumber  \!\!\!\!\!\!&-&\!\!\!\!\!\!\ 512s_{0|m}^m\beta^2b^8-1296s_ms^m\beta^3 -864s_m^is_i^m\beta^3 +112r_{00}s_0nb^4 -16r_{00}s_0nb -64r_{00}s_0nb^3 \\
\nonumber  \!\!\!\!\!\!&-&\!\!\!\!\!\!\ 64r_{00}s_0nb^5+64r_{00}s_0nb^2-32r_{00}r_0nb^4 -16r_{00}r_0nb-64r_{00}r_0nb^3 -64r_{00}r_0nb^5  -448s_0^2nb\beta\\
\nonumber  \!\!\!\!\!\!&-&\!\!\!\!\!\!\ 1408r_0s_0nb^3\beta-1024r_0s_0nb^5\beta-704r_0s_0n\beta b^2-512r_0s_0n\beta b^4-32r_{00}r_0nb^2-1408s_0^2nb^3\beta \\
\nonumber  \!\!\!\!\!\!&-&\!\!\!\!\!\!\ 1024s_0^2nb^5\beta-1184s_0^2n\beta b^2-2912s_0^2n\beta b^4-2560s_0^2n\beta b^6+576rr_{00}b^2\beta+384rr_{00}b^4\beta\\
\nonumber  \!\!\!\!\!\!&-&\!\!\!\!\!\!\  768rr_{00}b\beta-2304rr_{00}b^3\beta -1536rr_{00}b^5\beta-224r_0s_0n\beta +2528r_0s_0\beta b^2+3584r_0s_0\beta b^4
\\
\nonumber  \!\!\!\!\!\!&+&\!\!\!\!\!\!\ 1280r_0s_0b\beta+1280r_0s_0\beta b^3-256r_0s_0\beta b^5+1280r_0s_0\beta b^6+672r_{0m}s^m_0n\beta b^2 +68r_{00}s_0\\
\nonumber  \!\!\!\!\!\!&+&\!\!\!\!\!\!\ 1056r_{0m}s^m_0n\beta b^4+512r_{0m}s^m_0nb^6\beta-2r_{00|0}-32r_{00}r_0+2r_{00|0}n+24r_{00|0}b^4+32r_{00|0}b^6
\end{eqnarray*}
\begin{eqnarray*}
\nonumber  \!\!\!\!\!\!&+&\!\!\!\!\!\!\  1808s_0^2\beta-192r_0^2\beta-120r_{0m}r^m_0\beta-16r_{0m}s^m_0\beta-188s_{0|0}\beta -1296s_ms^m_0\beta^2 -2592s_m^is_i^m\beta^3b^2\\
\nonumber  \!\!\!\!\!\!&-&\!\!\!\!\!\!\ 1608s_mr^m_0\beta^2-1072r_ms^m_0\beta^2-848s_{0|m}^m\beta^2-72(\textbf{Ric}-{\bf \overline{Ric}})\beta -1152(\textbf{Ric}-{\bf \overline{Ric}})\beta b^4\\
\nonumber  \!\!\!\!\!\!&-&\!\!\!\!\!\!\ 480(\textbf{Ric}-{\bf \overline{Ric}})\beta b^2-1152(\textbf{Ric}-{\bf \overline{Ric}})\beta b^6+336s_{0|0}n\beta b^2+528s_{0|0}n\beta b^4+256s_{0|0}nb^6\beta   \\
\nonumber  \!\!\!\!\!\!&+&\!\!\!\!\!\!\  2048rs_0b^2\beta^2-4736rs_0b\beta^2-8192rs_0b^3\beta^2 +2544s_ms^m_0n\beta^2b^2+2496s_ms^m_0n\beta^2b^4 -2048rs_0b^5\beta^2\\
\nonumber  \!\!\!\!\!\!&+&\!\!\!\!\!\!\ 512s_ms^m_0nb^6\beta^2+512rs_0b^4\beta^2+576r_{00}r_m^mb^2\beta +864r_{00}r_m^mb^4\beta+384r_{00}r_m^mb^6\beta-384(\textbf{Ric}-{\bf \overline{Ric}})\beta b^8 \\
\nonumber  \!\!\!\!\!\!&+&\!\!\!\!\!\!\ 576r_{00|m}b^mb^2\beta+864r_{00|m}b^mb^4\beta+384r_{00|m}b^mb^6\beta-576r_{0m|0}b^mb^2\beta-864r_{0m|0}b^mb^4\beta \\
\nonumber  \!\!\!\!\!\!&-&\!\!\!\!\!\!\ 384r_{0m|0}b^mb^6\beta-1728s_m^is_i^m\beta^3b^4+576s_ms_0^m\beta b^2+1152s_ms_0^m\beta b^4+768s_ms_0^m\beta b^6 \\
\nonumber  \!\!\!\!\!\!&+&\!\!\!\!\!\!\  3552s_0r_m^mb^2\beta^2+3072s_0r_m^m\beta^2b^4+512s_0r_m^m\beta^2b^6+3552s_{0|m}b^mb^2\beta^2+3072s_{0|m}b^m\beta^2b^4 \\
\nonumber  \!\!\!\!\!\!&+&\!\!\!\!\!\!\ 512s_{0|m}b^m\beta^2b^6-1776s_{m|0}b^mb^2\beta^2-1536s_{m|0}b^m\beta^2b^4-256s_{m|0}b^m\beta^2b^6-1728s_ms^m\beta^3b^2
\\
d'_{10} \!\!\!\!&:=&\!\!\!\!\ -8(1+2b^2)^2(4s_m^is_i^m\beta b^4+4s_{0|m}^mb^4-2s_ms^m_0nb^2 +6s_mr^m_0b^2-4s_0r_m^mb^2+8s_ms^m\beta b^2\\
\nonumber  \!\!\!\!\!\!&-&\!\!\!\!\!\!\ 2s_ms^m_0b^2+16s_m^is_i^m\beta b^2+4r_ms^m_0b^2+2s_{m|0}b^mb^2-4s_{0|m}b^mb^2+4s_{0|m}^mb^2+16rs_0b+7s_m^is_i^m\beta\\
\nonumber  \!\!\!\!\!\!&+&\!\!\!\!\!\!\ 3s_mr^m_0+2r_ms^m_0+s_{m|0}b^m-4rs_0-2s_0r_m^m-2s_{0|m}b^m+22s_ms^m\beta+2s_ms^m_0+s_{0|m}^m-s_ms^m_0n)
\end{eqnarray*}

\section{Appendix 4: Coefficients in (\ref{a2})}
\begin{eqnarray*}
d^{''}_0\!\!\!\!&:=&\!\!\!\!\ -2592({\bf \overline{Ric}} -\textbf{Ric})\beta^7,\\
d^{''}_2\!\!\!\!&:=&\!\!\!\!\ -2\beta^3(-576c^2\beta^4 -144s_0^2\beta^2 +4428({\bf \overline{Ric}} -\textbf{Ric})\beta^2 +3456({\bf \overline{Ric}} -\textbf{Ric})\beta^2b^4 +8640({\bf \overline{Ric}} -\textbf{Ric})\beta^2b^2 \\
\!\!\!\!\!\!&+&\!\!\!\!\!\!\ 1728s_{0|0}nb^2\beta^2 -1152c s_0n\beta^3 +576c s_0\beta^3 +1944s_{0|0}n\beta^2 -4032s_{0|0}b^2\beta^2 -2160c_0\beta^3 -1728c_0\beta^2b^3
\\
\!\!\!\!\!\!&-&\!\!\!\!\!\!\ 4824s_{0|0}\beta^2 -1296s^m_{~0|m}\beta^3) +288\beta^7 c^2(11-8n) -24\beta^5(-192c_0\beta b^2 +96cs_0n\beta +144c_0n\beta b^2
\\
\!\!\!\!\!\!&+&\!\!\!\!\!\!\ 120cs_0\beta +198c_0n\beta -276c_0\beta -96c^2b\beta^2 +144c^2\beta^2)
\\
d^{''}_4\!\!\!\!&:=&\!\!\!\!\ 2\beta(1072c^2\beta^4 +2592cs_0\beta^3 +1728cs_0\beta^3 -1476s_0^2\beta^2  -1168({\bf \overline{Ric}}-\textbf{Ric})\beta^2 -256({\bf \overline{Ric}} -\textbf{Ric})\beta^2b^8
\\
\!\!\!\!\!\!&-&\!\!\!\!\!\!\ 2816({\bf \overline{Ric}}-\textbf{Ric})\beta^2b^6 -6432({\bf \overline{Ric}} -\textbf{Ric})\beta^2b^4 -4880({\bf \overline{Ric}} -\textbf{Ric})\beta^2b^2 +2432c\beta s_0n\beta^2b^2
\\
\!\!\!\!\!\!&-&\!\!\!\!\!\!\ 256s_{0|0}nb^6\beta^2 -2496s_{0|0}nb^2\beta^2 -1824s_{0|0}nb^4\beta^2 -1728s_ms^m_0n\beta^3b^2 -5536s_0^2n\beta^2b^2 -256s_0^2n\beta^2b^4 \\
\!\!\!\!\!\!&+&\!\!\!\!\!\!\ 1664c s_0n\beta^3 -2304c s_0\beta^3b^2 +256c^2b^4\beta^4 -2592c s_0\beta^3 +3552s_{0|0}b^4\beta^2 +512s_{0|0}\beta^2b^6 -878s_{0|0}n\beta^2 \\
\!\!\!\!\!\!&+&\!\!\!\!\!\!\ 5328s_{0|0}b^2\beta^2 +324s_m^is_i^m\beta^4 +1220c_0\beta^3 -432s_ms_0^m\beta^2 -1728s_0nc\beta^3 -1728s^m_{~0|m}b^m\beta^3 \\
\!\!\!\!\!\!&+&\!\!\!\!\!\!\ 864s_{m|0}b^m\beta^3 +3456s^m_{~0|m}\beta^3b^4 +6912s^m_{~0|m}b^2\beta^3 -576cb^2s_0\beta^3 -1512s_ms^m_0n\beta^3 +512c s_0n\beta^3b^4  \\
\!\!\!\!\!\!&+&\!\!\!\!\!\!\ 3168s_ms^m_0\beta^3b^2 +2592cs_0\beta^3b^2 +1728cs_0\beta^3b^2 -2308s_0^2n\beta^2 +3936s_0^2\beta^2b^2 +384s_0^2\beta^2b^4 \\
\!\!\!\!\!\!&+&\!\!\!\!\!\!\ 1408c^2b^2\beta^4 +3216c_0\beta^3b^2 +2112c_0\beta^3b^4 +256c_0 b^6\beta^3 -864s_ms_0^mb^2\beta^2 -1728s_0nc\beta^3b^2 \\
\!\!\!\!\!\!&-&\!\!\!\!\!\!\ 1728s^m_{~0|m}b^m\beta^3b^2+864s_{m|0}b^m\beta^3b^2 +2218s_{0|0}\beta^2+2880s_ms^m_0\beta^3+2916s^m_{~0|m}\beta^3) \\
\!\!\!\!\!\!&-&\!\!\!\!\!\!\ 2\beta^3(-2120c_0\beta -2160c^2\beta^2 +576c^2b^2\beta^2 +1024cs_0n\beta b^4-3200c^2nb^2\beta^2 -512c^2b\beta^2 b^4  \\
\!\!\!\!\!\!&+&\!\!\!\!\!\!\ 4288cs_0nb^2\beta+4152c_0n\beta b^2+2400c_0n\beta b^4 +256c_0n\beta b^6 -1024cs_0b^4\beta -5248cs_0b^2\beta +1384cs_0n\beta \\
\!\!\!\!\!\!&+&\!\!\!\!\!\!\ 4352c^2\beta^2 b^2 +512c^2\beta^2 b^4 -2768c^2b\beta^2 -1728c^2b^2\beta^2 +2160c^2n\beta^2 +2160c_mb^m\beta^2 +3824cs_0\beta  \\
\!\!\!\!\!\!&+&\!\!\!\!\!\!\ 5360c^2\beta^2 +1778c_0n\beta -256c_0b^6\beta-2112c_0b^4\beta -3936c_0\beta b^2 +1728c_mb^m\beta^2b^2 +1728c^2n\beta^2b^2)\\
\!\!\!\!\!\!&-&\!\!\!\!\!\!\  24\beta^5(-64c^2b^4-512c^2b^2-360c^2+392c^2nb^2+273c^2n+64c^2nb^4)
\\
d^{''}_6\!\!\!\!&:=&\!\!\!\!\ -1040s_0^2\beta +192c^2\beta^3 -72({\bf \overline{Ric}}-\textbf{Ric})\beta +480\textbf{Ric}\beta b^2-384({\bf \overline{Ric}}-\textbf{Ric})\beta b^8 -1152({\bf \overline{Ric}}-\textbf{Ric})\beta b^6 \\
\!\!\!\!\!\!&-&\!\!\!\!\!\!\ 1152({\bf \overline{Ric}}-\textbf{Ric})\beta b^4 -480\beta b^2 +704cs_0n\beta^2 b^2 +512c s_0n\beta^2 b^4 -512s_ms^m_0nb^6\beta^2 -256s_{0|0}nb^6\beta
\end{eqnarray*}
\begin{eqnarray*}
\!\!\!\!\!\!&-&\!\!\!\!\!\!\ 336s_{0|0}nb^2\beta -528s_{0|0}nb^4\beta -2048cs_0b^4\beta^2 -512cs_0b^6\beta^2 -2544s_ms^m_0n\beta^2b^2 -2496s_ms^m_0n\beta^2b^4
\\
\!\!\!\!\!\!&-&\!\!\!\!\!\!\ 3424s_0^2n\beta b^2 -6304s_0^2n\beta b^4-3584s_0^2n\beta b^6+224c s_0n\beta^2 -1760c s_0\beta^2 b^2-1664c s_0\beta^2 b^4-512cs_0\beta^2 b^6\\
\!\!\!\!\!\!&+&\!\!\!\!\!\!\  384c^2b^4\beta^3 -672c s_0\beta^2 +672s_{0|0}b^2\beta +256s_{0|0}\beta b^6 -68s_{0|0}n\beta -1072s_0nc\beta^2 +120c_0\beta^2 -96s_ms_0^m\beta \\
\!\!\!\!\!\!&+&\!\!\!\!\!\!\  512s^m_{~0|m}\beta^2b^8+7104s^m_{~0|m}b^4\beta^2 +4096s^m_{~0|m}\beta^2b^6 +4288s^m_{~0|m}b^2\beta^2 +1296s_ms^m\beta^3 \\
\!\!\!\!\!\!&+&\!\!\!\!\!\!\  864s_m^is_i^m\beta^3 -1072s^m_{~0|m}b^m\beta^2 +536s_{m|0}b^m\beta^2+720s_{0|0}b^4\beta -1184cb^2s_0\beta^2 -712s_ms^m_0n\beta^2 \\
\!\!\!\!\!\!&+&\!\!\!\!\!\!\  3504s_ms^m_0\beta^2b^2 +3072s_ms^m_0\beta^2b^4 +768s_ms^m_0b^6\beta^2 +5328cs_0b^2\beta^2 +4608cs_0\beta^2b^4 +768cs_0\beta^2b^6 \\
\!\!\!\!\!\!&+&\!\!\!\!\!\!\  3552cs_0b^2\beta^2 +3072cs_0\beta^2b^4 +512cs_0\beta^2b^6 -512s_0^2n\beta+704s_0^2\beta b^2+3776s_0^2\beta b^4+3328s_0^2\beta b^6\\
\!\!\!\!\!\!&+&\!\!\!\!\!\!\  576c^2b^2\beta^3+576c_0\beta^2b^2 +864c_0\beta^2 b^4+384c_0b^6\beta^2 -576s_ms_0^m\beta b^2-1152s_ms_0^m\beta b^4 -768s_ms_0^m\beta b^6 \\
\!\!\!\!\!\!&-&\!\!\!\!\!\!\  3552s_0ncb^2\beta^2 -3072s_0nc\beta^2b^4 -512s_0nc\beta^2b^6 -3552s^m_{~0|m}b^mb^2\beta^2 -3072s^m_{~0|m}b^m\beta^2b^4 \\
\!\!\!\!\!\!&-&\!\!\!\!\!\!\ 512s^m_{~0|m}b^m\beta^2b^6+1776s_{m|0}b^mb^2\beta^2 +1536s_{m|0}b^m\beta^2b^4 +256s_{m|0}b^m\beta^2b^6 +1728s_ms^m\beta^3b^2
\\
\!\!\!\!\!\!&+&\!\!\!\!\!\!\ 2592s_m^is_i^m\beta^3b^2+1728s_m^is_i^m\beta^3b^4-57+188s_{0|0}\beta +1296s_ms^m_0\beta^2 +848s^m_{~0|m}\beta^2  \\
\!\!\!\!\!\!&-&\!\!\!\!\!\!\ 1248cs_0n\beta b^4+1056c^2nb^2\beta^2+1608cs_0\beta^2 +1072cs_0\beta^2 +2\beta(159c_0\beta+1220c^2\beta^2 -1072c^2b^2\beta^2 \\
\!\!\!\!\!\!&+&\!\!\!\!\!\!\  576c^2b\beta^2 b^4 -1392cs_0nb^2\beta -630c_0n\beta b^2 -792c_0n\beta b^4 -288c_0n\beta b^6 +2304cs_0b^4\beta +2880cs_0b^2\beta \\
\!\!\!\!\!\!&-&\!\!\!\!\!\!\  114cs_0n\beta -1824c^2\beta^2 b^2-576c^2\beta^2 b^4 +420c^2b\beta^2  -1408c^2b^4\beta^2 -256c^2b^6\beta^2 +3216c^2b^2\beta^2
\\
\!\!\!\!\!\!&+&\!\!\!\!\!\!\  2112c^2b^4\beta^2 +256c^2b^6\beta^2 -1220c^2n\beta^2 -1220c_mb^m\beta^2 -1386cs_0\beta -1128c^2\beta^2  -153c_0n\beta  \\
\!\!\!\!\!\!&+&\!\!\!\!\!\!\ 96c_0b^4\beta +312c_0\beta b^2 -3216c_mb^m\beta^2b^2 -2112c_mb^m\beta^2b^4 -256c^2nb^6\beta^2 -256c_mb^mb^6\beta^2 -3216c^2n\beta^2b^2\\
\!\!\!\!\!\!&-&\!\!\!\!\!\!\  2112c^2n\beta^2b^4) -2\beta^3(-1440c^2b^4-4944c^2b^2+753c^2n -1149c^2 +1584c^2nb^4+2712c^2nb^2)
\\
d^{''}_8\!\!\!\!&:=&\!\!\!\!\ 8(1+2b^2)^2(4s_m^is_i^m\beta b^4+4s^m_{~0|m}b^4+16s_m^is_i^m\beta b^2+4s^m_{~0|m}b^2 -4s_0ncb^2 +2s_{m|0}b^mb^2\\
\nonumber \!\!\!\!\!\!&-&\!\!\!\!\!\!\ 2s_ms^m_0nb^2-4s^m_{~0|m}b^mb^2 -2s_ms^m_0b^2+6cs_0b^2+4cs_0b^2 +8s_ms^m\beta b^2 +s^m_{~0|m}+2s_ms^m_0\\
\nonumber \!\!\!\!\!\!&+&\!\!\!\!\!\!\ 3cs_0+2cs_0 -2s^m_{~0|m}b^m +s_{m|0}b^m-4cb^2s_0 -2s_0nc+7s_m^is_i^m\beta +22s_ms^m\beta-s_ms^m_0n)
\\
\nonumber \!\!\!\!\!\!&+&\!\!\!\!\!\!\ 120c^2\beta -120c_mb^m\beta +64cs_0b^4+112cs_0b^2-4cs_0n -12c_0nb^2-24c_0nb^4  -16c_0nb^6 -68cs_0 \\
\nonumber \!\!\!\!\!\!&-&\!\!\!\!\!\!\ 32c^2\beta -2c_0n -32c_0b^6  -24c_0b^4 -112cs_0nb^4-64cs_0nb^2+32c^2\beta nb^4 +32c^2\beta nb^2 -576c^2b^4\beta
\\
\!\!\!\!\!\!&-&\!\!\!\!\!\!\  384c^2b^6\beta +576c^2b^2\beta +864c^2b^4\beta +384c^2b^6\beta -32c^2\beta b^4-80c^2\beta b^2+8c^2\beta n -576c_mb^m\beta b^2
\\
\!\!\!\!\!\!&-&\!\!\!\!\!\!\  864c_mb^m\beta b^4 -384c^2nb^6\beta -384c_mb^mb^6\beta -192c^2b^2\beta-6c^2n\beta b^2-864c^2n\beta b^4 \\
\!\!\!\!\!\!&+&\!\!\!\!\!\!\ 2c_0 -120c^2n\beta +2\beta(264c^2b^4+492c^2b^2-24c^2n +36c^2 -180c^2nb^4-156c^2nb^2)
\end{eqnarray*}

\section{Appendix 5: Coefficients in (\ref{a5})}
\begin{eqnarray*}
t'_0\!\!\!\!&:=&\!\!\!\!\ -288[10 y^iy_jc^2\beta^2 +(12y^iy_jcs_0\beta +17y^iy_jc_0\beta +16y^iy_jc_0b^2\beta -16y^iy_jc^2\beta^2 +4y^iy_jc^2\beta^2) \\
\!\!\!\!\!\!&+&\!\!\!\!\!\!\ (32y^iy_js_0^2b -20y^iy_jcs_0\beta -8y^iy_js_{0|0}b^2 -6y^iy_js_0^2 -10y^iy_js_{0|0} -4y^icy_js_0\beta -8y^icy_jc\beta^2 \\
\!\!\!\!\!\!&-&\!\!\!\!\!\!\ 24y^i c_0y_j b^2\beta -24s_{0|0}^iy_jb^2\beta  +9(R^i_j-\overline{R}^i_j)\beta^2 -2979s^i_0s_{0j}\beta^2 -33y^ic_0y_j\beta -60s_{0|0}^iy_j\beta \\
\!\!\!\!\!\!&+&\!\!\!\!\!\!\ 32y^iy_jcs_0b\beta   +32y^icy_js_0b\beta+ 32y^ic^2y_jb\beta^2 )]\beta^5
\\
t'_2\!\!\!\!&:=&\!\!\!\!\ [-48\beta^5( -48\delta_j^ic^2\beta^2 +18b^iy_jc^2\beta -36y^ib_jc^2\beta +131y^iy_jc^2 +64y^iy_jc^2b^4 +222y^iy_jc^2b^2) \\
\!\!\!\!\!\!&-&\!\!\!\!\!\!\ 4\beta^3(1152\delta_j^ics_0b\beta^3 +1152\delta_j^ic^2b\beta^4 +1270y^iy_jcs_0\beta -1964y^iy_jc^2\beta^2 +2700y^iy_jc_0b^2\beta \\
\!\!\!\!\!\!&+&\!\!\!\!\!\!\ 1920y^iy_jc_0b^4\beta -720y^ib_jcs_0\beta^2 +288y^ib_jc^2\beta^3 -288y^ib_jc_0b^2\beta^2 -864y^is_jcb^2\beta^ 2+956y^iy_jc_0\beta \\
\!\!\!\!\!\!&+&\!\!\!\!\!\!\ 468b^iy_jc_0\beta^2 -36y^ib_jc_0\beta^2 +576y^ib_jc^2\beta^3 +720y^is_jc\beta^3 -576\delta_j^ics_0\beta^3 +960y^ic^2y_j\beta^2 \\
\!\!\!\!\!\!&+&\!\!\!\!\!\!\ 1728y^ic_jb^2\beta^3 -144b^ic^2y_j\beta^3 -1404y^is_jc\beta^2 -2672y^iy_jc^2b^2\beta^2 -3328y^iy_jcs_0b^3\beta +2368y^iy_jcs_0b^2\beta
\\
\!\!\!\!\!\!&-&\!\!\!\!\!\!\ 976y^iy_jcs_0b\beta -976y^iy_jc^2b\beta^2 -3328y^iy_jc^2b^3\beta^2 -1024y^iy_jcs_0b^5\beta +1024y^iy_jcs_0b^4\beta +576\delta_j^ic^2\beta^4
\end{eqnarray*}
\begin{eqnarray*}
\!\!\!\!\!\!&-&\!\!\!\!\!\!\ 864\delta_j^ic_0b^2\beta^3 +2376y^ic_j\beta^3 -1188\delta_j^ic_0\beta^3 -144b^iy_jc^2\beta^3 +288b^iy_jc_0b^2\beta^2 +256y^iy_jc_0b^6\beta \\
\!\!\!\!\!\!&-&\!\!\!\!\!\!\ 2304y^icb_jcb\beta^3  +144b^iy_jcs_0\beta^2 -1024y^iy_jc^2b^5\beta^2 -512y^iy_jc^2b^4\beta^2 -1152b^iy_jcs_0b\beta^2 -1152b^iy_jc^2b\beta^3
\\\
\!\!\!\!\!\!&-&\!\!\!\!\!\!\ 1152y^ib_jcs_0b\beta^2 -1152y^ib_jc^2b\beta^3 +336y^icy_jcb^2\beta^2 -2304y^is_jcb\beta^3 +288c^2y^iy_jb^2\beta^2 +468c^2y^iy_j\beta^2 \\
\!\!\!\!\!\!&+&\!\!\!\!\!\!\ 486\beta^2s^i_0y_jc) +\beta(9216y^icb_js_0b\beta^4 +2976y^iy_js_0^2b^2\beta^2 -7104y^iy_js_0^2b\beta^2 -6144y^iy_js_0^2b^3\beta^2 \\
\!\!\!\!\!\!&+&\!\!\!\!\!\!\ 7152y^iy_jcs_0\beta^3 -1152y^iy_js_ms^m_0b^2\beta^3 +1104y^iy_js_{0|0}b^2\beta^2 +384y^iy_js_{0|0}b^4\beta^2 -4608b^iy_js_{0|0}b^2\beta^3 \\
\!\!\!\!\!\!&+&\!\!\!\!\!\!\ 9216y^ib_js_0^2b\beta^3 -2304y^ib_jcs_0\beta^4 -2304y^ib_js_{0|0}b^2\beta^3 +3456y^is_js_0b^2\beta^3 +1296s^i_ks^k_0y_j\beta^4 +1992y^iy_js_0^2\beta^2 \\
\!\!\!\!\!\!&-&\!\!\!\!\!\!\ 576y^iy_js_ms^m_0\beta^3 +564y^iy_js_{0|0}\beta^2 -4464b^iy_js_{0|0}\beta^3 -864b^iy_js_0^2\beta^3 -2880y^ib_js_{0|0}\beta^3  -2304y^icb_js_0\beta^4\\
\!\!\!\!\!\!&-&\!\!\!\!\!\!\ 1728y^is_js_0\beta^4 +1152y^is_jc\beta^5 -4608\delta_j^is_0^2b\beta^4 -2304\delta_j^ics_0\beta^5 +3456\delta_j^is_{0|0}b^2\beta^4 +3040y^icy_js_0\beta^3 \\
\!\!\!\!\!\!&+&\!\!\!\!\!\!\ 5536y^ic^2y_j\beta^4 +16608y^ic_0y_jb^2\beta^3 +9600y^ic_0y_j\beta^3b^4 +1024y^ic_0y_jb^6\beta^3 -3456y^ics_{0j}b^2\beta^4 \\
\!\!\!\!\!\!&+&\!\!\!\!\!\!\ 6912y^ics_{j0}b^2\beta^4 +3456y^is_{j|0}b^2\beta^4 +23040s_{0|0}^iy_jb^2\beta^3 +18432s_{0|0}^iy_j\beta^3b^4 +3072s_{0|0}^iy_jb^6\beta^3 \\
\!\!\!\!\!\!&-&\!\!\!\!\!\!\ 6912y^is_{0|j}b^2\beta^4  -576b^icy_js_0\beta^4 -1152b^ic^2y_j\beta^5 -3456b^ic_0y_jb^2\beta^4  +5616y^is_js_0\beta^3  -4608\delta_j^ics_0b\beta^5 \\
\!\!\!\!\!\!&-&\!\!\!\!\!\!\ 8856(R^i_j-\overline{R}^i_j)\beta^4 +2931336s^i_0s_{0j}\beta^4 +1296s_{j|0}^i\beta^5 -2592s_{j|0}^i\beta^5 +12096y^iy_jcs_0b^2\beta^3 \\
\!\!\!\!\!\!&-&\!\!\!\!\!\!\ 7104y^iy_jcs_0b\beta^3 -6144y^iy_jcs_0b^3\beta^3 +5719680s^i_0s_{0j}b^2\beta^4 +2287872s^i_0s_{0j}\beta^4b^4 -4320cy^icy_j\beta^4 \\
\!\!\!\!\!\!&-&\!\!\!\!\!\!\ 7776y^is_{0|j}\beta^4 -4320b^ic_0y_j\beta^4 -1296s_{0|0}^ib_j\beta^4 +3888\delta_j^is_{0|0}\beta^4 +7112y^ic_0y_j\beta^3 -3888y^ics_{0j}\beta^4 \\
\!\!\!\!\!\!&+&\!\!\!\!\!\!\ 7776y^ics_{j0}\beta^4 +3888y^is_{j|0}\beta^4 +7620s_{0|0}^iy_j\beta^3  -576b^iy_jcs_0\beta^4 +9216b^iy_js_0^2b\beta^3  +3072y^iy_jcs_0b^4\beta^3 \\
\!\!\!\!\!\!&+&\!\!\!\!\!\!\ 9216b^iy_jcs_0b\beta^4 +9216y^ib_jcs_0b\beta^4 -22144y^icy_js_0b\beta^3 -22144y^ic^2y_jb\beta^4 +5440y^icy_js_0b^2\beta^3 \\
\!\!\!\!\!\!&-&\!\!\!\!\!\!\ 25600y^icy_js_0b^3\beta^3+6400y^ic^2y_jb^2\beta^4 -25600y^ic^2y_jb^3\beta^4 +1024y^icy_js_0\beta^3b^4 -4096y^icy_js_0b^5\beta^3 \\
\!\!\!\!\!\!&+&\!\!\!\!\!\!\ 1024y^ic^2y_jb^4\beta^4 -4096y^ic^2y_jb^5\beta^4 +4608b^icy_js_0b\beta^4 +4608b^ic^2y_jb\beta^5 +4608y^is_js_0b\beta^4 \\
\!\!\!\!\!\!&-&\!\!\!\!\!\!\ 4608y^is_jcb\beta^5 -4608cy^iy_js_0b^2\beta^3 -3456c^2y^iy_jb^2\beta^4 -4464cy^iy_js_0\beta^3  -1944\beta^3s^i_0y_js_0  \\
\!\!\!\!\!\!&+&\!\!\!\!\!\!\ -6912(R^i_j-\overline{R}^i_j)\beta^4b^4  -17280(R^i_j-\overline{R}^i_j)\beta^4b^2) ]
\\
t'_4\!\!\!\!&:=&\!\!\!\!\  [-4\beta^3(-384y^ib_jc^2b^4\beta -648y^ib_jc^2b^2\beta +54b^iy_jc^2\beta -72b^ib_jc^2\beta^2 -894y^ib_jc^2\beta -384\delta_j^ic^2b^4\beta^2 \\
\!\!\!\!\!\!&-&\!\!\!\!\!\!\ 2352\delta_j^ic^2b^2\beta^2 -1638\delta_j^ic^2\beta^2 +504b^iy_jc^2b^2\beta +207y^iy_jc^2 +1128y^iy_jc^2b^4 +1050y^iy_jc^2b^2) \\
\!\!\!\!\!\!&+&\!\!\!\!\!\!\ \beta(-4536\beta^2s^i_0y_jcb^2 -2592\beta^2s^i_0y_jcb^4 +8576\delta_j^ics_0b^2\beta^3 -11072\delta_j^ics_0b\beta^3 -12800\delta_j^ics_0b^3\beta^3 \\
\!\!\!\!\!\!&-&\!\!\!\!\!\!\ 6400\delta_j^ic^2b^2\beta^4 -11072\delta_j^ic^2b\beta^4 +24y^iy_jcs_0\beta +1068y^iy_jc^2\beta^2 -1524y^iy_jc_0b^2\beta -2208y^iy_jc_0b^4\beta  \\
\!\!\!\!\!\!&+&\!\!\!\!\!\!\ 4096y^ic^2b_jb^5\beta^3 +25600y^icb_jb^3\beta^3  -288b^ib_jc^2s_0\beta^3 -576b^ib_jc^2\beta^4 +1152b^ib_jc_0b^2\beta^3 +512y^ib_jc_0b^6\beta^2
\\
\!\!\!\!\!\!&+&\!\!\!\!\!\!\ 6056y^ib_jcs_0\beta^2 -2320y^ib_jc^2\beta^3 +2496y^ib_jc_0b^2\beta^2 +2880y^ib_jc_0b^4\beta^2 +2160y^is_jcb^2\beta^2 +2880g^iy_jc\beta^3b^2 \\
\!\!\!\!\!\!&-&\!\!\!\!\!\!\ 1728b^is_jc\beta^3b^2 -4608b^ic^2b_jb\beta^4 -330y^iy_jc_0\beta -684b^iy_jc_0\beta^2  +1872b^ib_jc_0\beta^3 +268y^ib_jc_0\beta^2 \\
\!\!\!\!\!\!&+&\!\!\!\!\!\!\ 3456\delta_j^icb^2\beta^4 +1152b^ic^2b_j\beta^4 -5536y^ib_jc^2\beta^3 +1944s^i_0b_jc\beta^3 +512\delta_j^ic_0b^6\beta^3 +2768\delta_j^ics_0\beta^3 \\
\!\!\!\!\!\!&+&\!\!\!\!\!\!\ 720b^iy_jc^2\beta^3 -1440b^iy_jc_0b^2\beta^2 -576b^iy_jc_0b^4\beta^2 -960y^iy_jc_0b^6\beta +22144y^ib_jc^2b\beta^3 -6400y^ib_jc^2\beta^3b^2 \\
\!\!\!\!\!\!&-&\!\!\!\!\!\!\ 1024y^icb_jc\beta^3b^4  +1440b^is_jc\beta^4 -6784y^is_jc\beta^3 -16608y^ic_jb^2\beta^3 -9600y^ic_j\beta^3b^4 -1024y^ic_jb^6\beta^3 \\
\!\!\!\!\!\!&+&\!\!\!\!\!\!\ 3456b^ic_jb^2\beta^4 -840y^iy_jc^2\beta^2 +2592b^iy_jc^2\beta^3 +3024y^iy_jc^2b^2\beta^2 +4608y^iy_jcs_0b^3\beta -1584y^iy_jcs_0b^2\beta \\
\!\!\!\!\!\!&+&\!\!\!\!\!\!\ 1200y^iy_jcs_0b\beta +1200y^iy_jc^2b\beta^2 +4608y^iy_jc^2b^3\beta^2 +2052y^is_jc\beta^2 +2304g^iy_jc\beta^3 -12800\delta_j^ic^2b^3\beta^4 \\
\!\!\!\!\!\!&+&\!\!\!\!\!\!\ 2048\delta_j^ics_0b^4\beta^3 -2048\delta_j^ic^2b^5\beta^4 -2048\delta_j^ics_0b^5\beta^3 -1620\beta^2s^i_0y_jc +3840y^iy_jcs_0b^5\beta \\
\!\!\!\!\!\!&-&\!\!\!\!\!\!\ 3264y^iy_jcs_0b^4\beta -5536\delta_j^ic^2\beta^4 +8304\delta_j^ic_0b^2\beta^3 +4800\delta_j^ic_0b^4\beta^3 +3556\delta_j^ic_0\beta^3 -3456b^is_jc\beta^3 \\
\!\!\!\!\!\!&+&\!\!\!\!\!\!\ 4320b^ic_j\beta^4 -7112y^ic_j\beta^3 +1872cy^ib_jc\beta^3 -1440cy^iy_jcb^2\beta^2 -576cy^iy_jc\beta^2b^4 +3816b^iy_jcs_0\beta^2
\\
\!\!\!\!\!\!&+&\!\!\!\!\!\!\ 3840y^iy_jc^2b^5\beta^2 +1920y^iy_jc^2b^4\beta^2 +288b^iy_jcs_0b^2\beta^2 +5760b^iy_jcs_0b\beta^2 +4608b^iy_jcs_0b^3\beta^2 \\
\!\!\!\!\!\!&+&\!\!\!\!\!\!\ 576b^iy_jc^2b^2\beta^3 +5760b^iy_jc^2b\beta^3 +4608b^iy_jc^2b^3\beta^3 -4608b^ib_jcs_0b\beta^3 -4608b^ib_jc^2b\beta^4 \\
\!\!\!\!\!\!&+&\!\!\!\!\!\!\ 2336y^ib_jcs_0b^2\beta^2 +7168y^ib_jcs_0b\beta^2 -512y^ib_jcs_0b^3\beta^2 -4288y^ib_jc^2b^2\beta^3 +7168y^ib_jc^2b\beta^3
\\
\!\!\!\!\!\!&-&\!\!\!\!\!\!\ 512y^ib_jc^2b^3\beta^3 +4320\delta_j^ic\beta^4 -2048y^ib_jcs_0b^5\beta^2 +2048y^ib_jcs_0b^4\beta^2 -2048y^ib_jc^2b^5\beta^3
\end{eqnarray*}
\begin{eqnarray*}
\!\!\!\!\!\!&-&\!\!\!\!\!\!\ 1024y^ib_jc^2b^4\beta^3 -1024\delta_j^ic^2b^4\beta^4 +288y^ic^2y_jb^2\beta^2 +768y^iy_jc^2\beta^2b^4 +1152b^iy_jc^2\beta^3b^2 \\
\!\!\!\!\!\!&-&\!\!\!\!\!\!\ 4608b^is_jcb\beta^4 +22144y^is_jcb\beta^3 -6880y^is_jc\beta^3b^2 -1024y^is_jc\beta^3b^4 +25600y^is_jcb^3\beta^3 +4096y^is_jcb^5\beta^3 \\
\!\!\!\!\!\!&+&\!\!\!\!\!\!\ 1152c^2y^ib_j\beta^3b^2 -684c^2y^iy_j\beta^2)-(648s^i_ks^k_j\beta^5 -773216s^i_0s_{0j}\beta^3 -2916s_{j|0}^i\beta^4 +5832s_{j|0}^i\beta^4 \\
\!\!\!\!\!\!&+&\!\!\!\!\!\!\ 2720y^is_js_0b^2\beta^3 +512y^is_js_0\beta^3b^4+2384y^is_js_0\beta^3  -1864192s^i_0s_{0j}b^6\beta^3 -169472s^i_0s_{0j}\beta^3b^8 \\
\!\!\!\!\!\!&+&\!\!\!\!\!\!\ 2440cy^icy_j\beta^3 -3456b^is_{0|j}\beta^4 -3456\delta_j^is_0\beta^4 -3024\delta_j^is_ms^m_0\beta^4 +1756y^ics_{0j}\beta^3 \\
\!\!\!\!\!\!&+&\!\!\!\!\!\!\ 1728s_{0|0}^ib_j\beta^3 -1728b^ics_{0j}\beta^4 +1000\delta_j^is_0^2\beta^3 -3512y^ics_{j0}\beta^3 -1756y^is_{j|0}\beta^3 -664s_{0|0}^iy_j\beta^2\\
\!\!\!\!\!\!&-&\!\!\!\!\!\!\  1756\delta_j^is_{0|0}\beta^3 +3512y^is_{0|j}\beta^3 +1728b^is_{j|0}\beta^4 +1512y^is_ms^m_0\beta^4 -612y^ic_0y_j\beta^2 +6912s_{j|0}^i\beta^4b^4 \\
\!\!\!\!\!\!&+&\!\!\!\!\!\!\ +2440b^ic_0y_j\beta^3 -6912s_{j|0}^ib^2\beta^4 -3456s_{j|0}^i\beta^4b^4 +3456b^ics_{j0}\beta^4 +13824s_{j|0}^ib^2\beta^4 +1728g^icy_j\beta^4 \\
\!\!\!\!\!\!&-&\!\!\!\!\!\!\ 3648y^is_{j|0}\beta^3b^4 -512y^is_{j|0}b^6\beta^3 -3680s_{0|0}^iy_jb^2\beta^2 +3456b^is_js_0\beta^3 -1476y^is_js_0\beta^2 +1024y^is_{0|j}b^6\beta^3 \\
\!\!\!\!\!\!&+&\!\!\!\!\!\!\ 368b^icy_js_0\beta^3 +564y^ib_js_{0|0}\beta^2 -236y^iy_js_0^2\beta +316y^iy_js_ms^m_0\beta^2 -4y^iy_js_{0|0}\beta -7008s_{0|0}^iy_j\beta^2b^4 \\
\!\!\!\!\!\!&-&\!\!\!\!\!\!\ 5120s_{0|0}^iy_j\beta^2b^6 -1024s_{0|0}^iy_j\beta^2b^8 +1728b^is_{j|0}b^2\beta^4 +5184s_{0|0}^ib_jb^2\beta^3 +3456s_{0|0}^ib_j\beta^3b^4 \\
\!\!\!\!\!\!&-&\!\!\!\!\!\!\ 1024y^ics_{j0}b^6\beta^3 +2144b^iy_jc^2\beta^4 -2736b^ib_js_{0|0}\beta^3 +9984y^is_{0|j}b^2\beta^3 +7296y^is_{0|j}\beta^3b^4 \\
\!\!\!\!\!\!&-&\!\!\!\!\!\!\ 576y^icy_js_0\beta^2 -840y^iy_jc^2\beta^3 -1152b^icb_js_0\beta^4 +3328y^icb_js_0\beta^3 +936y^ib_js_ms^m_0\beta^3 +6432b^ic_0y_jb^2\beta^3 \\
\!\!\!\!\!\!&+&\!\!\!\!\!\!\ 4224b^ic_0y_j\beta^3b^4 +512b^ic_0y_jb^6\beta^3 -576g^iy_js_0\beta^3 -2520y^ic_0y_jb^2\beta^2 -3168y^ic_0y_j\beta^2b^4 -1152y^ic_0y_j\beta^2b^6 \\
\!\!\!\!\!\!&+&\!\!\!\!\!\!\ 4992y^ics_{0j}b^2\beta^3 +3648y^ics_{0j}\beta^3b^4 +512y^ics_{0j}b^6\beta^3 -9984y^ics_{j0}b^2\beta^3 -7296y^ics_{j0}\beta^3b^4 \\
\!\!\!\!\!\!&+&\!\!\!\!\!\!\ 1728y^is_ms^m_0b^2\beta^4 -512\delta_j^is_0^2b^4\beta^3 +2048\delta_j^is_0^2b^5\beta^3 -512\delta_j^is_{0|0}b^6\beta^3 -3456\delta_j^is_ms^m_0b^2\beta^4 \\
\!\!\!\!\!\!&-&\!\!\!\!\!\!\ 704\delta_j^is_0^2b^2\beta^3 +6656\delta_j^is_0^2b\beta^3 +9728\delta_j^is_0^2b^3\beta^3 +3328\delta_j^ic\beta s_0\beta^3 -3456b^is_{0|j}b^2\beta^4 +3456b^ics_{j0}b^2\beta^4 \\
\!\!\!\!\!\!&-&\!\!\!\!\!\!\ 1664y^is_jc\beta^4 -3456\delta_j^is_0b^2\beta^4 -1336b^iy_js_0^2\beta^2 +2088b^iy_js_ms^m_0\beta^3 +952b^iy_js_{0|0}\beta^2 -1392y^ib_js_0^2\beta^2 \\
\!\!\!\!\!\!&+&\!\!\!\!\!\!\ 1728g^icy_jb^2\beta^4 -1728b^ics_{0j}b^2\beta^4 -4992y^is_{j|0}b^2\beta^3 -1328y^iy_jc s_0b^2\beta^2 +128y^iy_jcs_0b\beta^2 \\
\!\!\!\!\!\!&+&\!\!\!\!\!\!\ 128y^iy_jcs_0b^3\beta^2 -256y^iy_jcs_0b^5\beta^2 -1088y^iy_jcs_0b^4\beta^2 +1536y^ib_jcs_0b^4\beta^3 +6912b^ib_jc\beta s_0b\beta^3 \\
\!\!\!\!\!\!&+&\!\!\!\!\!\!\ 4800y^ib_jcs_0b^2\beta^3 -7104y^ib_jcs_0b\beta^3 -6144y^ib_jcs_0b^3\beta^3 -2048b^iy_jcs_0b^5\beta^3 +512b^iy_jcs_0b^4\beta^3 \\
\!\!\!\!\!\!&+&\!\!\!\!\!\!\ -9728b^iy_jcs_0b^3\beta^3 +1280b^iy_jcs_0b^2\beta^3 -6080b^iy_jcs_0b\beta^3 -3456s^i_ks^k_0y_j\beta^3b^4 -1944s^i_0b_js_0\beta^3
\\
\!\!\!\!\!\!&-&\!\!\!\!\!\!\ 4320s^i_ks^k_0y_jb^2\beta^3 +2592\beta^2s^i_0y_js_0b^4 -864b^is_js_0\beta^4 +3216cy^iy_js_0b^2\beta^2 +2880cy^iy_js_0\beta^2b^4 \\
\!\!\!\!\!\!&+&\!\!\!\!\!\!\ 512cy^iy_js_0b^6\beta^2 -2880cy^ib_js_0b^2\beta^3 +2592\beta^2s^i_0y_js_0b^2 -2880y^is_js_0b^2\beta^2 +1728b^is_js_0b^2\beta^3 \\
\!\!\!\!\!\!&+&\!\!\!\!\!\!\ 512b^icy_js_0\beta^3b^4 -2048b^icy_js_0b^5\beta^3 +576y^ib_js_ms^m_0b^2\beta^3 +4864\delta_j^ics_0b^2\beta^4 +6656\delta_j^ics_0b\beta^4
\\
\!\!\!\!\!\!&+&\!\!\!\!\!\!\ 9728\delta_j^ics_0b^3\beta^4 +2048\delta_j^ics_0b^5\beta^4 +48y^iy_js_{0|0}b^4\beta +512b^iy_js_{0|0}\beta^2b^6 -13312y^icb_js_0b\beta^3 \\
\!\!\!\!\!\!&-&\!\!\!\!\!\!\ 256b^iy_js_0^2b^2\beta^2 -6080b^iy_js_0^2b\beta^2 -9728b^iy_js_0^2b^3\beta^2 +3216b^iy_js_{0|0}b^2\beta^2 +2880b^iy_js_{0|0}b^4\beta^2 \\
\!\!\!\!\!\!&+&\!\!\!\!\!\!\ 512b^icy_jc\beta b^4\beta^3 -2048b^icy_jcb^5\beta^4 +4864y^icb_js_0b^2\beta^3 +960y^ib_js_0^2b^2\beta^2 -7104y^ib_js_0^2b\beta^2 \\
\!\!\!\!\!\!&+&\!\!\!\!\!\!\ 2048y^is_jcb^5\beta^4 +64y^iy_js_{0|0}\beta b^6 +1104y^ib_js_{0|0}b^2\beta^2 +384y^ib_js_{0|0}b^4\beta^2 -2048y^is_js_0b^5\beta^3 \\
\!\!\!\!\!\!&-&\!\!\!\!\!\!\ 512y^is_jcb^4\beta^4 -11264b^iy_jc^2b^3\beta^4 +1464y^iy_js_ms^m_0b^2\beta^2 +1920y^iy_js_ms^m_0b^4\beta^2 -256y^iy_js_0^2b^5\beta \\
\!\!\!\!\!\!&-&\!\!\!\!\!\!\ 32y^iy_js_0^2b^4\beta -8576b^icy_js_0b\beta^3 -1728y^icy_js_0\beta^2b^4 +4608y^icy_js_0b^5\beta^2 -1152y^iy_jc^2b^4\beta^3 \\
\!\!\!\!\!\!&+&\!\!\!\!\!\!\ 4608y^iy_jc^2b^5\beta^3 -488y^iy_js_0^2b^2\beta +128y^iy_js_0^2b\beta +128y^iy_js_0^2b^3\beta -392y^iy_jcs_0\beta^2 -8576b^iy_jc^2b\beta^4
\\
\!\!\!\!\!\!&+&\!\!\!\!\!\!\ 1856b^icy_js_0b^2\beta^3 -11264b^icy_js_0b^3\beta^3 +2816b^iy_jc^2b^2\beta^4 +512y^iy_js_ms^m_0b^6\beta^2 +8448y^icy_js_0b^3\beta^2 \\
\!\!\!\!\!\!&+&\!\!\!\!\!\!\ 8448y^iy_jc^2b^3\beta^3 -19456y^icb_js_0b^3\beta^3  +1024y^icb_js_0\beta^3b^4 -4096y^icb_js_0b^5\beta^3 -2880b^ib_js_{0|0}b^2\beta^3\\
\!\!\!\!\!\!&-&\!\!\!\!\!\!\ 6656y^is_js_0b\beta^3 +6656y^is_jcb\beta^4 -9728y^is_js_0b^3\beta^3 -2432y^is_jcb^2\beta^4 +9728y^is_jcb^3\beta^4 -6144y^ib_js_0^2b^3\beta^2 \\
\!\!\!\!\!\!&+&\!\!\!\!\!\!\ 2160y^ib_jcs_0\beta^3 +512b^iy_jcs_0\beta^3 +2880b^iy_js_ms^m_0b^2\beta^3 -1152g^iy_js_0b^2\beta^3 +6912b^ib_js_0^2b\beta^3 \\
\!\!\!\!\!\!&+&\!\!\!\!\!\!\ 2304b^is_js_0b\beta^4 -2304b^is_jcb\beta^5 -2112y^iy_jc^2\beta^3b^2 +1024\delta_j^ics_0b^4\beta^4 +3360y^icy_js_0b\beta^2 +3360y^icy_jc\beta b\beta^2 \\
\!\!\!\!\!\!&-&\!\!\!\!\!\!\ 2448y^icy_js_0b^2\beta^2-2048b^iy_js_0^2b^5\beta^2+512b^iy_js_0^2b^4\beta^2 +4608b^icb_js_0b\beta^4 +1728cy^is_jb^2\beta^4
\\
\!\!\!\!\!\!&+&\!\!\!\!\!\!\ 952cy^iy_js_0\beta^2 -2736cy^ib_js_0\beta^3 +648\beta^2s^i_0y_js_0 +6432c^2y^iy_jb^2\beta^3 +4224c^2y^iy_j\beta^3b^4 +512c^2y^iy_jb^6\beta^3
\end{eqnarray*}
\begin{eqnarray*}
\!\!\!\!\!\!&-&\!\!\!\!\!\!\ 1152y^is_js_0\beta^2b^4 -1296s^i_ks^k_0y_j\beta^3 +1728cy^is_j\beta^4 +648s^i_ks^k_0b_j\beta^4 -3230560s^i_0s_{0j}b^2\beta^3
\\
\!\!\!\!\!\!&-&\!\!\!\!\!\!\ 4257984s^i_0s_{0j}\beta^3b^4 +2336(R^i_j-\overline{R}^i_j)\beta^3 +9760(R^i_j-\overline{R}^i_j)\beta^3b^2 +12864(R^i_j-\overline{R}^i_j)\beta^3b^4\\
\!\!\!\!\!\!&+&\!\!\!\!\!\!\ 5632(R^i_j-\overline{R}^i_j)\beta^3b^6 +512(R^i_j-\overline{R}^i_j)\beta^3b^8)]
\\
t'_6\!\!\!\!&:=&\!\!\!\!\ [\beta(864b^ib_jc^2b^2\beta^2 +576y^ib_jc^2b^4\beta -1512y^ib_jc^2b^2\beta -108b^iy_jc^2\beta  -2376b^ib_jc^2\beta^2 +162y^ib_jc^2\beta \\
\!\!\!\!\!\!&+&\!\!\!\!\!\!\ 5424\delta_j^ic^2b^2\beta^2 +3168\delta_j^ic^2b^4\beta^2 +1506\delta_j^ic^2\beta^2 -27y^iy_jc^2 -504y^iy_jc^2b^4 -324y^iy_jc^2b^2  -648b^iy_jc^2b^2\beta )\\
\!\!\!\!\!\!&-&\!\!\!\!\!\!\ (648s^is_j\beta^3+864s^i_ks^k_j\beta^3 -23832s^i_0s_{0j}\beta -424s_{j|0}^i\beta^2 +848s_{j|0}^i\beta^2 -4s_{0|0}^iy_j -4y^ic_0y_j+120cy^icy_j\beta \\
\!\!\!\!\!\!&-&\!\!\!\!\!\!\ 12y^is_js_0-712\delta_j^is_ms^m_0\beta^2 +68y^ics_{0j}\beta +112s_{0|0}^ib_j\beta -648b^is_ms^m_0\beta^3 -536b^ics_{0j}\beta^2-136y^ics_{j0}\beta \\
\!\!\!\!\!\!&-&\!\!\!\!\!\!\ 68y^is_{j|0}\beta +136y^is_{0|j}\beta +356y^is_ms^m_0\beta^2 -2144s_{j|0}^ib^2\beta^2 -64s_{0|0}^iy_jb^8-128s_{0|0}^iy_jb^6-96s_{0|0}^iy_jb^4\\
\!\!\!\!\!\!&-&\!\!\!\!\!\!\ 32s_{0|0}^iy_jb^2 -8y^icy_js_0-8y^icy_jc\beta -48y^ic_0y_jb^4-24y^ic_0y_jb^2 -32y^ic_0y_jb^6+7104s_{j|0}^i\beta^2b^4 +4y^iy_js_ms^m_0\\
\!\!\!\!\!\!&+&\!\!\!\!\!\!\ 8b^iy_js_{0|0} +40b^iy_js_0^2+4y^ib_js_{0|0}+20y^ib_js_0^2+536b^is_{j|0}\beta^2 -1072\delta_j^is_0\beta^2 +4096s_{j|0}^i\beta^2b^6 +512s_{j|0}^i\beta^2b^8 \\
\!\!\!\!\!\!&+&\!\!\!\!\!\!\ 120b^ic_0y_j\beta +80\delta_j^is_0^2\beta -1072b^is_{0|j}\beta^2 -3552s_{j|0}^i\beta^2b^4 -2048s_{j|0}^i\beta^2b^6-256s_{j|0}^i\beta^2b^8 +4288s_{j|0}^ib^2\beta^2 \\
\!\!\!\!\!\!&+&\!\!\!\!\!\!\ 1072b^ics_{j0}\beta^2 +536s^icy_j\beta^2 -68\delta_j^is_{0|0}\beta+144y^is_js_0b^4 +864b^is_js_0\beta +136y^is_js_0\beta -80s^iy_js_0\beta \\
\!\!\!\!\!\!&+&\!\!\!\!\!\!\ 336y^ics_{0j}b^2\beta +528y^ics_{0j}\beta b^4 +256y^ics_{0j}\beta b^6 -672y^ics_{j0}b^2\beta -1056y^ics_{j0}\beta b^4 -512y^ics_{j0}\beta b^6 \\
\!\!\!\!\!\!&+&\!\!\!\!\!\!\ 32y^iy_js_ms^m_0b^6 +416y^ib_js_0^2b^2 +1248y^is_ms^m_0\beta^2b^4 +256y^is_ms^m_0\beta^2b^6 -256\delta_j^is_{0|0}\beta b^6 \\
\!\!\!\!\!\!&+&\!\!\!\!\!\!\ 1024\delta_j^is_0^2b^5\beta -352\delta_j^is_0^2b^4\beta -512\delta_j^is_ms^m_0b^6\beta^2 +448\delta_j^is_0^2b\beta +1408\delta_j^is_0^2b^3\beta +320y^ib_js_0^2b^4+56y^ib_jc\beta s_0\\
\!\!\!\!\!\!&-&\!\!\!\!\!\!\ 512b^is_{0|j}\beta^2b^6 -864b^is_ms^m_0b^2\beta^3 +1536s^icy_j\beta^2b^4 +512b^ics_{j0}\beta^2b^6 +1776s^icy_jb^2\beta^2 -1776b^ics_{0j}b^2\beta^2 \\
\!\!\!\!\!\!&-&\!\!\!\!\!\!\ 112y^is_jc\beta^2 -128y^ib_js_0^2b^3 -3072\delta_j^is_0\beta^2b^4 +152b^iy_js_ms^m_0\beta +40y^ib_js_ms^m_0\beta +256s^icy_j\beta^2b^6 \\
\!\!\!\!\!\!&-&\!\!\!\!\!\!\ 72s^ib_js_0\beta^2 -1536b^ics_{0j}\beta^2b^4 -256y^is_{j|0}\beta b^6 +1024b^ib_jcs_0b^5\beta^2 -256b^ib_jcs_0b^4\beta^2 +3712b^ib_jcs_0b\beta^2 \\
\!\!\!\!\!\!&+&\!\!\!\!\!\!\ 5632b^ib_jcs_0b^3\beta^2 -256b^ib_jcs_0b^2\beta^2 -256s^i_ks^k_0y_j\beta b^8 -648s^i_0b_js_0\beta -416s^i_ks^k_0y_jb^2\beta -960s^i_ks^k_0y_j\beta b^4 \\
\!\!\!\!\!\!&-&\!\!\!\!\!\!\ 896s^i_ks^k_0y_j\beta b^6 +72(R^i_j-\overline{R}^i_j)\beta +1152(R^i_j-\overline{R}^i_j)\beta b^4 +1152(R^i_j-\overline{R}^i_j)\beta b^6
\\
\!\!\!\!\!\!&+&\!\!\!\!\!\!\ 480(R^i_j-\overline{R}^i_j)\beta b^2 +384(R^i_j-\overline{R}^i_j)\beta b^8 -2592s^i_0b_js_0b^2\beta -2592s^i_0b_js_0\beta b^4 -1440cy^ib_js_0b^2\beta \\
\!\!\!\!\!\!&-&\!\!\!\!\!\!\ 1344cy^ib_js_0\beta b^4 -256cy^ib_js_0b^6\beta -64y^ib_js_{0|0}b^6-160\delta_j^is_0^2b^2\beta -256b^iy_js_0^2b^5+64b^iy_js_{0|0}b^6 \\
\!\!\!\!\!\!&+&\!\!\!\!\!\!\ 16b^iy_jc\beta s_0 +64b^iy_js_0^2b^4 -64b^iy_js_0^2b +112b^iy_js_0^2b^2 -256b^iy_js_0^2b^3+48b^iy_js_{0|0}b^2 +96b^iy_js_{0|0}b^4
\\
\!\!\!\!\!\!&+&\!\!\!\!\!\!\ 32y^icy_js_0b -56y^icy_js_0b^2 +128y^icy_js_0b^3 +32y^iy_jc^2\beta b -32y^ic^2y_jc\beta b^2 +128y^iy_jc^2\beta b^3 -336y^is_{j|0}b^2\beta \\
\!\!\!\!\!\!&-&\!\!\!\!\!\!\ 528y^is_{j|0}\beta b^4 -512\delta_j^is_0b^6\beta^2 -1184b^icb_js_0\beta^2 +192b^iy_jc^2\beta^2 +1056y^is_{0|j}\beta b^4 +512y^is_{0|j}\beta b^6 \\
\!\!\!\!\!\!&-&\!\!\!\!\!\!\ 3552\delta_j^is_0b^2\beta^2+ 24y^iy_js_ms^m_0b^2 +48y^iy_js_ms^m_0b^4 +592b^is_jc\beta^3 +1272y^is_ms^m_0b^2\beta^2 +224\delta_j^ic s_0\beta^2 \\
\!\!\!\!\!\!&-&\!\!\!\!\!\!\ 2544\delta_j^is_ms^m_0b^2\beta^2 -2496\delta_j^is_ms^m_0b^4\beta^2 -336\delta_j^is_{0|0}b^2\beta -528\delta_j^is_{0|0}b^4\beta +224y^icb_js_0\beta -128y^ib_js_0^2b \\
\!\!\!\!\!\!&+&\!\!\!\!\!\!\ 3072b^ics_{j0}\beta^2b^4 -256b^ics_{0j}\beta^2b^6 +3552b^ics_{j0}b^2\beta^2 +384b^ic_0y_j\beta b^6 +1536cy^is_j\beta^2b^4 +256cy^is_j\beta^2b^6 \\
\!\!\!\!\!\!&-&\!\!\!\!\!\!\ 416cy^ib_js_0\beta +1728s^i_ks^k_0b_jb^2\beta^2 +1728s^i_ks^k_0b_j\beta^2b^4 +48cy^iy_js_0b^2+96cy^iy_js_0b^4 +64cy^iy_js_0b^6 \\
\!\!\!\!\!\!&+&\!\!\!\!\!\!\ 576c^2y^iy_jb^2\beta +864c^2y^iy_j\beta b^4 +384c^2y^iy_j\beta b^6 +1776cy^is_jb^2\beta^2 +2592b^is_js_0b^2\beta +1728b^is_js_0\beta b^4 \\
\!\!\!\!\!\!&-&\!\!\!\!\!\!\ 384s^iy_js_0b^2\beta +256y^ib_jc\beta s_0b^5 +320y^ib_jc\beta s_0b^4 -128y^ib_jc\beta s_0b  -192y^ib_js_ms^m_0b^2\beta -672y^ib_js_ms^m_0b^4\beta \\
\!\!\!\!\!\!&+&\!\!\!\!\!\!\ 1664b^ib_js_0^2b^2\beta +3712b^ib_js_0^2b\beta +704\delta_j^ic\beta s_0b^2\beta +448\delta_j^ics_0b\beta^2 +1408\delta_j^ics_0b^3\beta^2 -2048y^icb_js_0b^5\beta \\
\!\!\!\!\!\!&-&\!\!\!\!\!\!\ 256b^ib_js_{0|0}\beta b^6 +1024b^ib_js_0^2b^5\beta -256b^ib_js_0^2b^4\beta +2368b^is_js_0b\beta^2  -2368b^is_jcb\beta^3 +1024b^is_jc\beta^3b^2
\\
\!\!\!\!\!\!&-&\!\!\!\!\!\!\ 576s^iy_js_0\beta b^4 +272y^ib_jc\beta s_0b^2 -128y^ib_jc\beta s_0b^3 +1408y^is_jcb^3\beta^2 -256y^ib_js_ms^m_0b^6\beta -1536b^icy_js_0b^5\beta \\
\!\!\!\!\!\!&+&\!\!\!\!\!\!\ 240b^icy_js_0b^2\beta -768b^icy_js_0b\beta +2304b^icy_js_0b^3\beta +576b^iy_jc^2b^2\beta^2 -768b^iy_jc^2b\beta^2 -2304b^iy_jc^2b^3\beta^2 \\
\!\!\!\!\!\!&+&\!\!\!\!\!\!\ 384b^icy_js_0\beta b^4 +384b^iy_jc^2b^4\beta^2 -1536b^iy_jc^2b^5\beta^2 -2016b^ib_js_ms^m_0b^2\beta^2 +64b^iy_jc\beta s_0b^4 \\
\!\!\!\!\!\!&-&\!\!\!\!\!\!\ 448y^is_js_0b\beta -1408y^is_js_0b^3\beta -352y^is_jcb^2\beta^2 +448y^is_jcb\beta^2 -1440b^ib_js_{0|0}b^2\beta -1344b^ib_js_{0|0}b^4\beta
\end{eqnarray*}
\begin{eqnarray*}
\!\!\!\!\!\!&-&\!\!\!\!\!\!\ 256b^iy_jc\beta s_0b^5 -64b^iy_jc\beta s_0b +64b^iy_jc\beta s_0b^2 -256b^iy_jc\beta s_0b^3 +232y^is_js_0b^2\beta +64y^is_js_0\beta b^4 \\
\!\!\!\!\!\!&-&\!\!\!\!\!\!\ 1024y^is_js_0b^5\beta -256y^is_jcb^4\beta^2 +1024y^is_jcb^5\beta^2 +672b^iy_js_ms^m_0b^2\beta +864b^iy_js_ms^m_0b^4\beta \\
\!\!\!\!\!\!&-&\!\!\!\!\!\!\ 256b^is_js_0\beta^2b^4 +288s^ib_js_0b^2\beta^2 +4736b^icb_js_0b\beta^2 -256s^iy_js_0b^6\beta +5632b^ib_js_0^2b^3\beta \\
\!\!\!\!\!\!&+&\!\!\!\!\!\!\ 80b^ib_jcs_0\beta^2 -4096b^is_jcb^3\beta^3 +1024b^is_js_0b^5\beta^2 +256b^is_jcb^4\beta^3 -1024b^is_jcb^5\beta^3 +4096b^is_js_0b^3\beta^2 \\
\!\!\!\!\!\!&+&\!\!\!\!\!\!\ 1024\delta_j^ics_0b^5\beta^2 +512\delta_j^ics_0b^4\beta^2 +256b^iy_js_ms^m_0b^6\beta +8192b^icb_js_0b^3\beta^2 -512b^icb_js_0\beta^2b^4 +2048b^icb_js_0b^5\beta^2
\\
\!\!\!\!\!\!&+&\!\!\!\!\!\!\ 704y^icb_js_0b^2\beta -896y^icb_js_0b\beta -2816y^icb_js_0b^3\beta +512y^icb_js_0\beta b^4 -2048b^icb_js_0b^2\beta^2 -1408b^is_js_0b^2\beta^2 \\
\!\!\!\!\!\!&-&\!\!\!\!\!\!\ 64s^i_ks^k_0y_j\beta +536cy^is_j\beta^2 +864s^is_jb^2\beta^3 +432s^i_ks^k_0b_j\beta^2 +2592s^i_ks^k_jb^2\beta^3 +1728s^i_ks^k_j\beta^3b^4 \\
\!\!\!\!\!\!&-&\!\!\!\!\!\!\ 158880s^i_0s_{0j}b^2\beta -381312s^i_0s_{0j}\beta b^4-381312s^i_0s_{0j}\beta b^6-127104s^i_0s_{0j}\beta b^8 +8cy^iy_js_0+2480b^ib_js_0^2\beta \\
\!\!\!\!\!\!&-&\!\!\!\!\!\!\ 1440b^ib_js_ms^m_0\beta^2 -416b^ib_js_{0|0}\beta +192y^is_js_0b^6-48b^icy_js_0\beta -48y^ib_js_{0|0}b^4+576b^ic_0y_jb^2\beta \\
\!\!\!\!\!\!&+&\!\!\!\!\!\!\ 864b^ic_0y_j\beta b^4+256y^ib_js_0^2b^5+672y^is_{0|j}b^2\beta -80y^icy_js_0b^4 +128y^icy_js_0b^5 -32y^icy_jc\beta b^4\\
\!\!\!\!\!\!&+&\!\!\!\!\!\!\ 128y^icy_jc\beta b^5 +1536b^is_{j|0}\beta^2b^4 +256b^is_{j|0}\beta^2b^6 -3552b^is_{0|j}b^2\beta^2 -3072b^is_{0|j}\beta^2b^4 -1144b^is_js_0\beta^2 \\
\!\!\!\!\!\!&+&\!\!\!\!\!\!\ 704s_{0|0}^ib_jb^2\beta +1536s_{0|0}^ib_j\beta b^4 +1280s_{0|0}^ib_j\beta b^6+256s_{0|0}^ib_j\beta b^8 +1776b^is_{j|0}b^2\beta^2) \\
\!\!\!\!\!\!&-&\!\!\!\!\!\!\ (-1620s^i_0b_jc\beta^2 +1944y^is_jcb^2\beta^2 +864y^is_jc\beta^2b^4  -3296b^is_jc\beta^3b^2 -512b^is_jc\beta^3b^4 +972y^is_jc\beta^2 \\
\!\!\!\!\!\!&-&\!\!\!\!\!\!\ 2440b^ic_j\beta^3 -306\delta_j^ic_0\beta^2 +2y^iy_jc_0 -2440\delta_j^ic\beta^3 +612y^ic_j\beta^2 +1692b^is_jc\beta^2 -36y^is_jc\beta -684b^ib_jc_0\beta^2\\
\!\!\!\!\!\!&-&\!\!\!\!\!\!\ 8y^iy_jc^2\beta +12b^iy_jc_0\beta +4y^iy_jcs_0 -592g^iy_jc\beta^2 +840\delta_j^ic^2\beta^3 -1260\delta_j^ic_0b^2\beta^2 -1584\delta_j^ic_0b^4\beta^2 \\
\!\!\!\!\!\!&-&\!\!\!\!\!\!\ 2144b^ib_jc^2\beta^3 -708b^iy_jc^2\beta^2 -512\delta_j^ic\beta^3b^6 -6432\delta_j^ic\beta^3b^2 +24y^iy_jc_0b^4 +12y^iy_jc_0b^2 +576b^is_jc\beta^5 \\
\!\!\!\!\!\!&-&\!\!\!\!\!\!\ 4992\delta_j^is_{0|0}b^2\beta^3 -3648\delta_j^is_{0|0}b^4\beta^3 +12y^iy_jc^2\beta +840y^ib_jc^2\beta^2 -6432b^ic_jb^2\beta^3 -4224b^ic_j\beta^3b^4 \\
\!\!\!\!\!\!&-&\!\!\!\!\!\!\ 512b^ic_jb^6\beta^3 +2520y^ic_jb^2\beta^2 +3168y^ic_j\beta^2b^4 +1152y^ic_j\beta^2b^6 +16y^iy_jc_0b^6 -576\delta_j^ic_0b^6\beta^2 -228\delta_j^ics_0\beta^2
\\
\!\!\!\!\!\!&-&\!\!\!\!\!\!\ 24y^ib_jc_0\beta -4224\delta_j^ic\beta^3b^4 +576g^ib_jc\beta^3 -24b^iy_jcs_0b^2\beta -192b^iy_jcs_0b\beta -384b^iy_jcs_0b^3\beta -48b^iy_jc^2b^2\beta^2 \\
\!\!\!\!\!\!&+&\!\!\!\!\!\!\ 768y^ib_jc^2b^4\beta^2 -480y^ib_jcs_0b\beta +384y^ib_jcs_0b^3\beta +912y^ib_jc^2b^2\beta^2 -480y^ib_jc^2b\beta^2 +384y^ib_jc^2b^3\beta^2 \\
\!\!\!\!\!\!&+&\!\!\!\!\!\!\ 1704y^ib_jcs_0b^2\beta +1536y^ib_jcs_0b^5\beta +96y^ib_jcs_0b^4\beta +1536y^ib_jc^2b^5\beta^2 +4608b^ib_jc^2b^3\beta^3 +576b^ib_jc^2b^2\beta^3 \\
\!\!\!\!\!\!&+&\!\!\!\!\!\!\ 5760b^ib_jc^2b\beta^3 +5760b^ib_jcs_0b\beta^2 +4608b^ib_jcs_0b^3\beta^2 +1728b^ib_jcs_0b^2\beta^2 -288b^iy_jcs_0b^4\beta -192b^iy_jc^2b\beta^2\\
\!\!\!\!\!\!&-&\!\!\!\!\!\!\ -384b^iy_jc^2b^3\beta^2 +432\beta s^i_0y_jcb^4 +288\beta s^i_0y_jcb^6 -3032b^is_jc\beta^3 -4536s^i_0b_jcb^2\beta^2 -2592s^i_0b_jc\beta^2b^4 \\
\!\!\!\!\!\!&+&\!\!\!\!\!\!\ 48cy^iy_jcb^2\beta +48cy^iy_jc\beta b^4 -1440cy^ib_jcb^2\beta^2 -576cy^ib_jc\beta^2b^4 +216\beta s^i_0y_jcb^2 +216y^is_jcb^2\beta
\\
\!\!\!\!\!\!&+&\!\!\!\!\!\!\ 1152b^is_jc\beta^2b^4 +3312b^is_jcb^2\beta^2 -2208g^iy_jcb^2\beta^2 -512g^iy_jcb^6\beta^2 -2496\delta_j^ics_0b^4\beta^2 +2304\delta_j^ic^2b^5\beta^3 \\
\!\!\!\!\!\!&+&\!\!\!\!\!\!\ 576y^is_jc\beta b^6 -576b^ib_jc_0b^4\beta^2 +2304\delta_j^ics_0b^5\beta^2 +6480b^ib_jcs_0\beta^2 +864y^is_jc\beta b^4+8576b^is_jcb\beta^3
\\
\!\!\!\!\!\!&+&\!\!\!\!\!\!\ 11264b^is_jcb^3\beta^3 +2048b^is_jcb^5\beta^3 +720b^ib_jc^2\beta^3 -1440b^ib_jc_0b^2\beta^2 -2304g^iy_jc\beta^2b^4 +228y^ib_jc^2\beta^2
\\
\!\!\!\!\!\!&-&\!\!\!\!\!\!\ 264y^ib_jc_0b^2\beta -24b^iy_jc^2\beta^2 +184y^iy_jcs_0b^4 -64y^iy_jcs_0b^5-32y^iy_jc^2\beta b^4 -64y^iy_jc^2\beta b^5 -624y^ib_jc_0b^4\beta \\
\!\!\!\!\!\!&+&\!\!\!\!\!\!\ 768b^icy_jcb^2\beta^2 +2112\delta_j^ic^2\beta b^2\beta^2 +1680\delta_j^ic^2\beta b\beta^2 -4608y^is_jcb^5\beta^2 +4224\delta_j^ic^2b^3\beta^3 \\
\!\!\!\!\!\!&-&\!\!\!\!\!\!\ 2784\delta_j^ics_0b^2\beta^2 +1680\delta_j^ics_0b\beta^2 +4224\delta_j^ics_0b^3\beta^2 -3360y^is_jcb\beta^2 -384y^ib_jc_0b^6\beta -16y^iy_jcs_0b \\
\!\!\!\!\!\!&+&\!\!\!\!\!\!\ 100y^iy_jcs_0b^2 -64y^iy_jcs_0b^3 -16y^iy_jc^2\beta b -32y^iy_jc^2\beta b^2 -64y^iy_jc^2\beta b^3 -144y^ib_jcs_0\beta +48b^iy_jc_0b^2\beta \\
\!\!\!\!\!\!&+&\!\!\!\!\!\!\ 48b^iy_jc_0b^4\beta -4608y^ib_jc^2b^5\beta^2 +1152g^ib_jc\beta^3b^2 +8576b^ib_jc^2b\beta^3 -8448y^is_jcb^3\beta^2 +1152\delta_j^ic^2b^4\beta^3 \\
\!\!\!\!\!\!&-&\!\!\!\!\!\!\ 108y^iy_jc^2b^2\beta -192y^iy_jc^2\beta b^4 +2048b^ib_jc^2b^5\beta^3 +2112y^ib_jc^2b^2\beta^2 -3360y^ib_jc^2b\beta^ 2-8448y^ib_jc^2b^3\beta^2 \\
\!\!\!\!\!\!&+&\!\!\!\!\!\!\ 1152y^ib_jc^2\beta^2b^4 -48b^iy_jcs_0\beta -2816b^ib_jc^2\beta^3b^2 -512b^ib_jc^2\beta^3b^4 +11264b^ib_jc^2b^3\beta^3 +12y^iy_jc^2\beta \\
\!\!\!\!\!\!&-&\!\!\!\!\!\!\ 684y^ib_jc^2\beta^2 +36\beta s^i_0y_jc)]
\\
t'_8\!\!\!\!&:=&\!\!\!\!\ [(4(2b^2+1))(-2s^i_ks^k_0b_j-22s^is_j\beta -14s^i_ks^k_j\beta -2cy^is_j-2b^is_{0|j}+4b^is_{0|j}+2b^ics_{0j}+2\delta_j^is_ms^m_0\\
\!\!\!\!\!\!&+&\!\!\!\!\!\!\ 12s_{0|j}^ib^4+6s_{0|j}^ib^2+8s_{0|j}^ib^6-24s_{j0}^ib^4-y^is_ms^m_0-2s^icy_j -4b^ics_{j0} +4\delta_j^is_0-16s_{j0}^ib^6\\
\!\!\!\!\!\!&-&\!\!\!\!\!\!\ 12s_{j0}^ib^2+8b^ics_{0j}b^4+8b^ics_{0j}b^2-16b^ics_{j0}b^4 -16b^ics_{j0}b^2-8s^icy_jb^4 -8s^icy_jb^2-8b^is_{0|j}b^2\\
\!\!\!\!\!\!&+&\!\!\!\!\!\!\ 22b^is_ms^m_0\beta+10b^is_js_0-8b^is_{0|j}b^4+2s^ib_js_0 +16\delta_j^is_0b^2+16\delta_j^is_0b^4+8b^icb_js_0
\end{eqnarray*}
\begin{eqnarray*}
\!\!\!\!\!\!&+&\!\!\!\!\!\!\ 8\delta_j^is_ms^m_0b^2 +8\delta_j^is_ms^m_0b^4+16b^is_{0|j}b^4-4b^is_jc\beta-4y^is_ms^m_0b^4 -4y^is_ms^m_0b^2+16b^ib_js_ms^m_0\\
\!\!\!\!\!\!&+&\!\!\!\!\!\!\ 16b^is_{0|j}b^2-32b^icb_js_0b-16b^is_js_0b-32b^is_js_0b^3 +40b^ib_js_ms^m_0b^2 +16b^ib_js_ms^m_0b^4 +8b^is_js_0b^2 \\
\!\!\!\!\!\!&+&\!\!\!\!\!\!\ 16b^is_jc\beta b-8b^is_jc\beta b^2+32b^is_jc\beta b^3+16b^icb_js_0b^2+52b^is_ms^m_0b^2\beta +16b^is_ms^m_0\beta b^4-16s^ib_js_0b^4\\
\!\!\!\!\!\!&-&\!\!\!\!\!\!\ 64b^icb_js_0b^3-2s_{j0}^i+s_{0|j}^i-4s^ib_js_0b^2-8cy^is_jb^4-52s^is_jb^2\beta -16s^is_j\beta b^4-60s^i_ks^k_jb^2\beta-72s^i_ks^k_j\beta b^4\\
\!\!\!\!\!\!&-&\!\!\!\!\!\!\ 16s^i_ks^k_j\beta b^6-8cy^is_jb^2-24s^i_ks^k_0b_jb^4-12s^i_ks^k_0b_jb^2-16s^i_ks^k_0b_jb^6) -(4y^ic_j-36s^i_0b_jc +36b^is_jc-2\delta_j^ic_0 \\
\!\!\!\!\!\!&-&\!\!\!\!\!\!\ 120b^ic_j\beta -24\delta_j^ic_0b^4 -12\delta_j^ic_0b^2 -4\delta_j^ics_0 +8\delta_j^ic^2\beta -16\delta_j^ic_0b^6 -312b^is_jc\beta -128y^is_jcb^3 \\
\!\!\!\!\!\!&-&\!\!\!\!\!\!\ 128y^is_jcb^5 -112\delta_j^ics_0b^4 +64\delta_j^ics_0b^5 +32\delta_j^ic^2\beta b^4 +64\delta_j^ic^2\beta b^5 +20y^is_jcb^2 +16\delta_j^ics_0b -64\delta_j^ics_0b^2 \\
\!\!\!\!\!\!&+&\!\!\!\!\!\!\ 64\delta_j^ics_0b^3 +8y^is_jcb^4-48s^iy_jcb^2 -120\delta_j^ic\beta -8s^iy_jc+8y^icb_jc +32y^ic_jb^6 -12b^icy_jc+8y^is_jc  \\
\!\!\!\!\!\!&-&\!\!\!\!\!\!\ 12b^ib_jc_0 +24y^ic_jb^2+144b^is_jcb^2 -384b^ic_j\beta b^6 +144b^is_jcb^4-576b^ic_jb^2\beta -864b^ic_j\beta b^4 -32y^ic^2b_jb\\
\!\!\!\!\!\!&+&\!\!\!\!\!\!\ 32y^ic^2b_jb^2-128y^ic^2b_jb^3+32y^ic^2b_jb^4-128y^ib_jc^2b^5 -32y^is_jcb +408b^ib_jcs_0 +24b^ib_jc^2\beta -48b^ib_jc_0b^4 \\
\!\!\!\!\!\!&-&\!\!\!\!\!\!\ 48b^ib_jc_0b^2 -192b^icb_jc\beta -64s^iy_jcb^6 -96s^iy_jcb^4  +56s^ib_jc\beta -576\delta_j^icb^2\beta -864\delta_j^ic\beta b^4-384\delta_j^ic\beta b^6 \\
\!\!\!\!\!\!&-&\!\!\!\!\!\!\ 216s^i_0b_jcb^2-432s^i_0b_jcb^4-288s^i_0b_jcb^6 -48cy^ib_jcb^2-48cy^ib_jcb^4 +288b^ib_jcs_0b^4+192b^ib_jcs_0b\\
\!\!\!\!\!\!&+&\!\!\!\!\!\!\ 528b^ib_jcs_0b^2+384b^ib_jcs_0b^3 +768b^is_jcb\beta+192b^ib_jc^2\beta b +48b^ib_jc^2\beta b^2+384b^ib_jc^2\beta b^3\\
\!\!\!\!\!\!&+&\!\!\!\!\!\!\ 2304b^is_jcb^3\beta +1536b^is_jcb^5\beta +432s^ib_jcb^2\beta +768s^ib_jc\beta b^4+256s^ib_jc\beta b^6 +768b^ib_jc^2b\beta\\
\!\!\!\!\!\!&+&\!\!\!\!\!\!\ 2304b^ib_jc^2b^3\beta -744b^is_jcb^2\beta -384b^icb_jc\beta b^4 -576b^icb_jcb^2\beta +1536b^ib_jc^2b^5\beta -384b^is_jc\beta b^4 \\
\!\!\!\!\!\!&+&\!\!\!\!\!\!\ -12y^ib_jc^2-24b^iy_jc^2b^2 +16\delta_j^ic^2\beta b+32\delta_j^ic^2\beta b^2 +64\delta_j^ic^2\beta b^3) -(-48\delta_j^ic^2\beta +12b^iy_jc^2 +9y^ib_jc^2 \\
\!\!\!\!\!\!&+&\!\!\!\!\!\!\ 24b^iy_jc^2b^2 +600b^ib_jc^2\beta +48y^ib_jc^2b^4+96y^ib_jc^2b^2 -360\delta_j^ic^2b^4\beta -312\delta_j^ic^2b^2\beta +120b^ib_jc^2b^2\beta)]
\end{eqnarray*}

\section{Appendix 6: Coefficients in (\ref{s3})}
\begin{eqnarray*}
A_0\!\!\!\!&:=&\!\!\!\!\ -288\beta^9r_{00}^2(-11+8n),
\\
A_1\!\!\!\!&:=&\!\!\!\!\ 48\beta^8(-54r_{00|0}\beta-104r_{00}^2b^2-295r_{00}^2+80r_{00}^2nb^2+220r_{00}^2n+36r_{00|0}n\beta)
\\
A_2\!\!\!\!&:=&\!\!\!\!\ 24\beta^7(108(\textbf{Ric}-{\bf \overline{Ric}})\beta^2-144r_{00|0}nb^2\beta+192r_{00}s_0nb\beta -384r_{00}r_0b\beta+192r_{00}r_0nb\beta-96r_{00}s_0n\beta \\
\!\!\!\!\!\!&-&\!\!\!\!\!\!\ 342r_{00|0}n\beta+192r_{00|0}\beta b^2+120r_{00}s_0\beta-144r_{00}r_0\beta +96r_{00}r_0n\beta-384r_{00}s_0b\beta +492r_{00|0}\beta \\
\!\!\!\!\!\!&+&\!\!\!\!\!\!\ 1144r_{00}^2+928r_{00}^2b^2+64r_{00}^2b^4-865r_{00}^2n-64r_{00}^2nb^4-712r_{00}^2nb^2)
\\
A_3\!\!\!\!&:=&\!\!\!\!\ -12\beta^6(162\beta^2\theta-384r_{00}s_0nb^2\beta +1696r_{00}s_0nb\beta+512r_{00}s_0nb^3\beta +256r_{00}r_0nb^2\beta+512r_{00}s_0b^2\beta\\
\!\!\!\!\!\!&-&\!\!\!\!\!\!\ 720r_{00}s_0n\beta-768r_{00}s_0b^3\beta-3168r_{00}s_0b\beta -320r_{00}r_0b^2\beta+848r_{00}r_0n\beta-768r_{00}r_0b^3\beta\\
\!\!\!\!\!\!&-&\!\!\!\!\!\!\ 3168r_{00}r_0b\beta-1272r_{00|0}nb^2\beta-192r_{00|0}nb^4\beta+1044(\textbf{Ric}-{\bf \overline{Ric}})\beta^2 +1696r_{00}r_0nb\beta\\
\!\!\!\!\!\!&+&\!\!\!\!\!\!\ 512r_{00}r_0nb^3\beta +576(\textbf{Ric}-{\bf \overline{Ric}})\beta^2b^2 +512r_{00}^2b^4 +3648r_{00}^2b^2 -1923r_{00}^2n+1916r_{00|0}\beta \\
\!\!\!\!\!\!&+&\!\!\!\!\!\!\ 144r_{0m}r^m_0\beta^2+288r_{0m}s^m_0\beta^2+360s_{0|0}\beta^2+224r_{00|0}\beta b^4+592r_{00}s_0\beta-1360r_{00}r_0\beta \\
\!\!\!\!\!\!&-&\!\!\!\!\!\!\ 1398r_{00|0}n\beta+1568r_{00|0}\beta b^2-2696r_{00}^2nb^2 -544r_{00}^2nb^4 -288r_{0m}s^m_0n\beta^2-144s_{0|0}n\beta^2 \\
\!\!\!\!\!\!&+&\!\!\!\!\!\!\ 2560r_{00}^2-144r_{00}r_m^m\beta^2-144r_{00|m}b^m\beta^2+144r_{0m|0}b^m\beta^2-162n\beta^2\theta)
\\
A_4\!\!\!\!&:=&\!\!\!\!\ 2\beta^5(3726\beta^2\theta-8896r_{00}s_0nb^2\beta-1024r_{00}s_0n\beta b^4+18976r_{00}s_0nb\beta +12544r_{00}s_0nb^3\beta \\
\!\!\!\!\!\!&+&\!\!\!\!\!\!\ 1024r_{00}s_0n\beta b^5+6272r_{00}r_0nb^2\beta+512r_{00}r_0n\beta b^4+1024r_{00}s_0b^4\beta +10432r_{00}s_0b^2\beta \\
\!\!\!\!\!\!&-&\!\!\!\!\!\!\ 6568r_{00}s_0n\beta -1024r_{00}s_0b^5\beta-15616r_{00}s_0b^3\beta -32608r_{00}s_0b\beta -512r_{00}r_0b^4\beta-8192r_{00}r_0b^2\beta\\
\!\!\!\!\!\!&+&\!\!\!\!\!\!\ 9488r_{00}r_0n\beta-1024r_{00}r_0b^5\beta-15616r_{00}r_0b^3\beta-32608r_{00}r_0b\beta-14232r_{00|0}nb^2\beta -4704r_{00|0}nb^4\beta \\
\!\!\!\!\!\!&-&\!\!\!\!\!\!\ 256r_{00|0}nb^6\beta+2304s_0^2nb\beta^2 +2304rr_{00}b\beta^2 +1152r_0s_0n\beta^2 -8064r_0s_0b\beta^2 -3456r_{0m}s^m_0nb^2\beta^2 \\
\!\!\!\!\!\!&-&\!\!\!\!\!\!\ 1728s_{0|0}nb^2\beta^2-2592n\beta^2b^2\theta-1728r_{00}r_m^mb^2\beta^2-1728r_{00|m}b^mb^2\beta^2
+1728r_{0m|0}b^mb^2\beta^2
\end{eqnarray*}
\begin{eqnarray*}
\!\!\!\!\!\!&+&\!\!\!\!\!\!\  13068(\textbf{Ric}-{\bf \overline{Ric}})\beta^2+18976r_{00}r_0nb\beta +12544r_{00}r_0nb^3\beta +1024r_{00}r_0n\beta b^5+2304r_0s_0nb\beta^2 \\
\!\!\!\!\!\!&+&\!\!\!\!\!\!\  15552(\textbf{Ric}-{\bf \overline{Ric}})\beta^2b^2+3456(\textbf{Ric}-{\bf \overline{Ric}})\beta^2b^4 +324\beta^2\sigma+5280r_{00}^2b^4 +25296r_{00}^2b^2-8049r_{00}^2n \\
\!\!\!\!\!\!&+&\!\!\!\!\!\!\ 12584r_{00|0}\beta+144s_0^2\beta^2 +576r_0^2\beta^2+3888r_{0m}r^m_0\beta^2 +6624r_{0m}s^m_0\beta^2 +9144s_{0|0}\beta^2 \\
\!\!\!\!\!\!&+&\!\!\!\!\!\!\ 1296s_{0|m}^m\beta^3+4800r_{00|0}\beta b^4+256r_{00|0}b^6\beta+640r_{00}s_0\beta -16496r_{00}r_0\beta-9698r_{00|0}n\beta \\
\!\!\!\!\!\!&+&\!\!\!\!\!\!\ 15840r_{00|0}\beta b^2-6912s_0^2b\beta^2-576rr_{00}\beta^2 -2304r_0^2b\beta^2+288r_0s_0\beta^2 -17112r_{00}^2nb^2-5808r_{00}^2nb^4 \\
\!\!\!\!\!\!&+&\!\!\!\!\!\!\ 1728r_{0m}r^m_0b^2\beta^2+2880r_{0m}s^m_0b^2\beta^2-7344r_{0m}s^m_0n\beta^2+4032s_{0|0}b^2\beta^2
-3672s_{0|0}n\beta^2+11085r_{00}^2 \\
\!\!\!\!\!\!&-&\!\!\!\!\!\!\ 3888r_{00}r_m^m\beta^2-3888r_{00|m}b^m\beta^2+3888r_{0m|0}b^m\beta^2-324n\beta^2\sigma-3726n\beta^2\theta
+2592\beta^2b^2\theta)
\\
A_5\!\!\!\!&:=&\!\!\!\!\ -3\beta^4(4050\beta^2\theta-9472r_{00}s_0nb^2\beta-2560r_{00}s_0n\beta b^4+13024r_{00}s_0nb\beta +14080r_{00}s_0nb^3\beta \\
\!\!\!\!\!\!&+&\!\!\!\!\!\!\ 2560r_{00}s_0n\beta b^5+7040r_{00}r_0nb^2\beta+1280r_{00}r_0n\beta b^4+2048r_{00}s_0b^4\beta +9536r_{00}s_0b^2\beta \\
\!\!\!\!\!\!&-&\!\!\!\!\!\!\  3520r_{00}s_0n\beta-1536r_{00}s_0b^5\beta-13056r_{00}s_0b^3\beta -20256r_{00}s_0b\beta -1280r_{00}r_0b^4\beta-9728r_{00}r_0b^2\beta \\
\!\!\!\!\!\!&+&\!\!\!\!\!\!\ 6512r_{00}r_0n\beta-1536r_{00}r_0b^5\beta -13056r_{00}r_0b^3\beta-20256r_{00}r_0b\beta-9768r_{00|0}nb^2\beta-5280r_{00|0}nb^4\beta \\
\!\!\!\!\!\!&-&\!\!\!\!\!\!\ 640r_{00|0}nb^6\beta-256s_0^2n\beta^2b^2 +6016s_0^2nb\beta^2+2048s_0^2nb^3\beta^2 -512rr_{00}b^2\beta^2 +6400rr_{00}b\beta^2 \\
\!\!\!\!\!\!&+&\!\!\!\!\!\!\ 2048rr_{00}b^3\beta^2+3008r_0s_0n\beta^2+512r_0s_0\beta^2b^2-20352r_0s_0b\beta^2 -6144r_0s_0b^3\beta^2-9024r_{0m}s^m_0nb^2\beta^2 \\
\!\!\!\!\!\!&-&\!\!\!\!\!\!\ 1536r_{0m}s^m_0nb^4\beta^2-4512s_{0|0}nb^2\beta^2-768s_{0|0}nb^4\beta^2 -6048n\beta^2b^2\theta-4800r_{00}r_m^mb^2\beta^2 \\
\!\!\!\!\!\!&-&\!\!\!\!\!\!\ 768r_{00}r_m^mb^4\beta^2-4800r_{00|m}b^mb^2\beta^2-768r_{00|m}b^mb^4\beta^2+4800r_{0m|0}b^mb^2\beta^2
\\
\!\!\!\!\!\!&+&\!\!\!\!\!\!\ 768r_{0m|0}b^mb^4\beta^2+10316(\textbf{Ric}-{\bf \overline{Ric}})\beta^2 +13024r_{00}r_0nb\beta +14080r_{00}r_0nb^3\beta +2560r_{00}r_0n\beta b^5 \\
\!\!\!\!\!\!&+&\!\!\!\!\!\!\  1024r_0s_0n\beta^2b^2+6016r_0s_0nb\beta^2 +2048r_0s_0nb^3\beta^2 +1024(\textbf{Ric}-{\bf \overline{Ric}})\beta^2b^6+19776(\textbf{Ric}-{\bf \overline{Ric}})\beta^2b^2 \\
\!\!\!\!\!\!&+&\!\!\!\!\!\!\ 9600(\textbf{Ric}-{\bf \overline{Ric}})\beta^2b^4+612\beta^2\sigma +3488r_{00}^2b^4 -1728n\beta^2b^4\theta +12728r_{00}^2b^2 -2461r_{00}^2n \\
\!\!\!\!\!\!&+&\!\!\!\!\!\!\ 5710r_{00|0}\beta-1648s_0^2\beta^2 +1600r_0^2\beta^2+4944r_{0m}r^m_0\beta^2 +6880r_{0m}s^m_0\beta^2+10856s_{0|0}\beta^2
\\
\!\!\!\!\!\!&+&\!\!\!\!\!\!\ 1152s_ms^m_0\beta^3+864s_mr^m_0\beta^3+576r_ms^m_0\beta^3+3744s_{0|m}^m\beta^3 +4416r_{00|0}\beta b^4+512r_{00|0}b^6\beta \\
\!\!\!\!\!\!&-&\!\!\!\!\!\!\ 4240r_{00}s_0\beta-12464r_{00}r_0\beta-4670r_{00|0}n\beta+9504r_{00|0}\beta b^2+544s_0^2n\beta^2 +640s_0^2\beta^2b^2-16384s_0^2b\beta^2 \\
\!\!\!\!\!\!&-&\!\!\!\!\!\!\ 5120s_0^2b^3\beta^2 -1600rr_{00}\beta^2+512r_0^2b^2\beta^2 -6400r_0^2b\beta^2 -2048r_0^2b^3\beta^2 -416r_0s_0\beta^2 \\
\!\!\!\!\!\!&-&\!\!\!\!\!\!\  7432r_{00}^2nb^2-3760r_{00}^2nb^4 +4800r_{0m}r^m_0b^2\beta^2 +768r_{0m}r^m_0b^4\beta^2 +6208r_{0m}s^m_0b^2\beta^2
\\
\!\!\!\!\!\!&+&\!\!\!\!\!\!\ 1024r_{0m}s^m_0b^4\beta^2-8736r_{0m}s^m_0n\beta^2+10304s_{0|0}b^2\beta^2+1664s_{0|0}b^4\beta^2 -4368s_{0|0}n\beta^2 \\
\!\!\!\!\!\!&-&\!\!\!\!\!\!\ 576s_ms^m_0n\beta^3+3593r_{00}^2-4944r_{00}r_m^m\beta^2-4944r_{00|m}b^m\beta^2 +4944r_{0m|0}b^m\beta^2 -216s_ms_0^m\beta^2
\\
\!\!\!\!\!\!&-&\!\!\!\!\!\!\ 576s_0r_m^m\beta^3-576s_{0|m}b^m\beta^3+288s_{m|0}b^m\beta^3+2304s_{0|m}^mb^2\beta^3 +1728\beta^2b^4\theta-612n\beta^2\sigma \\
\!\!\!\!\!\!&-&\!\!\!\!\!\!\ 4050n\beta^2\theta+6048\beta^2b^2\theta+576\beta^2b^2\sigma-576n\beta^2b^2\sigma)
\\
A_6\!\!\!\!&:=&\!\!\!\!\ \beta^3(-2304n\beta^2b^6\theta+11061\beta^2\theta -1728n\beta^2b^4\sigma -24672r_{00}s_0nb^2\beta -11712r_{00}s_0n\beta b^4 \\
\!\!\!\!\!\!&+&\!\!\!\!\!\!\ 24336r_{00}s_0nb\beta+38016r_{00}s_0nb^3\beta+11520r_{00}s_0n\beta b^5+19008r_{00}r_0nb^2\beta +5760r_{00}r_0n\beta b^4 \\
\!\!\!\!\!\!&+&\!\!\!\!\!\!\ 6912r_{00}s_0b^4\beta+21600r_{00}s_0b^2\beta -4980r_{00}s_0n\beta -768r_{00}s_0b^5\beta - 20352r_{00}s_0b^3\beta \\
\!\!\!\!\!\!&-&\!\!\!\!\!\!\  33744r_{00}s_0b\beta-5760r_{00}r_0b^4\beta -28224r_{00}r_0b^2\beta+12168r_{00}r_0n\beta -768r_{00}r_0b^5\beta-20352r_{00}r_0b^3\beta \\
\!\!\!\!\!\!&-&\!\!\!\!\!\!\ 33744r_{00}r_0b\beta-18252r_{00|0}nb^2\beta-14256r_{00|0}nb^4\beta-2880r_{00|0}nb^6\beta+1216s_0^2n\beta^2b^2 \\
\!\!\!\!\!\!&-&\!\!\!\!\!\!\ 512s_0^2n\beta^2b^4 +28928s_0^2nb\beta^2+22016s_0^2nb^3\beta^2 +2048s_0^2nb^5\beta^2 -5888rr_{00}b^2\beta^2\\
\!\!\!\!\!\!&-&\!\!\!\!\!\!\ 512rr_{00}b^4\beta^2 +33152rr_{00}b\beta^2+23552rr_{00}b^3\beta^2 +2048rr_{00}b^5\beta^2 +14464r_0s_0n\beta^2 \\
\!\!\!\!\!\!&-&\!\!\!\!\!\!\ 1344r_0s_0\beta^2b^2 +768r_0s_0\beta^2b^4-94016r_0s_0b\beta^2 -60416r_0s_0b^3\beta^2 -5120r_0s_0b^5\beta^2
\\
\!\!\!\!\!\!&-&\!\!\!\!\!\!\ 43392r_{0m}s^m_0nb^2\beta^2-16512r_{0m}s^m_0nb^4\beta^2 -1024r_{0m}s^m_0n\beta^2b^6 -21696s_{0|0}nb^2\beta^2-8256s_{0|0}nb^4\beta^2 \\
\!\!\!\!\!\!&-&\!\!\!\!\!\!\ 512s_{0|0}n\beta^2b^6+4608rs_0b\beta^3 -26352n\beta^2b^2\theta-3456s_ms^m_0n\beta^3b^2 -24864r_{00}r_m^mb^2\beta^2-8832r_{00}r_m^mb^4\beta^2 \\
\!\!\!\!\!\!&-&\!\!\!\!\!\!\ 512r_{00}r_m^mb^6\beta^2-24864r_{00|m}b^mb^2\beta^2-8832r_{00|m}b^mb^4\beta^2-512r_{00|m}b^mb^6\beta^2+24864r_{0m|0}b^mb^2\beta^2 \\
\!\!\!\!\!\!&+&\!\!\!\!\!\!\ 8832r_{0m|0}b^mb^4\beta^2+512r_{0m|0}b^mb^6\beta^2 -1728s_ms_0^mb^2\beta^2 -3456s_0r_m^m\beta^3b^2-3456s_{0|m}b^m\beta^3b^2
\end{eqnarray*}
\begin{eqnarray*}
\!\!\!\!\!\!&+&\!\!\!\!\!\!\ 1728s_{m|0}b^m\beta^3b^2 +22976(\textbf{Ric}-{\bf \overline{Ric}})\beta^2 +24336r_{00}r_0nb\beta +38016r_{00}r_0nb^3\beta+11520r_{00}r_0n\beta b^5 \\
\!\!\!\!\!\!&+&\!\!\!\!\!\!\  11008r_0s_0n\beta^2b^2 +1024r_0s_0n\beta^2b^4+28928r_0s_0nb\beta^2 +22016r_0s_0nb^3\beta^2 +2048r_0s_0nb^5\beta^2 \\
\!\!\!\!\!\!&+&\!\!\!\!\!\!\ 11776(\textbf{Ric}-{\bf \overline{Ric}})\beta^2b^6+512(\textbf{Ric}-{\bf \overline{Ric}})\beta^2b^8 +62752(\textbf{Ric}-{\bf \overline{Ric}})\beta^2b^2+49728(\textbf{Ric}-{\bf \overline{Ric}})\beta^2b^4 \\
\!\!\!\!\!\!&+&\!\!\!\!\!\!\ 2214\beta^2\sigma+6672r_{00}^2b^4-16416n\beta^2b^4\theta +19272r_{00}^2b^2-2256r_{00}^2n +7554r_{00|0}\beta-20616s_0^2\beta^2 \\
\!\!\!\!\!\!&+&\!\!\!\!\!\!\  8288r_0^2\beta^2+15688r_{0m}r^m_0\beta^2+16912r_{0m}s^m_0\beta^2+31988s_{0|0}\beta^2+12672s_ms^m_0\beta^3 +10368s_mr^m_0\beta^3\\
\!\!\!\!\!\!&+&\!\!\!\!\!\!\ +6912r_ms^m_0\beta^3+20520s_{0|m}^m\beta^3+8640r_{00|0}\beta b^4+1536r_{00|0}b^6\beta -18036r_{00}s_0\beta -26064r_{00}r_0\beta \\
\!\!\!\!\!\!&-&\!\!\!\!\!\!\ 6570r_{00|0}n\beta+14832r_{00|0}\beta b^2+5560s_0^2n\beta^2 -2112s_0^2\beta^2b^2 +768s_0^2\beta^2b^4 -69568s_0^2b\beta^2 -43264s_0^2b^3\beta^2 \\
\!\!\!\!\!\!&-&\!\!\!\!\!\!\  4096s_0^2b^5\beta^2-8288rr_{00}\beta^2 +5888r_0^2b^2\beta^2+512r_0^2b^4\beta^2 -33152r_0^2b\beta^2-23552r_0^2b^3\beta^2 \\
\!\!\!\!\!\!&-&\!\!\!\!\!\!\ 2048r_0^2b^5\beta^2-8640r_0s_0\beta^2 -9240r_{00}^2nb^2-6504r_{00}^2nb^4 +24864r_{0m}r^m_0b^2\beta^2 +8832r_{0m}r^m_0b^4\beta^2 \\
\!\!\!\!\!\!&+&\!\!\!\!\!\!\  512r_{0m}r^m_0b^6\beta^2 +21984r_{0m}s^m_0b^2\beta^2 +7680r_{0m}s^m_0b^4\beta^2 +512r_{0m}s^m_0b^6\beta^2 +1728\beta^2b^4\sigma\\
\!\!\!\!\!\!&-&\!\!\!\!\!\!\ 25688r_{0m}s^m_0n\beta^2+48288s_{0|0}b^2\beta^2+17088s_{0|0}b^4\beta^2+1024s_{0|0}b^6\beta^2 -12844s_{0|0}n\beta^2 \\
\!\!\!\!\!\!&-&\!\!\!\!\!\!\  1152rs_0\beta^3-6480s_ms^m_0n\beta^3 +6336s_ms^m_0\beta^3b^2 +5184s_mr^m_0\beta^3b^2 +3456r_ms^m_0\beta^3b^2+3504r_{00}^2 \\
\!\!\!\!\!\!&-&\!\!\!\!\!\!\ 15688r_{00}r_m^m\beta^2-15688r_{00|m}b^m\beta^2+15688r_{0m|0}b^m\beta^2 -2160s_ms_0^m\beta^2 -6912s_0r_m^m\beta^3 \\
\!\!\!\!\!\!&-&\!\!\!\!\!\!\ 6912s_{0|m}b^m\beta^3 +3456s_{m|0}b^m\beta^3+27648s_{0|m}^mb^2\beta^3+6912s_{0|m}^mb^4\beta^3 +648s^i_ms^m_i\beta^4 \\
\!\!\!\!\!\!&+&\!\!\!\!\!\!\  16416\beta^2b^4\theta-2214n\beta^2\sigma -11061n\beta^2\theta+26352\beta^2b^2\theta+4320\beta^2b^2\sigma+2304\beta^2b^6\theta-4320n\beta^2b^2\sigma)
\\
A_7\!\!\!\!&:=&\!\!\!\!\ -\beta^2(-6528n\beta^2 b^6\theta+6171\beta^2\theta-3744n\beta^2b^4\sigma -12864r_{00}s_0nb^2\beta -9408r_{00}s_0n\beta b^4+9408r_{00}s_0nb\beta \\
\!\!\!\!\!\!&+&\!\!\!\!\!\!\ 19776r_{00}s_0nb^3\beta+8832r_{00}s_0n\beta b^5+9888r_{00}r_0nb^2\beta+4416r_{00}r_0n\beta b^4+3840r_{00}s_0b^4\beta \\
\!\!\!\!\!\!&+&\!\!\!\!\!\!\ 10896r_{00}s_0b^2\beta-1488r_{00}s_0n\beta+5760r_{00}s_0b^5\beta-1152r_{00}s_0b^3\beta -11520r_{00}s_0b\beta-4416r_{00}r_0b^4\beta \\
\!\!\!\!\!\!&-&\!\!\!\!\!\!\  16080r_{00}r_0b^2\beta+4704r_{00}r_0n\beta+5760r_{00}r_0b^5\beta -1152r_{00}r_0b^3\beta -11520r_{00}r_0b\beta-7056r_{00|0}nb^2\beta \\
\!\!\!\!\!\!&-&\!\!\!\!\!\!\ 7416r_{00|0}nb^4\beta-2208r_{00|0}nb^6\beta+8128s_0^2n\beta^2b^2 +2944s_0^2n\beta^2b^4 +24608s_0^2nb\beta^2+31232s_0^2nb^3\beta^2 \\
\!\!\!\!\!\!&+&\!\!\!\!\!\!\  6656s_0^2nb^5\beta^2-9088rr_{00}b^2\beta^2 -1792rr_{00}b^4\beta^2 +30784rr_{00}b\beta^2 +36352rr_{00}b^3\beta^2 +7168rr_{00}b^5\beta^2 \\
\!\!\!\!\!\!&+&\!\!\!\!\!\!\ 12304r_0s_0n\beta^2-13664r_0s_0\beta^2b^2-2432r_0s_0\beta^2b^4-76512r_0s_0b\beta^2 -75264r_0s_0b^3\beta^2 -13824r_0s_0b^5\beta^2 \\
\!\!\!\!\!\!&-&\!\!\!\!\!\!\  36912r_{0m}s^m_0nb^2\beta^2 -23424r_{0m}s^m_0nb^4\beta^2 -3328r_{0m}s^m_0n\beta^2b^6 -18456s_{0|0}nb^2\beta^2-11712s_{0|0}nb^4\beta^2
\\
\!\!\!\!\!\!&-&\!\!\!\!\!\!\ 1664s_{0|0}n\beta^2b^6-1536rs_0\beta^3b^2 +16896rs_0b\beta^3+6144rs_0\beta^3b^3 -20712n\beta^2b^2\theta -11808s_ms^m_0n\beta^3b^2 \\
\!\!\!\!\!\!&-&\!\!\!\!\!\!\ 2304s_ms^m_0n\beta^3b^4-23088r_{00}r_m^mb^2\beta^2-13632r_{00}r_m^mb^4\beta^2 -1792r_{00}r_m^mb^6\beta^2-23088r_{00|m}b^mb^2\beta^2 \\
\!\!\!\!\!\!&-&\!\!\!\!\!\!\ 13632r_{00|m}b^mb^4\beta^2-1792r_{00|m}b^mb^6\beta^2+23088r_{0m|0}b^mb^2\beta^2 +13632r_{0m|0}b^mb^4\beta^2\\
\!\!\!\!\!\!&+&\!\!\!\!\!\!\ 1792r_{0m|0}b^mb^6\beta^2 -5184s_ms_0^mb^2\beta^2 -1728s_ms_0^mb^4\beta^2 -12672s_0r_m^m\beta^3b^2-2304s_0r_m^m\beta^3b^4
\\
\!\!\!\!\!\!&-&\!\!\!\!\!\!\ 12672s_{0|m}b^m\beta^3b^2-2304s_{0|m}b^m\beta^3b^4+6336s_{m|0}b^m\beta^3b^2+1152s_{m|0}b^m\beta^3b^4 +1728s^i_ms^m_i\beta^4b^2 \\
\!\!\!\!\!\!&+&\!\!\!\!\!\!\ 11116(\textbf{Ric}-{\bf \overline{Ric}})\beta^2 +9408r_{00}r_0nb\beta+19776r_{00}r_0nb^3\beta +8832r_{00}r_0n\beta b^5+15616r_0s_0n\beta^2b^2 \\
\!\!\!\!\!\!&+&\!\!\!\!\!\!\ 3328r_0s_0n\beta^2b^4 +24608r_0s_0nb\beta^2 +31232r_0s_0nb^3\beta^2 +6656r_0s_0nb^5\beta^2 +18176(\textbf{Ric}-{\bf \overline{Ric}})\beta^2b^6 \\
\!\!\!\!\!\!&+&\!\!\!\!\!\!\ 1792(\textbf{Ric}-{\bf \overline{Ric}})\beta^2b^8+40352(\textbf{Ric}-{\bf \overline{Ric}})\beta^2b^2 +46176(\textbf{Ric}-{\bf \overline{Ric}})\beta^2b^4 -384n\beta^2b^8\theta -768n\beta^2b^6\sigma \\
\!\!\!\!\!\!&+&\!\!\!\!\!\!\ 1473\beta^2\sigma+2760r_{00}^2b^4 -20880n\beta^2b^4\theta +6264r_{00}^2b^2-450r_{00}^2n+2172r_{00|0}\beta-31588s_0^2\beta^2 \\
\!\!\!\!\!\!&+&\!\!\!\!\!\!\ 7696r_0^2\beta^2+10088r_{0m}r^m_0\beta^2+7872r_{0m}s^m_0\beta^2+19068s_{0|0}\beta^2 +18816s_ms^m_0\beta^3 +17064s_mr^m_0\beta^3
\\
\!\!\!\!\!\!&+&\!\!\!\!\!\!\ 11376r_ms^m_0\beta^3+20688s_{0|m}^m\beta^3+2376r_{00|0}\beta b^4+384r_{00|0}b^6\beta-11712r_{00}s_0\beta -11472r_{00}r_0\beta \\
\!\!\!\!\!\!&-&\!\!\!\!\!\!\ 2004r_{00|0}n\beta+4464r_{00|0}\beta b^2+7396s_0^2n\beta^2-19264s_0^2\beta^2b^2 -4480s_0^2\beta^2b^4-50592s_0^2b\beta^2-40320s_0^2b^3\beta^2 \\
\!\!\!\!\!\!&-&\!\!\!\!\!\!\ 6144s_0^2b^5\beta^2 -7696rr_{00}\beta^2+9088r_0^2b^2\beta^2 +1792r_0^2b^4\beta^2 -30784r_0^2b\beta^2 -36352r_0^2b^3\beta^2 \\
\!\!\!\!\!\!&-&\!\!\!\!\!\!\ 7168r_0^2b^5\beta^2-13856r_0s_0\beta^2-2400r_{00}^2nb^2-2244r_{00}^2nb^4 +23088r_{0m}r^m_0b^2\beta^2+13632r_{0m}r^m_0b^4\beta^2 \\
\!\!\!\!\!\!&+&\!\!\!\!\!\!\ 1792r_{0m}r^m_0b^6\beta^2 +10512r_{0m}s^m_0b^2\beta^2+4608r_{0m}s^m_0b^4\beta^2 +768r_{0m}s^m_0b^6\beta^2-15176r_{0m}s^m_0n\beta^2 \\
\!\!\!\!\!\!&+&\!\!\!\!\!\!\ 39888s_{0|0}b^2\beta^2+22752s_{0|0}b^4\beta^2+3072s_{0|0}b^6\beta^2-7588s_{0|0}n\beta^2-4224rs_0\beta^3
-9864s_ms^m_0n\beta^3
\end{eqnarray*}
\begin{eqnarray*}
\!\!\!\!\!\!&+&\!\!\!\!\!\!\ 20544s_ms^m_0\beta^3b^2+3840s_ms^m_0\beta^3b^4 +19008s_mr^m_0\beta^3b^2 +3456s_mr^m_0\beta^3b^4+12672r_ms^m_0\beta^3b^2 \\
\!\!\!\!\!\!&+&\!\!\!\!\!\!\ 2304r_ms^m_0\beta^3b^4 +714r_{00}^2+864s_ms^m\beta^4-10088r_{00}r_m^m\beta^2-10088r_{00|m}b^m\beta^2+10088r_{0m|0}b^m\beta^2 \\
\!\!\!\!\!\!&-&\!\!\!\!\!\!\ 2808s_ms_0^m\beta^2-11376s_0r_m^m\beta^3-11376s_{0|m}b^m\beta^3+5688s_{m|0}b^m\beta^3+45504s_{0|m}^mb^2\beta^3
+25344s_{0|m}^mb^4\beta^3 \\
\!\!\!\!\!\!&+&\!\!\!\!\!\!\ 3072s_{0|m}^m\beta^3b^6+2484s^i_ms^m_i\beta^4 +3744\beta^2b^4\sigma +20880\beta^2b^4\theta +768\beta^2b^6\sigma -1473n\beta^2\sigma \\
\!\!\!\!\!\!&-&\!\!\!\!\!\!\  6171n\beta^2\theta+20712\beta^2b^2\theta+4464\beta^2b^2\sigma+6528\beta^2b^6\theta +384\beta^2b^8\theta-4464n\beta^2b^2\sigma)
\\
A_8\!\!\!\!&:=&\!\!\!\!\ 2\beta(-3552n\beta^2b^6\theta+1083\beta^2\theta-1608n\beta^2b^4\sigma -2064r_{00}s_0nb^2\beta -2136r_{00}s_0n\beta b^4+1104r_{00}s_0nb\beta \\
\!\!\!\!\!\!&+&\!\!\!\!\!\!\ 2976r_{00}s_0nb^3\beta+1824r_{00}s_0n\beta b^5+1488r_{00}r_0nb^2\beta+912r_{00}r_0n\beta b^4+672r_{00}s_0b^4\beta+2040r_{00}s_0b^2\beta \\
\!\!\!\!\!\!&-&\!\!\!\!\!\!\ 156r_{00}s_0n\beta+2976r_{00}s_0b^5\beta+1536r_{00}s_0b^3\beta -1200r_{00}s_0b\beta -912r_{00}r_0b^4\beta-2712r_{00}r_0b^2\beta \\
\!\!\!\!\!\!&+&\!\!\!\!\!\!\ 552r_{00}r_0n\beta+2976r_{00}r_0b^5\beta +1536r_{00}r_0b^3\beta-1200r_{00}r_0b\beta-828r_{00|0}nb^2\beta-1116r_{00|0}nb^4\beta\\
\!\!\!\!\!\!&-&\!\!\!\!\!\!\ 456r_{00|0}nb^6\beta  +5360s_0^2n\beta^2b^2 +5168s_0^2n\beta^2b^4 +6016s_0^2nb\beta^2 +11200s_0^2nb^3\beta^2 +4096s_0^2nb^5\beta^2 \\
\!\!\!\!\!\!&-&\!\!\!\!\!\!\ 3616rr_{00}b^2\beta^2-1216rr_{00}b^4\beta^2+8320rr_{00}b\beta^2+14464rr_{00}b^3\beta^2 +4864rr_{00}b^5\beta^2+3008r_0s_0n\beta^2 \\
\!\!\!\!\!\!&-&\!\!\!\!\!\!\ 9840r_0s_0\beta^2b^2-5376r_0s_0\beta^2b^4-17920r_0s_0b\beta^2 -22144r_0s_0b^3\beta^2 -6016r_0s_0b^5\beta^2 -9024r_{0m}s^m_0nb^2\beta^2 \\
\!\!\!\!\!\!&-&\!\!\!\!\!\!\ 8400r_{0m}s^m_0nb^4\beta^2 -2048r_{0m}s^m_0n\beta^2b^6-4512s_{0|0}nb^2\beta^2-4200s_{0|0}nb^4\beta^2-1024s_{0|0}n\beta^2b^6 \\
\!\!\!\!\!\!&-&\!\!\!\!\!\!\ 2560rs_0\beta^3b^2 +12352rs_0b\beta^3 +10240rs_0\beta^3b^3-4776n\beta^2b^2\theta +1024rs_0\beta^3b^5 -256rs_0\beta^3b^4
\\
\!\!\!\!\!\!&+&\!\!\!\!\!\!\ 64\beta^2b^8\sigma-7896s_ms^m_0n\beta^3b^2-3552s_ms^m_0n\beta^3b^4 -6240r_{00}r_m^mb^2\beta^2 -5424r_{00}r_m^mb^4\beta^2
\\
\!\!\!\!\!\!&-&\!\!\!\!\!\!\ 1216r_{00}r_m^mb^6\beta^2-6240r_{00|m}b^mb^2\beta^2-5424r_{00|m}b^mb^4\beta^2-1216r_{00|m}b^mb^6\beta^2+6240r_{0m|0}b^mb^2\beta^2 \\
\!\!\!\!\!\!&+&\!\!\!\!\!\!\ 5424r_{0m|0}b^mb^4\beta^2+1216r_{0m|0}b^mb^6\beta^2 -2880s_ms_0^mb^2\beta^2 -2304s_ms_0^mb^4\beta^2-9264s_0r_m^m\beta^3b^2 \\
\!\!\!\!\!\!&-&\!\!\!\!\!\!\ 3840s_0r_m^m\beta^3b^4-9264s_{0|m}b^m\beta^3b^2-3840s_{0|m}b^m\beta^3b^4+4632s_{m|0}b^m\beta^3b^2 +1920s_{m|0}b^m\beta^3b^4 \\
\!\!\!\!\!\!&+&\!\!\!\!\!\!\ 3024s^i_ms^m_i\beta^4b^2 +1756(\textbf{Ric}-{\bf \overline{Ric}})\beta^2+1104r_{00}r_0nb\beta +2976r_{00}r_0nb^3\beta+1824r_{00}r_0n\beta b^5 \\
\!\!\!\!\!\!&+&\!\!\!\!\!\!\ 5600r_0s_0n\beta^2b^2+2048r_0s_0n\beta^2b^4+6016r_0s_0nb\beta^2 +11200r_0s_0nb^3\beta^2+4096r_0s_0nb^5\beta^2 \\
\!\!\!\!\!\!&+&\!\!\!\!\!\!\ 7232(\textbf{Ric}-{\bf \overline{Ric}})\beta^2b^6+1216(\textbf{Ric}-{\bf \overline{Ric}})\beta^2b^8 +8096(\textbf{Ric}-{\bf \overline{Ric}})\beta^2b^2 +12480(\textbf{Ric}-{\bf \overline{Ric}})\beta^2b^4 \\
\!\!\!\!\!\!&-&\!\!\!\!\!\!\ 480n\beta^2b^8\theta-64n\beta^2b^8\sigma-704n\beta^2b^6\sigma +292\beta^2\sigma+384s_mr^m_0\beta^3b^6+336r_{00}^2b^4-6876n\beta^2b^4\theta\\
\!\!\!\!\!\!&+&\!\!\!\!\!\!\  588r_{00}^2b^2-27r_{00}^2n+198r_{00|0}\beta-12176s_0^2\beta^2+2080r_0^2\beta^2 +2024r_{0m}r^m_0\beta^2+1040r_{0m}s^m_0\beta^2
\\
\!\!\!\!\!\!&+&\!\!\!\!\!\!\ 3556s_{0|0}\beta^2+7368s_ms^m_0\beta^3+7500s_mr^m_0\beta^3+5000r_ms^m_0\beta^3+6316s_{0|m}^m\beta^3-36r_{00|0}\beta b^4\\
\!\!\!\!\!\!&-&\!\!\!\!\!\!\ 144r_{00|0}b^6\beta-2100r_{00}s_0\beta-1560r_{00}r_0\beta-192r_{00|0}n\beta+360r_{00|0}\beta b^2+2448s_0^2n\beta^2-14208s_0^2\beta^2b^2 \\
\!\!\!\!\!\!&-&\!\!\!\!\!\!\ 8352s_0^2\beta^2b^4-10240s_0^2b\beta^2 -5120s_0^2b^3\beta^2+2176s_0^2b^5\beta^2 -2080rr_{00}\beta^2+3616r_0^2b^2\beta^2+1216r_0^2b^4\beta^2
\\
\!\!\!\!\!\!&-&\!\!\!\!\!\!\ 8320r_0^2b\beta^2-14464r_0^2b^3\beta^2-4864r_0^2b^5\beta^2 -5168r_0s_0\beta^2-180r_{00}^2nb^2-216r_{00}^2nb^4+6240r_{0m}r^m_0b^2\beta^2 \\
\!\!\!\!\!\!&+&\!\!\!\!\!\!\ 5424r_{0m}r^m_0b^4\beta^2+1216r_{0m}r^m_0b^6\beta^2 +192r_{0m}s^m_0b^2\beta^2 -1344r_{0m}s^m_0b^4\beta^2 -320r_{0m}s^m_0b^6\beta^2\\
\!\!\!\!\!\!&-&\!\!\!\!\!\!\ 2776r_{0m}s^m_0n\beta^2+9456s_{0|0}b^2\beta^2+7464s_{0|0}b^4\beta^2+1664s_{0|0}b^6\beta^2
-1388s_{0|0}n\beta^2-3088rs_0\beta^3 \\
\!\!\!\!\!\!&-&\!\!\!\!\!\!\  3956s_ms^m_0n\beta^3+12792s_ms^m_0\beta^3b^2 +5376s_ms^m_0\beta^3b^4 +13896s_mr^m_0\beta^3b^2+5760s_mr^m_0\beta^3b^4 \\
\!\!\!\!\!\!&+&\!\!\!\!\!\!\ 9264r_ms^m_0\beta^3b^2+3840r_ms^m_0\beta^3b^4 +39r_{00}^2+1512s_ms^m\beta^4 -2024r_{00}r_m^m\beta^2-2024r_{00|m}b^m\beta^2 \\
\!\!\!\!\!\!&+&\!\!\!\!\!\!\ 2024r_{0m|0}b^m\beta^2-912s_ms_0^m\beta^2-5000s_0r_m^m\beta^3-5000s_{0|m}b^m\beta^3 +2500s_{m|0}b^m\beta^3+20000s_{0|m}^mb^2\beta^3 \\
\!\!\!\!\!\!&+&\!\!\!\!\!\!\ 18528s_{0|m}^mb^4\beta^3+5120s_{0|m}^m\beta^3b^6+1944s^i_ms^m_i\beta^4 -384s_ms_0^m\beta^2b^6+1608\beta^2b^4\sigma+6876\beta^2b^4\theta \\
\!\!\!\!\!\!&+&\!\!\!\!\!\!\ 256s_{0|m}^m\beta^3b^8-1408s_0^2\beta^2b^6 +704\beta^2b^6\sigma-292n\beta^2\sigma -1083n\beta^2\theta+4776\beta^2b^2\theta \\
\!\!\!\!\!\!&+&\!\!\!\!\!\!\ 1220\beta^2b^2\sigma+3552\beta^2b^6\theta +512s_0^2\beta^2b^7+480\beta^2b^8\theta -256s_0r_m^m\beta^3b^6-256s_{0|m}b^m\beta^3b^6 \\
\!\!\!\!\!\!&+&\!\!\!\!\!\!\ 128s_{m|0}b^m\beta^3b^6-1220n\beta^2b^2\sigma +384s_ms^m_0b^6\beta^3 +864s_ms^mb^2\beta^4 +864s^i_ms^m_i\beta^4b^4 \\
\!\!\!\!\!\!&+&\!\!\!\!\!\!\  1280s_0^2nb^6\beta^2-256s_ms^m_0n\beta^3b^6-640r_0s_0b^6\beta^2+256r_ms^m_0\beta^3b^6)
\\
A_9\!\!\!\!&:=&\!\!\!\!\ 3744n\beta^2b^6\theta+128s_ms_0^mb^8\beta^2 -468\beta^2\theta +1368n\beta^2b^4\sigma+768r_{00}s_0nb^2\beta+1056r_{00}s_0n\beta b^4 \\
\!\!\!\!\!\!&-&\!\!\!\!\!\!\ 288r_{00}s_0nb\beta-960r_{00}s_0nb^3\beta-768r_{00}s_0n\beta b^5-480r_{00}r_0nb^2\beta -384r_{00}r_0n\beta b^4-384r_{00}s_0b^4\beta
\end{eqnarray*}
\begin{eqnarray*}
\!\!\!\!\!\!&-&\!\!\!\!\!\!\ 1008r_{00}s_0b^2\beta+48r_{00}s_0n\beta -2304r_{00}s_0b^5\beta-1152r_{00}s_0b^3\beta +288r_{00}s_0b\beta+384r_{00}r_0b^4\beta\\
\!\!\!\!\!\!&+&\!\!\!\!\!\!\ 1008r_{00}r_0b^2\beta-144r_{00}r_0n\beta -2304r_{00}r_0b^5\beta -1152r_{00}r_0b^3\beta +288r_{00}r_0b\beta+216r_{00|0}nb^2\beta \\
\!\!\!\!\!\!&+&\!\!\!\!\!\!\ 360r_{00|0}nb^4\beta+192r_{00|0}nb^6\beta-5888s_0^2n\beta^2b^2 -10288s_0^2n\beta^2b^4-3392s_0^2nb\beta^2-8576s_0^2nb^3\beta^2 \\
\!\!\!\!\!\!&-&\!\!\!\!\!\!\ 4736s_0^2nb^5\beta^2 +3136rr_{00}b^2\beta^2 +1600rr_{00}b^4\beta^2-5248rr_{00}b\beta^2 -12544rr_{00}b^3\beta^2 -6400rr_{00}b^5\beta^2 \\
\!\!\!\!\!\!&-&\!\!\!\!\!\!\ 1696r_0s_0n\beta^2+11744r_0s_0\beta^2b^2 +11840r_0s_0\beta^2b^4+9792r_0s_0b\beta^2 +12288r_0s_0b^3\beta^2+2688r_0s_0b^5\beta^2\\
\!\!\!\!\!\!&+&\!\!\!\!\!\!\ 5088r_{0m}s^m_0nb^2\beta^2+6432r_{0m}s^m_0nb^4\beta^2 +2368r_{0m}s^m_0n\beta^2b^6 +2544s_{0|0}nb^2\beta^2+3216s_{0|0}nb^4\beta^2 \\
\!\!\!\!\!\!&+&\!\!\!\!\!\!\ 1184s_{0|0}n\beta^2b^6+6528rs_0\beta^3b^2 -18432rs_0b\beta^3 -26112rs_0\beta^3b^3 +2592n\beta^2b^2\theta -6144rs_0\beta^3b^5 \\
\!\!\!\!\!\!&+&\!\!\!\!\!\!\ 1536rs_0\beta^3b^4-192\beta^2b^8\sigma +10560s_ms^m_0n\beta^3b^2 +8160s_ms^m_0n\beta^3b^4 +3936r_{00}r_m^mb^2\beta^2   +4704r_{00}r_m^mb^4\beta^2 \\
\!\!\!\!\!\!&+&\!\!\!\!\!\!\ 1600r_{00}r_m^mb^6\beta^2+3936r_{00|m}b^mb^2\beta^2+4704r_{00|m}b^mb^4\beta^2 +1600r_{00|m}b^mb^6\beta^2-3936r_{0m|0}b^mb^2\beta^2
\\
\!\!\!\!\!\!&-&\!\!\!\!\!\!\ 4704r_{0m|0}b^mb^4\beta^2-1600r_{0m|0}b^mb^6\beta^2+2944s_ms_0^mb^2\beta^2 +4224s_ms_0^mb^4\beta^2+13824s_0r_m^m\beta^3b^2 \\
\!\!\!\!\!\!&+&\!\!\!\!\!\!\ 9792s_0r_m^m\beta^3b^4+13824s_{0|m}b^m\beta^3b^2+9792s_{0|m}b^m\beta^3b^4-6912s_{m|0}b^m\beta^3b^2
-4896s_{m|0}b^m\beta^3b^4 \\
\!\!\!\!\!\!&-&\!\!\!\!\!\!\ 8352s^i_ms^m_i\beta^4b^2 -700(\textbf{Ric}-{\bf \overline{Ric}})\beta^2 -288r_{00}r_0nb\beta -960r_{00}r_0nb^3\beta -768r_{00}r_0n\beta b^5 \\
\!\!\!\!\!\!&-&\!\!\!\!\!\!\ 4288r_0s_0n\beta^2b^2-2368r_0s_0n\beta^2b^4 -3392r_0s_0nb\beta^2-8576r_0s_0nb^3\beta^2 -4736r_0s_0nb^5\beta^2\\
\!\!\!\!\!\!&-&\!\!\!\!\!\!\ 6272(\textbf{Ric}-{\bf \overline{Ric}})\beta^2b^6-1600(\textbf{Ric}-{\bf \overline{Ric}})\beta^2b^8 -3968(\textbf{Ric}-{\bf \overline{Ric}})\beta^2b^2-^2b^4 +864n\beta^2b^8\theta
\\
\!\!\!\!\!\!&+&\!\!\!\!\!\!\ 192n\beta^2b^8\sigma+960n\beta^2b^6\sigma-138\beta^2\sigma -2304s_mr^m_0\beta^3b^6-72r_{00}^2b^4+4968n\beta^2b^4\theta-96r_{00}^2b^2+3r_{00}^2n
\\
\!\!\!\!\!\!&-&\!\!\!\!\!\!\ 42r_{00|0}\beta+10152s_0^2\beta^2 -1312r_0^2\beta^2-992r_{0m}r^m_0\beta^2 -288r_{0m}s^m_0\beta^2-1632s_{0|0}\beta^2-6656s_ms^m_0\beta^3\\
\!\!\!\!\!\!&-&\!\!\!\!\!\!\ 7632s_mr^m_0\beta^3-5088r_ms^m_0\beta^3-4800s_{0|m}^m\beta^3+168r_{00|0}\beta b^4+192r_{00|0}b^6\beta+816r_{00}s_0\beta+480r_{00}r_0\beta\\
\!\!\!\!\!\!&+&\!\!\!\!\!\!\ 42r_{00|0}n\beta-48r_{00|0}\beta7872(\textbf{Ric}-{\bf \overline{Ric}})\beta b^2-1704s_0^2n\beta^2 +17600s_0^2\beta^2b^2+18016s_0^2\beta^2b^4+4672s_0^2b\beta^2\\
\!\!\!\!\!\!&-&\!\!\!\!\!\!\ 4608s_0^2b^3\beta^2-11904s_0^2b^5\beta^2 +1312rr_{00}\beta^2-3136r_0^2b^2\beta^2 -1600r_0^2b^4\beta^2 +5248r_0^2b\beta^2 +12544r_0^2b^3\beta^2 \\
\!\!\!\!\!\!&+&\!\!\!\!\!\!\  6400r_0^2b^5\beta^2+4032r_0s_0\beta^2+24r_{00}^2nb^2 +36r_{00}^2nb^4-3936r_{0m}r^m_0b^2\beta^2-4704r_{0m}r^m_0b^4\beta^2\\
\!\!\!\!\!\!&-&\!\!\!\!\!\!\ 1600r_{0m}r^m_0b^6\beta^2+1440r_{0m}s^m_0b^2\beta^2 +4032r_{0m}s^m_0b^4\beta^2 +1728r_{0m}s^m_0b^6\beta^2+1232r_{0m}s^m_0n\beta^2\\
\!\!\!\!\!\!&-&\!\!\!\!\!\!\ 5184s_{0|0}b^2\beta^2-5040s_{0|0}b^4\beta^2-1536s_{0|0}b^6\beta^2+616s_{0|0}n\beta^2 +4608rs_0\beta^3+3632s_ms^m_0n\beta^3
\\
\!\!\!\!\!\!&-&\!\!\!\!\!\!\ 15552s_ms^m_0\beta^3b^2-10560s_ms^m_0\beta^3b^4 -20736s_mr^m_0\beta^3b^2 -14688s_mr^m_0\beta^3b^4-13824r_ms^m_0\beta^3b^2\\
\!\!\!\!\!\!&-&\!\!\!\!\!\!\ 3r_{00}^2-4176s_ms^m\beta^4+992r_{00}r_m^m\beta^2+992r_{00|m}b^m\beta^2-992r_{0m|0}b^m\beta^2 +632s_ms_0^m\beta^2+5088s_0r_m^m\beta^3\\
\!\!\!\!\!\!&+&\!\!\!\!\!\!\ 5088s_{0|m}b^m\beta^3-2544s_{m|0}b^m\beta^3-20352s_{0|m}^mb^2\beta^3-27648s_{0|m}^mb^4\beta^3 -13056s_{0|m}^m\beta^3b^6 -3228s^i_ms^m_i\beta^4 \\
\!\!\!\!\!\!&+&\!\!\!\!\!\!\ 1792s_ms_0^m\beta^2b^6-768s^i_ms^m_i\beta^4b^6 -1368\beta^2b^4\sigma -4968\beta^2b^4\theta-1536s_{0|m}^m\beta^3b^8 +512s_0^2b^8\beta^2
\\
\!\!\!\!\!\!&+&\!\!\!\!\!\!\ 6784s_0^2\beta^2b^6-960\beta^2b^6\sigma +138n\beta^2\sigma+468n\beta^2\theta -2592\beta^2b^2\theta-744\beta^2b^2\sigma-3744\beta^2b^6\theta \\
\!\!\!\!\!\!&-&\!\!\!\!\!\!\ 2560s_0^2\beta^2b^7-864\beta^2b^8\theta+1536s_0r_m^m\beta^3b^6+1536s_{0|m}b^m\beta^3b^6 -768s_{m|0}b^m\beta^3b^6 +744n\beta^2b^2\sigma \\
\!\!\!\!\!\!&-&\!\!\!\!\!\!\ 1792s_ms^m_0b^6\beta^3-5472s_ms^mb^2\beta^4
-5472s^i_ms^m_i\beta^4b^4 -1152s_ms^m\beta^4b^4 -6016s_0^2nb^6\beta^2\\
\!\!\!\!\!\!&+&\!\!\!\!\!\!\ 1408s_ms^m_0n\beta^3b^6+3200r_0s_0b^6\beta^2-1536r_ms^m_0\beta^3b^6
-512s_0^2nb^8\beta^2-9792r_ms^m_0\beta^3b^4
\\
A_{10}\!\!\!\!&:=&\!\!\!\!\ -256s_ms_0^mb^8\beta+128s^i_ms^m_ib^8\beta^3 +256s_ms^m\beta^3b^6 +2176s^i_ms^m_i\beta^3b^6 +3264s_ms^m\beta^3b^4 -64r_{00}s_0nb^2 \\
\!\!\!\!\!\!&-&\!\!\!\!\!\!\  112r_{00}s_0nb^4+64r_{00}s_0nb^3 +64r_{00}s_0nb^5 +16r_{00}s_0nb+32r_{00}r_0nb^2 +32r_{00}r_0nb^4+64r_{00}r_0nb^3\\
\!\!\!\!\!\!&+&\!\!\!\!\!\!\ 64r_{00}r_0nb^5+16r_{00}r_0nb+1440s_0^2n\beta b^2+3904s_0^2n\beta b^4+4352s_0^2n\beta b^6+512s_0^2nb\beta+1664s_0^2nb^3\beta \\
\!\!\!\!\!\!&+&\!\!\!\!\!\!\ 1280s_0^2nb^5\beta-704rr_{00}b^2\beta -512rr_{00}b^4\beta+896rr_{00}b\beta+2816rr_{00}b^3\beta +2048rr_{00}b^5\beta+256r_0s_0n\beta\\
\!\!\!\!\!\!&-&\!\!\!\!\!\!\ 3168r_0s_0\beta b^2-4992r_0s_0\beta b^4-1280r_0s_0\beta b^3-1472r_0s_0b\beta-2560r_0s_0\beta b^6+1024r_0s_0\beta b^5\\
\!\!\!\!\!\!&-&\!\!\!\!\!\!\ 768r_{0m}s^m_0nb^2\beta-1248r_{0m}s^m_0nb^4\beta-640r_{0m}s^m_0nb^6\beta-384s_{0|0}nb^2\beta-624s_{0|0}nb^4\beta \\
\!\!\!\!\!\!&-&\!\!\!\!\!\!\ 320s_{0|0}nb^6\beta +7424rs_0b\beta^2+15872rs_0b^3\beta^2 -1664rs_0b^4\beta^2+6656rs_0b^5\beta^2 -3744s_ms^m_0n\beta^2b^2 \\
\!\!\!\!\!\!&-&\!\!\!\!\!\!\ 4320s_ms^m_0n\beta^2b^4-1344s_ms^m_0nb^6\beta^2 -3968rs_0b^2\beta^2 -672r_{00}r_m^mb^2\beta -1056r_{00}r_m^mb^4\beta\\
\!\!\!\!\!\!&-&\!\!\!\!\!\!\ 512r_{00}r_m^mb^6\beta-672r_{00|m}b^mb^2\beta-1056r_{00|m}b^mb^4\beta-512r_{00|m}b^mb^6\beta+672r_{0m|0}b^mb^2\beta\\
\!\!\!\!\!\!&+&\!\!\!\!\!\!\ 1056r_{0m|0}b^mb^4\beta+512r_{0m|0}b^mb^6\beta-704s_ms_0^m\beta b^2-1536s_ms_0^m\beta b^4-1280s_ms_0^m\beta b^6-5568s_0r_m^mb^2\beta^2
\end{eqnarray*}
\begin{eqnarray*}
\!\!\!\!\!\!&-&\!\!\!\!\!\!\ 5952s_0r_m^m\beta^2b^4-1664s_0r_m^m\beta^2b^6-5568s_{0|m}b^mb^2\beta^2-5952s_{0|m}b^m\beta^2b^4-1664s_{0|m}b^m\beta^2b^6\\
\!\!\!\!\!\!&+&\!\!\!\!\!\!\ 2784s_{m|0}b^mb^2\beta^2+2976s_{m|0}b^m\beta^2b^4+832s_{m|0}b^m\beta^2b^6 +6528s_ms^m\beta^3b^2+5824s^i_ms^m_i\beta^3b^2
\\
\!\!\!\!\!\!&+&\!\!\!\!\!\!\ 6528s^i_ms^m_i\beta^3b^4+1112s_{0|m}^m\beta^2+80(\textbf{Ric}-{\bf \overline{Ric}})\beta+832r_0s_0n\beta b^2+640r_0s_0n\beta b^4+512r_0s_0nb\beta\\
\!\!\!\!\!\!&+&\!\!\!\!\!\!\ 1664r_0s_0nb^3\beta+1280r_0s_0nb^5\beta+1408(\textbf{Ric}-{\bf \overline{Ric}})\beta b^6+512(\textbf{Ric}-{\bf \overline{Ric}})\beta b^8+544(\textbf{Ric}-{\bf \overline{Ric}})\beta b^2\\
\!\!\!\!\!\!&+&\!\!\!\!\!\!\ 1344(\textbf{Ric}-{\bf \overline{Ric}})\beta b^4-68r_{00}s_0-32r_{00}r_0-24r_{00|0}b^4-2r_{00|0}n-32r_{00|0}b^6 -2200s_0^2\beta+224r_0^2\beta\\
\!\!\!\!\!\!&+&\!\!\!\!\!\!\ 136r_{0m}r^m_0\beta +16r_{0m}s^m_0\beta+212s_{0|0}\beta+1760s_ms^m_0\beta^2 +2256s_mr^m_0\beta^2+1504r_ms^m_0\beta^2\\
\!\!\!\!\!\!&+&\!\!\!\!\!\!\ 64r_{00}s_0b^4+112r_{00}s_0b^2-4r_{00}s_0n +320r_{00}s_0b^5+128r_{00}s_0b^3 -16r_{00}s_0b-32r_{00}r_0b^4-80r_{00}r_0b^2
\\
\!\!\!\!\!\!&+&\!\!\!\!\!\!\ 8r_{00}r_0n+320r_{00}r_0b^5+128r_{00}r_0b^3-16r_{00}r_0b-12r_{00|0}nb^2-24r_{00|0}nb^4-16r_{00|0}nb^6+296s_0^2n\beta \\
\!\!\!\!\!\!&-&\!\!\!\!\!\!\ 4928s_0^2\beta b^2-7680s_0^2\beta b^4+2048s_0^2\beta b^7+2816s_0^2\beta b^3-576s_0^2b\beta -5120s_0^2\beta b^6+6656s_0^2\beta b^5 \\
\!\!\!\!\!\!&-&\!\!\!\!\!\!\ 224rr_{00}\beta+704r_0^2b^2\beta+512r_0^2b^4\beta -896r_0^2b\beta-2816r_0^2b^3\beta -2048r_0^2b^5\beta-800r_0s_0\beta+672r_{0m}r^m_0b^2\beta \\
\!\!\!\!\!\!&+&\!\!\!\!\!\!\ 1056r_{0m}r^m_0b^4\beta+512r_{0m}r^m_0b^6\beta-480r_{0m}s^m_0\beta b^2 -1536r_{0m}s^m_0\beta b^4-1024r_{0m}s^m_0\beta b^6 \\
\!\!\!\!\!\!&-&\!\!\!\!\!\!\ 152r_{0m}s^m_0n\beta+768s_{0|0}\beta b^2+816s_{0|0}\beta b^4+256s_{0|0}\beta b^6-76s_{0|0}n\beta-1856rs_0\beta^2-960s_ms^m_0n\beta^2
\\
\!\!\!\!\!\!&+&\!\!\!\!\!\!\ 4896s_ms^m_0\beta^2b^2+4224s_ms^m_0\beta^2b^4 +1216s_ms^m_0b^6\beta^2 +8352s_mr^m_0b^2\beta^2 +8928s_mr^m_0\beta^2b^4
\\
\!\!\!\!\!\!&+&\!\!\!\!\!\!\ 2496s_mr^m_0\beta^2b^6+5568r_ms^m_0b^2\beta^2 +5952r_ms^m_0\beta^2b^4 +1664r_ms^m_0\beta^2b^6 +57\beta\theta+18\beta\sigma +2r_{00|0}\\
\!\!\!\!\!\!&+&\!\!\!\!\!\!\ 2912s_ms^m\beta^3-136r_{00}r_m^m\beta-136r_{00|m}b^m\beta+136r_{0m|0}b^m\beta -112s_ms_0^m\beta-1504s_0r_m^m\beta^2 \\
\!\!\!\!\!\!&-&\!\!\!\!\!\!\ 1504s_{0|m}b^m\beta^2+752s_{m|0}b^m\beta^2+6016s_{0|m}^mb^2\beta^2+11136s_{0|m}^mb^4\beta^2+7936s_{0|m}^mb^6\beta^2+1664s_{0|m}^m\beta^2b^8 \\
\!\!\!\!\!\!&+&\!\!\!\!\!\!\  1544s^i_ms^m_i\beta^3+1024s_0^2nb^8\beta-384n\beta b^2\theta-120n\beta b^2\sigma-936n\beta b^4\theta-288n\beta b^4\sigma-960n\beta b^6\theta\\
\!\!\!\!\!\!&-&\!\!\!\!\!\!\ 288n\beta b^6\sigma-336n\beta b^8\theta-96n\beta b^8\sigma-57n\beta\theta-18n\beta\sigma+384\beta b^2\theta+120\beta b^2\sigma+936\beta b^4\theta\\
\!\!\!\!\!\!&+&\!\!\!\!\!\!\ 288\beta b^4\sigma+960\beta b^6\theta+288\beta b^6\sigma+336\beta b^8\theta+96\beta b^8\sigma-1024s_0^2b^8\beta
\\
A_{11}\!\!\!\!&:=&\!\!\!\!\ -4(\textbf{Ric}-{\bf \overline{Ric}})+24nb^4\sigma+32nb^6\sigma+8nb^2\sigma+48nb^8\theta +16nb^8\sigma+96nb^6\theta+72nb^4\theta-64r_0s_0nb^2\\
\!\!\!\!\!\!&-&\!\!\!\!\!\!\ 64r_0s_0nb^4-128r_0s_0nb^3-128r_0s_0nb^5-32r_0s_0nb+1152rs_0\beta b^2+768rs_0\beta b^4-320s^i_ms^m_ib^8\beta^2\\
\!\!\!\!\!\!&-&\!\!\!\!\!\!\ 640s_ms^m\beta^2b^6 -3072rs_0\beta b^5-4608rs_0\beta b^3-1536rs_0b\beta+672s_ms^m_0n\beta b^2+1056s_ms^m_0n\beta b^4\\
\!\!\!\!\!\!&+&\!\!\!\!\!\!\ 512s_ms^m_0n\beta b^6+1152s_0r_m^m\beta b^2+1728s_0r_m^m\beta b^4+768s_0r_m^m\beta b^6+1152s_{0|m}b^m\beta b^2+1728s_{0|m}b^m\beta b^4\\
\!\!\!\!\!\!&+&\!\!\!\!\!\!\ 768s_{0|m}b^m\beta b^6-576s_{m|0}b^m\beta b^2-864s_{m|0}b^m\beta b^4-384s_{m|0}b^m\beta b^6-3648s_ms^mb^2\beta^2-3264s_ms^m\beta^2b^4\\
\!\!\!\!\!\!&-&\!\!\!\!\!\!\ 2176s^i_ms^m_ib^2\beta^2-3648s^i_ms^m_i\beta^2b^4 -2176s^i_ms^m_i\beta^2b^6 -128(\textbf{Ric}-{\bf \overline{Ric}})b^6 -64(\textbf{Ric}-{\bf \overline{Ric}})b^8
\\
\!\!\!\!\!\!&-&\!\!\!\!\!\!\ 32(\textbf{Ric}-{\bf \overline{Ric}})b^2 -96(\textbf{Ric}-{\bf \overline{Ric}})b^4+24nb^2\theta+512s_0^2b^2 +1120s_0^2b^4-20s_0^2n+1152s_0^2b^6\\
\!\!\!\!\!\!&-&\!\!\!\!\!\!\ 384s_0^2b^3-1152s_0^2b^5-512s_0^2b^7+32s_0^2b+512s_0^2b^8 -64r_0^2b^2-64r_0^2b^4 +64r_0^2b+256r_0^2b^3+256r_0^2b^5\\
\!\!\!\!\!\!&+&\!\!\!\!\!\!\ 64r_0s_0-48r_{0m}r^m_0b^2-96r_{0m}r^m_0b^4 -64r_{0m}r^m_0b^6+48r_{0m}s^m_0b^2+192r_{0m}s^m_0b^4 +8s_ms_0^m\\
\!\!\!\!\!\!&+&\!\!\!\!\!\!\ 192r_{0m}s^m_0b^6+8r_{0m}s^m_0n-48s_{0|0}b^2-48s_{0|0}b^4+4s_{0|0}n-256s_ms^m_0\beta -360s_mr^m_0\beta-240r_ms^m_0\beta\\
\!\!\!\!\!\!&-&\!\!\!\!\!\!\ 144s_{0|m}^m\beta-3\theta-\sigma-8r_{0m}r^m_0-128s_0^2nb^2 -496s_0^2nb^4-896s_0^2nb^6 -128s_0^2nb^3-128s_0^2nb^5\\
\!\!\!\!\!\!&+&\!\!\!\!\!\!\ 64rr_{00}b^2+64rr_{00}b^4-64rr_{00}b-256rr_{00}b^3 -256rr_{00}b^5+320r_0s_0b^2 +704r_0s_0b^4-16r_0s_0n\\
\!\!\!\!\!\!&+&\!\!\!\!\!\!\ 96r_0s_0b+48r_{0m}s^m_0nb^2+96r_{0m}s^m_0nb^4+64r_{0m}s^m_0nb^6+24s_{0|0}nb^2 +48s_{0|0}nb^4+32s_{0|0}nb^6\\
\!\!\!\!\!\!&+&\!\!\!\!\!\!\ 136s_ms^m_0n\beta-768s_ms^m_0\beta b^2-576s_ms^m_0\beta b^4-128s_ms^m_0\beta b^6-1728s_mr^m_0\beta b^2-2592s_mr^m_0\beta b^4\\
\!\!\!\!\!\!&-&\!\!\!\!\!\!\ 1152s_mr^m_0\beta b^6-1152r_ms^m_0\beta b^2-1728r_ms^m_0\beta b^4-768r_ms^m_0\beta b^6-1088s_ms^m\beta^2+48r_{00}r_m^mb^2\\
\!\!\!\!\!\!&+&\!\!\!\!\!\!\ 96r_{00}r_m^mb^4+64r_{00}r_m^mb^6+48r_{00|m}b^mb^2+96r_{00|m}b^mb^4+64r_{00|m}b^mb^6-48r_{0m|0}b^mb^2-96r_{0m|0}b^mb^4\\
\!\!\!\!\!\!&-&\!\!\!\!\!\!\ 64r_{0m|0}b^mb^6+192s_ms_0^mb^4+64s_ms_0^mb^2 +256s_ms_0^mb^6+128s_ms_0^mb^8 +240s_0r_m^m\beta+240s_{0|m}b^m\beta\\
\!\!\!\!\!\!&-&\!\!\!\!\!\!\ 120s_{m|0}b^m\beta-960s_{0|m}^m\beta b^2-8b^2\sigma-72b^4\theta-24b^4\sigma-96b^6\theta -32b^6\sigma-48b^8\theta-16b^8\sigma+16rr_{00}\\
\!\!\!\!\!\!&+&\!\!\!\!\!\!\ 8r_{00}r_m^m+8r_{00|m}b^m-8r_{0m|0}b^m+196s_0^2-16r_0^2-12s_{0|0}-2304s_{0|m}^m\beta b^4-2304s_{0|m}^m\beta b^6-768s_{0|m}^m\beta b^8\\
\!\!\!\!\!\!&-&\!\!\!\!\!\!\ 428s^i_ms^m_i\beta^2+3n\theta-24b^2\theta+n\sigma -32s_0^2nb-512s_0^2nb^8 +640r_0s_0b^6-384r_0s_0b^5 +384rs_0\beta,
\end{eqnarray*}
\begin{eqnarray*}
A_{12}\!\!\!\!&:=&\!\!\!\!\ 8(1+2b^2)^2(8s^i_ms^m_i\beta b^4+4s_{0|m}^mb^4-2s_ms^m_0b^2-4s_{0|m}b^mb^2+6s_mr^m_0b^2+16s_ms^m\beta b^2+4s_{0|m}^mb^2\\
\!\!\!\!\!\!&+&\!\!\!\!\!\!\ 4r_ms^m_0b^2+20s^i_ms^m_i\beta b^2-4s_0r_m^mb^2-2s_ms^m_0nb^2+2s_{m|0}b^mb^2+16rs_0b+3s_mr^m_0+2r_ms^m_0\\
\!\!\!\!\!\!&+&\!\!\!\!\!\!\ s_{0|m}^m-4rs_0-2s_0r_m^m-2s_{0|m}b^m+8s^i_ms^m_i\beta+26s_ms^m\beta-s_ms^m_0n +s_{m|0}b^m+2s_ms^m_0),
\\
A_{13}\!\!\!\!&:=&\!\!\!\!\ -4(1+2b^2)^3(2s^i_ms^m_ib^2+s^i_ms^m_i+4s_ms^m).
\end{eqnarray*}

\section{Appendix 7: Coefficients in (\ref{s5})}
\begin{eqnarray*}
A'_0 \!\!\!\!&:=&\!\!\!\!\  24\beta^9(108(\textbf{Ric}-{\bf \overline{Ric}})
\\
A'_2 \!\!\!\!&:=&\!\!\!\!\ -288\beta^9c^2(-11+8n)+24\beta^7(-144c_0nb^2\beta +192cs_0nb\beta -384c^2b\beta^2 +192c^2nb\beta -96cs_0n\beta \\
\!\!\!\!\!\!&-&\!\!\!\!\!\!\ 342c_0n\beta +192c_0\beta b^2 +120cs_0\beta -144c^2\beta^2 +96c^2n\beta^2 -384cs_0b\beta +492c_0\beta) +2\beta^5(3726\beta^2\theta \\
\!\!\!\!\!\!&-&\!\!\!\!\!\!\ 2592n\beta^2b^2\theta +1728c_0b^2\beta^3 +13068(\textbf{Ric}-{\bf \overline{Ric}})\beta^2 +2304cs_0nb\beta^3 +15552(\textbf{Ric}-{\bf \overline{Ric}})\beta^2b^2  \\
\!\!\!\!\!\!&+&\!\!\!\!\!\!\ 3456(\textbf{Ric}-{\bf \overline{Ric}})\beta^2b^4 +324\beta^2\sigma +144s_0^2\beta^2 +576c^2\beta^4 +9144s_{0|0}\beta^2 +1296s_{0|m}^m\beta^3 -6912s_0^2b\beta^2 \\
\!\!\!\!\!\!&-&\!\!\!\!\!\!\ 2304c^2b\beta^4 +288cs_0\beta^3 +4032s_{0|0}b^2\beta^2 -3672s_{0|0}n\beta^2 +3888c_0\beta^3 \\
\!\!\!\!\!\!&-&\!\!\!\!\!\!\ 324n\beta^2\sigma -3726n\beta^2\theta +2592\beta^2b^2\theta)
\\
A'_4 \!\!\!\!&:=&\!\!\!\!\ 24\beta^7(1144c^2+928c^2b^2+64c^2b^4-865c^2n -64c^2nb^4-712c^2nb^2) +\beta^3(-2304n\beta^2b^6\theta +11061\beta^2\theta\\
\!\!\!\!\!\!&-&\!\!\!\!\!\!\ 1728n\beta^2b^4\sigma +1216s_0^2n\beta^2b^2 -512s_0^2n\beta^2b^4 +28928s_0^2nb\beta^2 +22016s_0^2nb^3\beta^2 + 2048s_0^2nb^5\beta^2 \\
\!\!\!\!\!\!&+&\!\!\!\!\!\!\ 14464cs_0n\beta^3 -1344cs_0\beta^3b^2 +768cs_0\beta^3b^4 -94016cs_0b\beta^3 -60416cs_0b^3\beta^2-5120c\beta s_0b^5\beta^2\\
\!\!\!\!\!\!&-&\!\!\!\!\!\!\ 21696s_{0|0}nb^2\beta^3 -8256s_{0|0}nb^4\beta^2 -512s_{0|0}n\beta^2b^6 +4608cb^2s_0b\beta^3 -26352n\beta^2b^2\theta -3456s_ms^m_0n\beta^3b^2 \\
\!\!\!\!\!\!&+&\!\!\!\!\!\!\ 24864c_0b^2\beta^3 +8832c_0b^4\beta^3 +512c_0b^6\beta^3 -1728s_ms_0^mb^2\beta^2 -3456s_0cn\beta^3b^2 -3456s_{0|m}b^m\beta^3b^2 \\
\!\!\!\!\!\!&+&\!\!\!\!\!\!\ 1728s_{m|0}b^m\beta^3b^2 +22976(\textbf{Ric}-{\bf \overline{Ric}})\beta^2 +11008cs_0n\beta^3b^2 +1024cs_0n\beta^3b^4 +28928cs_0nb\beta^3 \\
\!\!\!\!\!\!&+&\!\!\!\!\!\!\ 22016cs_0nb^3\beta^3 +2048cs_0nb^5\beta^3 +11776(\textbf{Ric}-{\bf \overline{Ric}})\beta^2b^6 +512(\textbf{Ric}-{\bf \overline{Ric}})\beta^2b^8 \\
\!\!\!\!\!\!&+&\!\!\!\!\!\!\ 62752(\textbf{Ric}-{\bf \overline{Ric}})\beta^2b^2 +49728(\textbf{Ric}-{\bf \overline{Ric}})\beta^2b^4 +2214\beta^2\sigma -16416n\beta^2b^4\theta -20616s_0^2\beta^2
\\
\!\!\!\!\!\!&+&\!\!\!\!\!\!\ 8288c^2\beta^4 +31988s_{0|0}\beta^2 +12672s_ms^m_0\beta^3 +10368cs_0\beta^3 +6912cs_0\beta^3 +20520s_{0|m}^m\beta^3 +5560s_0^2n\beta^2 \\
\!\!\!\!\!\!&-&\!\!\!\!\!\!\ 2112s_0^2\beta^2b^2 +768s_0^2\beta^2b^4 -69568s_0^2b\beta^2 -43264s_0^2b^3\beta^2 -4096s_0^2b^5\beta^2 +5888c^2b^2\beta^4 \\
\!\!\!\!\!\!&+&\!\!\!\!\!\!\ 512c^2b^4\beta^4 -33152c^2b\beta^4 -23552c^2b^3\beta^4 -2048c^2b^5\beta^4 -8640cs_0\beta^3 +48288s_{0|0}b^2\beta^2 \\
\!\!\!\!\!\!&+&\!\!\!\!\!\!\ 17088s_{0|0}b^4\beta^2 +1024s_{0|0}b^6\beta^2 -12844s_{0|0}n\beta^2 -1152cb^2s_0\beta^3 -6480s_ms^m_0n\beta^3 +6336s_ms^m_0\beta^3b^2 \\
\!\!\!\!\!\!&+&\!\!\!\!\!\!\ 5184cs_0\beta^3b^2 +3456cs_0\beta^3b^2 +15688c_0\beta^3 -2160s_ms_0^m\beta^2 -6912s_0cn\beta^3 -6912s_{0|m}b^m\beta^3 \\
\!\!\!\!\!\!&+&\!\!\!\!\!\!\ 3456s_{m|0}b^m\beta^3 +27648s_{0|m}^mb^2\beta^3 +6912s_{0|m}^mb^4\beta^3 +648s^i_ms^m_i\beta^4 +1728\beta^2b^4\sigma \\
\!\!\!\!\!\!&+&\!\!\!\!\!\!\ 16416\beta^2b^4\theta -2214n\beta^2\sigma -11061n\beta^2\theta +26352\beta^2b^2\theta +4320\beta^2b^2\sigma +2304\beta^2b^6\theta  -4320n\beta^2b^2\sigma) \\
\!\!\!\!\!\!&+&\!\!\!\!\!\!\ 2\beta^5( -8896cs_0nb^2\beta -1024cs_0n\beta b^4 +18976cs_0nb\beta +12544cs_0nb^3\beta +1024cs_0n\beta b^5 +6272c^2nb^2\beta^2 \\
\!\!\!\!\!\!&+&\!\!\!\!\!\!\ 512c^2n\beta^2 b^4+1024cs_0b^4\beta +10432cs_0b^2\beta -6568cs_0n\beta -1024cs_0b^5\beta -15616cs_0b^3\beta -32608cs_0b\beta \\
\!\!\!\!\!\!&-&\!\!\!\!\!\!\ 512c^2b^4\beta^2 -8192c^2b^2\beta^2 +9488c^2n\beta^2 -1024c^2b^5\beta^2 -15616c^2b^3\beta^2 -32608c^2b\beta^2 -14232c_0nb^2\beta \\
\!\!\!\!\!\!&-&\!\!\!\!\!\!\ 4704c_0nb^4\beta -256c_0nb^6\beta +2304s_0^2nb\beta^2 +2304c^2b^3\beta^2 +1152cs_0n\beta^3 -8064cs_0b\beta^3 -1728s_{0|0}nb^2\beta^2 \\
\!\!\!\!\!\!&-&\!\!\!\!\!\!\ 1728c^2nb^2\beta^2 -1728c_mb^mb^2\beta^2 +18976c^2nb\beta^2 +12544c^2nb^3\beta^2 +1024c^2n\beta^2 b^5 +12584c_0\beta \\
\!\!\!\!\!\!&+&\!\!\!\!\!\!\ 3888c^2\beta^2 +4800c_0\beta b^4 +256c_0b^6\beta +640cs_0\beta -16496c^2\beta^2 -9698c_0n\beta +15840c_0\beta b^2 -576c^2b^2\beta^2 \\
\!\!\!\!\!\!&+&\!\!\!\!\!\!\ 1728c^2b^2\beta^2 -3888c^2n\beta^2 -3888c_mb^m\beta^2 -3888c^2n\beta^2 -3888c_mb^m\beta^2)
\\
A'_6 \!\!\!\!&:=&\!\!\!\!\ 2\beta(-3552n\beta^2b^6\theta +1083\beta^2\theta -1608n\beta^2b^4\sigma +5360s_0^2n\beta^2b^2 +5168s_0^2n\beta^2b^4 +6016s_0^2nb\beta^2 \\
\!\!\!\!\!\!&+&\!\!\!\!\!\!\ 11200s_0^2nb^3\beta^2 +4096s_0^2nb^5\beta^2 +3008cs_0n\beta^3 -9840cs_0\beta^3b^2 -5376cs_0\beta^3b^4 -17920cs_0b\beta^3
\end{eqnarray*}
\begin{eqnarray*}
\!\!\!\!\!\!&-&\!\!\!\!\!\!\ 22144cs_0b^3\beta^3 -6016cs_0b^5\beta^3 -4512s_{0|0}nb^2\beta^2 -4200s_{0|0}nb^4\beta^2 -1024s_{0|0}n\beta^2b^6 -2560cb^4s_0\beta^3 \\
\!\!\!\!\!\!&+&\!\!\!\!\!\!\ 12352cb^2s_0b\beta^3 +10240cs_0\beta^3b^5 -4776n\beta^2b^2\theta +1024cs_0\beta^3b^7 -256cs_0\beta^3b^6 +64\beta^2b^8\sigma \\
\!\!\!\!\!\!&-&\!\!\!\!\!\!\ 7896s_ms^m_0n\beta^3b^2 -3552s_ms^m_0n\beta^3b^4 +6240c_0b^2\beta^3 +5424c_0b^4\beta^3 +1216c_0b^6\beta^3 -2880s_ms_0^mb^2\beta^2\\
\!\!\!\!\!\!&-&\!\!\!\!\!\!\ 2304s_ms_0^mb^4\beta^2 -9264s_0cn\beta^3b^2 -3840s_0cn\beta^3b^4 -9264s_{0|m}b^m\beta^3b^2 -3840s_{0|m}b^m\beta^3b^4
\\
\!\!\!\!\!\!&+&\!\!\!\!\!\!\ 4632s_{m|0}b^m\beta^3b^2 +1920s_{m|0}b^m\beta^3b^4 +3024s^i_ms^m_i\beta^4b^2 +1756(\textbf{Ric}-{\bf \overline{Ric}})\beta^2 +5600cs_0n\beta^3b^2
\\
\!\!\!\!\!\!&+&\!\!\!\!\!\!\ 2048cs_0n\beta^3b^4 +6016cs_0nb\beta^3 +11200cs_0nb^3\beta^3 +4096cs_0nb^5\beta^3 +7232(\textbf{Ric}-{\bf \overline{Ric}})\beta^2b^6 \\
\!\!\!\!\!\!&+&\!\!\!\!\!\!\ 1216(\textbf{Ric}-{\bf \overline{Ric}})\beta^2b^8 +8096(\textbf{Ric}-{\bf \overline{Ric}})\beta^2b^2 +12480(\textbf{Ric}-{\bf \overline{Ric}})\beta^2b^4 -480n\beta^2b^8\theta \\
\!\!\!\!\!\!&-&\!\!\!\!\!\!\ 64n\beta^2b^8\sigma -704n\beta^2b^6\sigma +292\beta^2\sigma +384cs_0\beta^3b^6 -6876n\beta^2b^4\theta -12176s_0^2\beta^2 +2080c^2\beta^4 \\
\!\!\!\!\!\!&+&\!\!\!\!\!\!\ 3556s_{0|0}\beta^2 +7368s_ms^m_0\beta^3 +7500cs_0\beta^3 +5000cs_0\beta^3 +6316s_{0|m}^m\beta^3 +2448s_0^2n\beta^2 \\
\!\!\!\!\!\!&-&\!\!\!\!\!\!\ 14208s_0^2\beta^2b^2 -8352s_0^2\beta^2b^4 -10240s_0^2b\beta^2 -5120s_0^2b^3\beta^2 +2176s_0^2b^5\beta^2 +3616c^2b^2\beta^4 \\
\!\!\!\!\!\!&+&\!\!\!\!\!\!\ 1216c^2b^4\beta^4 -8320c^2b\beta^4 -14464c^2b^3\beta^4 -4864c^2b^5\beta^4 -5168cs_0\beta^3 +9456s_{0|0}b^2\beta^2 +7464s_{0|0}b^4\beta^2 \\
\!\!\!\!\!\!&+&\!\!\!\!\!\!\ 1664s_{0|0}b^6\beta^2 -1388s_{0|0}n\beta^2 -3088cb^2s_0\beta^3 -3956s_ms^m_0n\beta^3 +12792s_ms^m_0\beta^3b^2 +5376s_ms^m_0\beta^3b^4 \\
\!\!\!\!\!\!&+&\!\!\!\!\!\!\ 13896cs_0\beta^3b^2 +5760cs_0\beta^3b^4 +9264cs_0\beta^3b^2 +3840cs_0\beta^3b^4 +1512s_ms^m\beta^4 +2024c_0\beta^3 \\
\!\!\!\!\!\!&-&\!\!\!\!\!\!\ 912s_ms_0^m\beta^2 -5000s_0cn\beta^3 -5000s_{0|m}b^m\beta^3 +2500s_{m|0}b^m\beta^3 +20000s_{0|m}^mb^2\beta^3 +18528s_{0|m}^mb^4\beta^3\\
\!\!\!\!\!\!&+&\!\!\!\!\!\!\ +5120s_{0|m}^m\beta^3b^6 +1944s^i_ms^m_i\beta^4 -384s_ms_0^m\beta^2b^6 +1608\beta^2b^4\sigma +6876\beta^2b^4\theta +256s_{0|m}^m\beta^3b^8 \\
\!\!\!\!\!\!&-&\!\!\!\!\!\!\ 1408s_0^2\beta^2b^6 +704\beta^2b^6\sigma -292n\beta^2\sigma -1083n\beta^2\theta +4776\beta^2b^2\theta +1220\beta^2b^2\sigma
\\
\!\!\!\!\!\!&+&\!\!\!\!\!\!\ 3552\beta^2b^6\theta +512s_0^2\beta^2b^7 +480\beta^2b^8\theta -256s_0cn\beta^3b^6 -256s_{0|m}b^m\beta^3b^6 +128s_{m|0}b^m\beta^3b^6 \\
\!\!\!\!\!\!&-&\!\!\!\!\!\!\ 1220n\beta^2b^2\sigma +384s_ms^m_0b^6\beta^3 +864s_ms^mb^2\beta^4 +864s^i_ms^m_i\beta^4b^4 +1280s_0^2nb^6\beta^2 -256s_ms^m_0n\beta^3b^6 \\
\!\!\!\!\!\!&-&\!\!\!\!\!\!\ 640cs_0b^6\beta^3 +256cs_0\beta^3b^6) +\beta^3(-24672cs_0nb^2\beta -11712cs_0n\beta b^4 +24336cs_0nb\beta +38016cs_0nb^3\beta \\
\!\!\!\!\!\!&+&\!\!\!\!\!\!\ 11520cs_0n\beta b^5 +19008c^2nb^2\beta^2 +5760c^2n\beta^2 b^4 +6912cs_0b^4\beta +21600cs_0b^2\beta -4980cs_0n\beta -768cs_0b^5\beta \\
\!\!\!\!\!\!&-&\!\!\!\!\!\!\ 20352cs_0b^3\beta -33744cs_0b\beta-5760c^2b^4\beta^2 -28224c^2b^2\beta^2 +12168c^2n\beta^2 -768c^2b^5\beta^2 -20352c^2b^3\beta^2 \\
\!\!\!\!\!\!&-&\!\!\!\!\!\!\ 33744c^2b\beta^2 -18252c_0nb^2\beta -14256c_0nb^4\beta -2880c_0nb^6\beta -5888c^2b^4\beta^2 -512c^2b^6\beta^2 +33152c^2b^3\beta^2 \\
\!\!\!\!\!\!&+&\!\!\!\!\!\!\ +23552c^2b^5\beta^2 +2048c^2b^7\beta^2 -24864c^2nb^2\beta^2 -8832c^2nb^4\beta^2 -512c^2nb^6\beta^2 -24864c_mb^mb^2\beta^2 \\
\!\!\!\!\!\!&-&\!\!\!\!\!\!\ 8832c_mb^mb^4\beta^2 -512c_mb^mb^6\beta^2 +24336c^2nb\beta^2 +38016c^2nb^3\beta^2 +11520c^2n\beta^2 b^5 +7554c_0\beta
\\
\!\!\!\!\!\!&+&\!\!\!\!\!\!\ 15688c^2\beta^2 +8640c_0\beta b^4 +1536c_0b^6\beta -18036cs_0\beta -26064c^2\beta^2 -6570c_0n\beta +14832c_0\beta b^2 \\
\!\!\!\!\!\!&-&\!\!\!\!\!\!\ 8288cb^2c\beta^2 +24864c^2b^2\beta^2 +8832c^2b^4\beta^2 +512c^2b^6\beta^2 -15688c^2n\beta^2 -15688c_mb^m\beta^2) \\
\!\!\!\!\!\!&+&\!\!\!\!\!\!\ 2\beta^5(5280c^2b^4+25296c^2b^2-8049c^2n -17112c^2nb^2-5808c^2nb^4 +11085c^2 )
\\
A'_8 \!\!\!\!&:=&\!\!\!\!\ -256s_ms_0^mb^8\beta +128s^i_ms^m_ib^8\beta^3 +256s_ms^m\beta^3b^6 +2176s^i_ms^m_i\beta^3b^6 +3264s_ms^m\beta^3b^4 +1440s_0^2n\beta b^2\\
\!\!\!\!\!\!&+&\!\!\!\!\!\!\ 3904s_0^2n\beta b^4 +4352s_0^2n\beta b^6 +512s_0^2nb\beta +1664s_0^2nb^3\beta +1280s_0^2nb^5\beta +256cs_0n\beta^2  -3168cs_0\beta^2 b^2\\
\!\!\!\!\!\!&-&\!\!\!\!\!\!\ 4992cs_0\beta^2 b^4-1280cs_0\beta^2 b^3 -1472cs_0b\beta^2 -2560cs_0\beta^2 b^6+1024cs_0\beta^2 b^5 -384s_{0|0}nb^2\beta \\
\!\!\!\!\!\!&-&\!\!\!\!\!\!\ 624s_{0|0}nb^4\beta -320s_{0|0}nb^6\beta +7424cb^3s_0\beta^2 +15872cs_0b^5\beta^2 -1664cb^6s_0\beta^2 +6656cb^7s_0\beta^2 \\
\!\!\!\!\!\!&-&\!\!\!\!\!\!\ 3744s_ms^m_0n\beta^2b^2 -4320s_ms^m_0n\beta^2b^4 -1344s_ms^m_0nb^6\beta^2 -3968cs_0b^4\beta^2 +672c_0b^2\beta^2 \\
\!\!\!\!\!\!&+&\!\!\!\!\!\!\ 1056c_0b^4\beta^2 +512c_0b^6\beta^2 -704s_ms_0^m\beta b^2 -1536s_ms_0^m\beta b^4 -1280s_ms_0^m\beta b^6 -5568s_0cnb^2\beta^2 \\
\!\!\!\!\!\!&-&\!\!\!\!\!\!\ 5952s_0cn\beta^2b^4 -1664s_0cn\beta^2b^6 -5568s_{0|m}b^mb^2\beta^2 -5952s_{0|m}b^m\beta^2b^4 -1664s_{0|m}b^m\beta^2b^6 \\
\!\!\!\!\!\!&+&\!\!\!\!\!\!\ 2784s_{m|0}b^mb^2\beta^2 +2976s_{m|0}b^m\beta^2b^4 +832s_{m|0}b^m\beta^2b^6 +6528s_ms^m\beta^3b^2 +5824s^i_ms^m_i\beta^3b^2
\\
\!\!\!\!\!\!&+&\!\!\!\!\!\!\ 6528s^i_ms^m_i\beta^3b^4 +1112s_{0|m}^m\beta^2 +80(\textbf{Ric}-{\bf \overline{Ric}})\beta +832cs_0n\beta^2 b^2+640cs_0n\beta^2 b^4 +512cs_0nb\beta^2 \\
\!\!\!\!\!\!&+&\!\!\!\!\!\!\ 1664cs_0nb^3\beta^2 +1280cs_0nb^5\beta^2 +1408(\textbf{Ric}-{\bf \overline{Ric}})\beta b^6 +512(\textbf{Ric}-{\bf \overline{Ric}})\beta b^8 +544(\textbf{Ric}-{\bf \overline{Ric}})\beta b^2 \\
\!\!\!\!\!\!&+&\!\!\!\!\!\!\ 1344(\textbf{Ric}-{\bf \overline{Ric}})\beta b^4 -2200s_0^2\beta +224c^2\beta^3 +1760s_ms^m_0\beta^2 +2256cs_0\beta^2 +1504cs_0\beta^2 +296s_0^2n\beta \\
\!\!\!\!\!\!&-&\!\!\!\!\!\!\ 4928s_0^2\beta b^2 -7680s_0^2\beta b^4 +2048s_0^2\beta b^7 +2816s_0^2\beta b^3 -576s_0^2b\beta -5120s_0^2\beta b^6 +6656s_0^2\beta b^5 +704c^2\beta^3b^2 \\
\!\!\!\!\!\!&+&\!\!\!\!\!\!\ 512c^2\beta^3b^4 -896c^2\beta^3b -2816c^2\beta^3b^3 -2048c^2\beta^3b^5 -800cs_0\beta^2 +768s_{0|0}\beta b^2 +816s_{0|0}\beta b^4 \end{eqnarray*}
\begin{eqnarray*}
\!\!\!\!\!\!&+&\!\!\!\!\!\!\ 256s_{0|0}\beta b^6 -76s_{0|0}n\beta -1856cb^2s_0\beta^2 -960s_ms^m_0n\beta^2 +4896s_ms^m_0\beta^2b^2 +4224s_ms^m_0\beta^2b^4 \\
\!\!\!\!\!\!&+&\!\!\!\!\!\!\ 1216s_ms^m_0b^6\beta^2 +8352cs_0b^2\beta^2 +8928cs_0\beta^2b^4 +2496cs_0\beta^2b^6 +5568cs_0b^2\beta^2 +5952cs_0\beta^2b^4 \\
\!\!\!\!\!\!&+&\!\!\!\!\!\!\ 1664cs_0\beta^2b^6 +57\beta\theta +18\beta\sigma +2912s_ms^m\beta^3 +136c_0\beta^2 -112s_ms_0^m\beta -1504s_0cn\beta^2 -1504s_{0|m}b^m\beta^2 \\
\!\!\!\!\!\!&+&\!\!\!\!\!\!\ 752s_{m|0}b^m\beta^2 +6016s_{0|m}^mb^2\beta^2 +11136s_{0|m}^mb^4\beta^2 +7936s_{0|m}^mb^6\beta^2 +1664s_{0|m}^m\beta^2b^8
\\
\!\!\!\!\!\!&+&\!\!\!\!\!\!\ 1544s^i_ms^m_i\beta^3 +1024s_0^2nb^8\beta -384n\beta b^2\theta -120n\beta b^2\sigma -936n\beta b^4\theta -288n\beta b^4\sigma -960n\beta b^6\theta \\
\!\!\!\!\!\!&-&\!\!\!\!\!\!\ 288n\beta b^6\sigma -336n\beta b^8\theta -96n\beta b^8\sigma -57n\beta\theta -18n\beta\sigma +384\beta b^2\theta +120\beta b^2\sigma +936\beta b^4\theta \\
\!\!\!\!\!\!&+&\!\!\!\!\!\!\ 288\beta b^4\sigma +960\beta b^6\theta +288\beta b^6\sigma +336\beta b^8\theta +96\beta b^8\sigma -1024s_0^2b^8\beta +2\beta(-2064cs_0nb^2\beta \\
\!\!\!\!\!\!&-&\!\!\!\!\!\!\ 2136cs_0n\beta b^4 +1104cs_0nb\beta +2976cs_0nb^3\beta +1824cs_0n\beta b^5 +1488c^2nb^2\beta^2 +912c^2n\beta^2 b^4 \\
\!\!\!\!\!\!&+&\!\!\!\!\!\!\ 672cs_0b^4\beta +2040cs_0b^2\beta -156cs_0n\beta +2976cs_0b^5\beta +1536cs_0b^3\beta -1200cs_0b\beta -912c^2b^4\beta^2\\
\!\!\!\!\!\!&-&\!\!\!\!\!\!\ 2712c^2b^2\beta^2 +552c^2n\beta^2 +2976c^2b^5\beta^2 +1536c^2b^3\beta^2 -1200c^2b\beta^2 -828c_0nb^2\beta -1116c_0nb^4\beta \\
\!\!\!\!\!\!&-&\!\!\!\!\!\!\ 456c_0nb^6\beta -3616c^2b^4\beta^2 -1216c^2b^6\beta^2+8320c^2b^3\beta^2 +14464c^2b^5\beta^2 +4864c^2b^7\beta^2 -6240c^2nb^2\beta^2 \\
\!\!\!\!\!\!&-&\!\!\!\!\!\!\ 5424c^2nb^4\beta^2 -1216c^2nb^6\beta^2 -6240c_mb^mb^2\beta^2 -5424c_mb^mb^4\beta^2 -1216c_mb^mb^6\beta^2 +1104c^2nb\beta^2\\
\!\!\!\!\!\!&+&\!\!\!\!\!\!\  2976c^2nb^3\beta^2 +1824c^2n\beta^2 b^5 +198c_0\beta +2024c^2\beta^2 -36c_0\beta b^4 -144c_0b^6\beta -2100cs_0\beta -1560c^2\beta^2\\
\!\!\!\!\!\!&-&\!\!\!\!\!\!\  192c_0n\beta +360c_0\beta b^2 -2080c^2b^2\beta^2 +6240c^2b^2\beta^2 +5424c^2b^4\beta^2  +1216c^2b^6\beta^2 -2024c^2n\beta^2 \\
\!\!\!\!\!\!&-&\!\!\!\!\!\!\  2024c_mb^m\beta^2) +\beta^3(6672c^2b^4 +19272c^2b^2-2256c^2n -9240c^2nb^2-6504c^2nb^4 +3504c^2)
\\
A'_{10} \!\!\!\!&:=&\!\!\!\!\ 8(1+2b^2)^2\big(8s^i_ms^m_i\beta b^4+4s_{0|m}^mb^4-2s_ms^m_0b^2-4s_{0|m}b^mb^2+6cs_0b^2 +16s_ms^m\beta b^2+4s_{0|m}^mb^2\\
\!\!\!\!\!\!&+&\!\!\!\!\!\!\ 4cs_0b^2+20s^i_ms^m_i\beta b^2-4s_0cnb^2-2s_ms^m_0nb^2+2s_{m|0}b^mb^2 +16cb^2s_0b+3cs_0 +2cs_0+s_{0|m}^m-4cb^2s_0\\
\!\!\!\!\!\!&-&\!\!\!\!\!\!\ 2s_{0|m}b^m+8s^i_ms^m_i\beta +26s_ms^m\beta -s_ms^m_0n +s_{m|0}b^m+2s_ms^m_0) +2\beta(336c^2b^4 +588c^2b^2-27c^2n\\
\!\!\!\!\!\!&-&\!\!\!\!\!\!\ 180c^2nb^2-216c^2nb^4 +39c^2) -64cs_0nb^2 -112cs_0nb^4 +64cs_0nb^3+64cs_0nb^5 +16cs_0nb +32c^2\beta nb^2 \\
\!\!\!\!\!\!&+&\!\!\!\!\!\!\ 64c^2\beta nb^3 +64c^2\beta nb^5 +16c^2\beta nb -704cb^2cb^2\beta -512c^2b^2b^4\beta +896c^2b^3\beta +2816c^2b^5\beta +320cs_0b^5\\
\!\!\!\!\!\!&+&\!\!\!\!\!\!\ 2048c^2b^7\beta -672c^2nb^2\beta -1056c^2nb^4\beta -512c^2nb^6\beta -672c_mb^mb^2\beta -1056c_mb^mb^4\beta -512c_mb^mb^6\beta \\
\!\!\!\!\!\!&-&\!\!\!\!\!\!\ 68cs_0-32c^2\beta -24c_0b^4 -2c_0n -32c_0b^6 +136c^\beta +212s_{0|0}\beta +64cs_0b^4 +112cs_0b^2 -4cs_0n  \\
\!\!\!\!\!\!&+&\!\!\!\!\!\!\ 128cs_0b^3 -16cs_0b -32c^2\beta b^4 -80c^2\beta b^2 +8c^2\beta n +320c^2\beta b^5 +128c^2\beta b^3 -16c^2\beta b -12c_0nb^2 \\
\!\!\!\!\!\!&-&\!\!\!\!\!\!\ 24c_0nb^4 -16c_0nb^6 -224c^2b^2\beta +672c^2b^2\beta +1056c^2b^4\beta +512c^2b^6\beta b^2 \\
\!\!\!\!\!\!&+&\!\!\!\!\!\!\ 2c_0 -136c^2n\beta -136c_mb^m\beta-2s_0cn +32c^2\beta nb^4\big).
\end{eqnarray*}

\bigskip

\noindent
Akbar Tayebi\\
Faculty  of Science, Department of Mathematics\\
University of Qom \\
Qom. Iran\\
Email:\ akbar.tayebi@gmail.com\\\\
Tayebeh Tabatabaeifar\\
Department of Mathematics and Computer Sciences\\
Amirkabir University\\
Tehran. Iran\\
Email:\  t.tabaee@gmail.com

\end{document}